\documentclass[12pt]{article}
\usepackage{amsmath}
\usepackage{amssymb}
\usepackage{amsthm}
\usepackage[totalwidth=17cm, totalheight=25cm]{geometry}
\usepackage{graphicx}
\usepackage{mathabx}
\usepackage{tabularx}
	\newcolumntype{C}[1]{>{\centering\arraybackslash}m{#1}} 
	\newcolumntype{R}[1]{>{\raggedleft\arraybackslash}m{#1}} 
\newtheoremstyle{boldplain}% name
{9pt}%      Space above
{9pt}%      Space below
{\itshape}%         Body font
{}%         Indent amount (empty = no indent, \parindent = para indent)
{\bfseries}% Thm head font
{.}%        Punctuation after thm head
{.5em}%     Space after thm head: " " = normal interword space;
{\thmname{#1}\thmnumber{ #2}\thmnote{ (#3)}}%

\newtheoremstyle{bolddefinition}% name
{9pt}%      Space above
{9pt}%      Space below
{}%         Body font
{}%         Indent amount (empty = no indent, \parindent = para indent)
{\bfseries}% Thm head font
{.}%        Punctuation after thm head
{.5em}%     Space after thm head: " " = normal interword space;
{\thmname{#1}\thmnumber{ #2}\thmnote{ (#3)}}%

\theoremstyle{boldplain}
\newtheorem{add}[equation]{Addendum}

\newtheorem{cor}[equation]{Corollary}

\newtheorem{lem}[equation]{Lemma}
\newtheorem{lemma}[equation]{Lemma}
\newtheorem{prop}[equation]{Proposition}

\newtheorem{sublem}[equation]{Sublemma}
\newtheorem{thm}[equation]{Theorem}

\newtheorem{eqthm}[equation]{Equivalence Theorem}

\theoremstyle{bolddefinition}
\newtheorem{dfn}[equation]{Definition}
\newtheorem{definition}[equation]{Definition}

\newtheorem{rem}[equation]{Remark}

\newfont{\bigbf}{cmbx10 scaled\magstep1}
\normalbaselines
\numberwithin{equation}{section}

\def\no{\noindent}

\def\R{{\mathbb R}}

\def\N{{\mathbb N}}

\def\Z{{\mathbb Z}}

\def\al{\alpha}
\def\ga{\gamma}
\def\Ga{\Gamma}
\def\de{\delta}
\def\De{\Delta}
\def\eps{\epsilon}
\def\la{\lambda}
\def\La{\Lambda}
\def\si{\sigma}
\def\Si{\Sigma}
\def\ups{\upsilon}

\def\Om{\Omega}

\def\3{\ss}

\def\acct{\operatorname{acc}_{\tau_{mod}}}
\def\acts{\curvearrowright}
\def\amod{a_{mod}}

\def\B{\operatorname{B}}

\def\D{\partial}
\def\DF{\partial_{F\ddot u}}
\def\Ds{\partial_{\si_{mod}}}
\def\Dt{\partial_{\tau_{mod}}}
\def\Dtp{\partial_{\tau'_{mod}}}

\def\diam{\mathop{\hbox{diam}}}
\def\diamo{\diamondsuit}

\def\diamot{\diamondsuit_{\tau_{mod}}}
\def\diamoTh{\diamondsuit_{\Theta}}

\def\Fix{\operatorname{Fix}}
\def\Flag{\operatorname{Flag}}
\def\Flagn{\operatorname{Flag_{\nu_{mod}}}}
\def\Flags{\operatorname{Flag_{\si_{mod}}}}

\def\Flagt{\operatorname{Flag_{\tau_{mod}}}}
\def\Flagit{\Flag_{\iota\tau_{mod}}}
\def\Flagmt{\Flag_{-\tau_{mod}}}
\def\Flagpt{\Flag_{+\tau_{mod}}}
\def\Flagpmt{\Flag_{\pm\tau_{mod}}}

\def\Fmod{F_{mod}}

\def\geo{\partial_{\infty}}
\def\geot{\partial_{\infty}^{\tau_{mod}-reg}}

\def\Hb{\operatorname{Hb}}
\def\Hc{\operatorname{Hc}}
\def\Hs{\operatorname{Hs}}

\def\Homeo{\operatorname{Homeo}}
\def\id{\mathop{\hbox{id}}}

\def\inte{\operatorname{int}}
\def\interior{\operatorname{int}}
\def\Isom{\mathop{\hbox{Isom}}}
\def\Jac{\operatorname{Jac}}
\def\LaGa{\Lambda(\Ga)}
\def\Las{\Lambda_{\si_{mod}}}
\def\LasGa{\Lambda_{\si_{mod}}(\Ga)}
\def\Lat{\Lambda_{\tau_{mod}}}
\def\LatGa{\Lambda_{\tau_{mod}}(\Ga)}
\def\LatGacon{\Lambda^{con}_{\tau_{mod}}(\Ga)}

\def\lra{\longrightarrow}

\def\numod{\nu_{mod}}

\def\oa{\overrightarrow}
\def\ol{\overline}

\def\pihalf{\frac{\pi}{2}}

\def\2pithird{\frac{2\pi}{3}}

\def\rank{\mathop{\hbox{rank}}}

\def\Ra{\Rightarrow}

\def\simod{\si_{mod}}

\def\st{\operatorname{st}}
\def\stTh{\operatorname{st}_{\Theta}}

\def\ost{\operatorname{ost}}

\def\Stab{\operatorname{Stab}}

\def\tangle{\angle_{Tits}}
\def\taumod{\tau_{mod}}

\def\tits{\partial_{Tits}}

\def\Vert{\operatorname{Vert}}

\def\Wn{W_{\nu_{mod}}}
\def\Wt{W_{\tau_{mod}}}

\def\<{\langle}
\def\>{\rangle}

\def\BI{\begin{itemize}}
\def\EI{\end{itemize}}

\long\def\comment#1\endcomment{}

\title{Anosov subgroups: \\
Dynamical and geometric characterizations}
\author{Michael Kapovich, Bernhard Leeb, Joan Porti}
\date{March 3, 2017}

\begin{document}

\maketitle

\begin{abstract}
We study infinite covolume discrete subgroups
of higher rank semisimple Lie groups,
motivated by understanding basic properties of Anosov subgroups 
from various viewpoints 
(geometric, coarse geometric and dynamical).
The class of Anosov subgroups 
constitutes a natural generalization of convex cocompact subgroups of rank one Lie groups
to higher rank.
Our main goal is to give several new 
equivalent characterizations for this important class of discrete subgroups.
Our characterizations capture ``rank one behavior'' of Anosov subgroups 
and are direct generalizations 
of rank one equivalents to convex cocompactness.
Along the way,
we considerably simplify the original definition,
avoiding the geodesic flow. 
We also show that the Anosov condition can be relaxed further 
by requiring only non-uniform unbounded expansion 
along the (quasi)geodesics in the group. 
\end{abstract}

\tableofcontents

\section{Introduction}

This paper is devoted to studying 
basic properties of Anosov subgroups 
of semisimple Lie groups from various viewpoints 
(geometric, coarse geometric and dynamical).
The class of Anosov subgroups,
introduced by Labourie \cite{Labourie} 
and further extended by Guichard and Wienhard \cite{GW}, 
constitutes a natural generalization of convex cocompact subgroups of rank one Lie groups
to higher rank.
Our main goal here is to give several new 
equivalent characterizations for this important class of discrete subgroups,
including a considerable simplification of their original definition. 
For convex cocompact subgroups as well as for word hyperbolic groups, 
it is very fruitful to have different viewpoints and 
alternative definitions,
as they were developed by many authors 
starting with Ahlfors' work on geometric finiteness in the 60s,
and later by Beardon, Maskit, Marden, Thurston, Sullivan, Bowditch and others.
Besides a deeper understanding, 
it enables one to switch perspectives in a nontrivial way, 
adapted to the situation at hand. 
A main purpose of this paper is to demonstrate 
that much of this theory extends to Anosov subgroups,
and we hope that the concepts and results presented here 
will be useful for their further study. 
In our related work,
they lay the basis for the results on the Higher Rank Morse Lemma \cite{mlem},
compactifications of locally symmetric spaces for Anosov subgroups \cite{bordif},
the local-to-global principle and the construction of Morse-Schottky subgroups \cite{morse}.
We refer to the surveys \cite{anosov,manicures}
for more details on these developments. 

In rank one, 
among Kleinian groups and, more generally, among discrete subgroups of rank one Lie groups,
one distinguishes 
{\em geometrically finite} subgroups.
They form a large and flexible class of discrete subgroups 
which are strongly tied to the negatively curved symmetric spaces they act on.
Therefore they have especially good geometric, topological and dynamical properties 
and one can prove many interesting results about them.
The simplest are geometrically finite subgroups {\em without parabolics}, 
which are lie at the root of this paper.
They can be characterized in many (not obviously) equivalent ways: 
As convex cocompact subgroups,
as undistorted subgroups,
as subgroups with conical limit set, 
as subgroups which are expanding at their limit set,
and as intrinsically word hyperbolic subgroups 
with Gromov boundary equivariantly homeomorphic to their limit set,
to name some.

In higher rank,
a satisfying and sufficiently broad definition of geometric finiteness, with or without parabolics, 
remains yet to be found. 
Convex cocompactness turns out to be much too restrictive a condition:
it was shown by Kleiner and the second author \cite{convcoco}
that in higher rank only few subgroups are convex cocompact.
Undistortion by itself,
on the other hand,
is way too weak:
undistorted subgroups can even fail to be finitely presented.
Thus, one is forced to look for suitable replacements 
of these notions 
in higher rank. 
It turns out that some of the other equivalent characterizations 
of convex cocompactness in rank one
do admit useful modifications in higher rank,
which lead to the class of Anosov subgroups.
The Anosov condition is not too rigid and, 
at the same time, it imposes enough restrictions on the subgroups 
making it possible to analyze their geometric and dynamical properties. 
One way to think of Anosov subgroups is 
as geometrically finite subgroups without parabolics
which exhibit some {\em rank one behavior}.
Indeed, 
they are intrinsically word hyperbolic
and we will see that also extrinsically they display hyperbolic behavior 
in a variety of ways.

In this paper, we primarily consider 
four notions generalizing convex cocompactness to higher rank,
all equivalent to the Anosov condition,
see the Equivalence Theorem~\ref{thm:eqvi} below:

(i) asymptotic embeddedness

(ii) expansivity

(iii) conicality

(iv) Morse property

Whereas the conditions Anosov, (i) and (ii) are dynamical,
(iii) is a condition on the asymptotic geometry of the subgroup,
and (iv) is coarse geometric. 

\medskip
We now describe in more detail some of our concepts and results. 

Let $X=G/K$ be a symmetric space of noncompact type
and, for simplicity, let the semisimple Lie group $G$ be the connected component 
of its isometry group. 
Our approach to studying Anosov subgroups $\Ga<G$ begins with the observation 
that they satisfy a strong form of discreteness 
which we call {\em regularity}
and which is primarily responsible for their extrinsic ``rank one behavior'' 
alluded to above. 
Discreteness of a subgroup $\Ga<G$ means that 
for sequences $(\ga_n)$ of distinct elements 
the distance $d(x,\ga_nx)$ in $X$ diverges to infinity. 
For higher rank symmetric spaces 
there is a natural vector-valued refinement $d_{\De}$ 
of the Riemannian distance $d$,
which takes values in the euclidean Weyl chamber $\De$ of $X$. 
The regularity assumption on $\Ga$,
in its strongest form of {\em $\simod$-regularity}, 
means that $d_{\De}(x,\ga_nx)$ 
diverges away from the boundary of $\De$. 
We will work more generally with relaxations of this condition,
called $\taumod$-regularity,
associated with a face $\taumod$ of the model spherical Weyl chamber $\simod$, 
where one only requires divergence of $d_{\De}(x,\ga_nx)$ 
away from some of the faces of $\De$, depending on $\taumod$. 
To be precise,
think of $\simod$ as the visual boundary of the euclidean Weyl chamber,
$\simod\cong\geo\De$.
Given a face $\taumod\subseteq\simod$,
we define {\em $\taumod$-regularity} by requiring 
that $d_{\De}(x,\ga_nx)$ diverges away from the faces of $\De$
whose visual boundaries do not contain $\taumod$.
We will also need the stronger notion of {\em uniform $\taumod$-regularity}
where one requires the divergence to be linear in terms of $d(x,\ga_nx)$.
Most of the discussion in this paper 
will take place within the framework of $\taumod$-regular subgroups.

Classically,
the asymptotic behavior of discrete subgroups $\Ga<G$ 
is captured by their visual limit set $\La(\Ga)$ 
which is the accumulation set of their orbits $\Ga x\subset X$ 
in the visual boundary $\geo X$.
In our context of $\taumod$-regular subgroups,
the visual limit set is replaced by the {\em $\taumod$-limit set} $\LatGa$
contained in the partial flag manifold $\Flagt=G/P_{\taumod}$ 
and defined as the accumulation set of $\Ga$-orbits 
in the bordification $X\sqcup\Flagt$ of $X$,
equipped with the topology of {\em flag convergence} 
(see section~\ref{sec:conv}).
Here, $P_{\taumod}$ is a parabolic subgroup 
in the conjugacy class corresponding to $\taumod$.
The notion of $\taumod$-limit set extends to arbitrary discrete subgroups.

We call a $\taumod$-regular subgroup $\Ga<G$ 
{\em nonelementary}
if $|\LatGa|\geq3$,
and {\em antipodal}
if it satisfies the visibility condition
that any two distinct limit simplices in $\LatGa$ are antipodal.
The latter means that they can be connected by a geodesic in $X$ 
in the sense that the geodesic is asymptotic to interior points 
of the simplices. 
It is worth noting that
the action of a $\taumod$-regular antipodal subgroup on its $\taumod$-limit set
enjoys the classical convergence property,
which is a typical rank one phenomenon. 

Regularity, which is a condition on the asymptotic geometry 
of orbits in the symmetric space, 
can be converted into an equivalent dynamical condition 
about a certain contraction behavior of the subgroup on suitable flag manifolds 
(see Definition~\ref{def:contracting_sequence}), 
allowing one to switch between geometry and dynamics.
The contraction behavior here is a higher rank version 
of the classical convergence (dynamics) property in the theory of Kleinian groups. 
This yields an equivalent characterization of $\taumod$-regular subgroups
as {\em $\taumod$-convergence subgroups}
(see Definition~\ref{def:conv}). 
Also the limit sets, respectively, limit simplices 
can be defined purely dynamically
as the possible limits of contracting sequences in $\Ga$,
i.e.\ of sequences converging to constants on suitable open and dense subsets 
of the flag manifolds, see Definition~\ref{def:flagl}.

Much of the material in section~\ref{sec:regcontr} 
can be found in some form already in the work of Benoist, 
see \cite[\S 3]{Benoist},
in the setting of Zariski dense subgroups of reductive algebraic groups over local fields,
notably the notions of {\em regularity} and {\em contraction}, 
their essential equivalence, and the notion of {\em limit set}.
For the sake of completeness we give independent proofs in our setting
of discrete subgroups of semisimple Lie groups.
Also our methods are rather different.
We give here a geometric treatment 
and present the material 
in a form suitable to serve as a basis for the further development 
of our theory of discrete isometry groups acting on Riemannian symmetric spaces and euclidean buildings of higher rank,
such as in our papers \cite{coco15,morse,mlem,bordif}.

\medskip
We now (mostly) restrict to the class of $\taumod$-regular, equivalently,
$\taumod$-convergence subgroups
and introduce various geometric and dynamical conditions 
in the spirit of 
geometric finiteness.
We begin with three dynamical ones:

1.\ 
We say that a 
subgroup $\Ga<G$ 
is {\em $\taumod$-asymptotically embedded}
if it is an antipodal $\taumod$-convergence subgroup,
$\Ga$ is word hyperbolic 
and there exists a $\Ga$-equivariant homeomorphism
\begin{equation*}
\label{eq:mapalphatauintro}
\alpha: \geo \Gamma \stackrel{\cong}{\lra}
\LatGa\subset \Flagt
\end{equation*}
from its Gromov boundary onto its $\taumod$-limit set. 

This condition can be understood as a continuity at infinity property
for the orbit maps $o_x:\Ga\to\Ga x\subset X$:
By extending an orbit map $o_x$ to infinity by the boundary map $\al$,
one obtains a continuous map 
$$o_x\sqcup\al:\Ga\sqcup\geo\Ga\lra X\sqcup\Flagt$$ 
from the Gromov compactification of $\Ga$
(see Proposition~\ref{prop:asyembprp}). 

2.\ 
Our next condition
is inspired by Sullivan's notion of {\em expanding actions} \cite{Sullivan}.
Following Sullivan, 
we call a subgroup $\Ga< G$ {\em expanding at infinity}
if its action on the appropriate partial flag manifold 
is expanding at the limit set.
More precisely:

We call a $\taumod$-convergence subgroup $\Ga< G$
{\em $\taumod$-expanding at the limit set}
if for every limit flag in $\LatGa$
there exists a neighborhood $U$ in $\Flagt$ and an element $\ga\in\Ga$
which is uniformly expanding on $U$,
i.e.\ for some constant $c>1$ and all $\tau_1,\tau_2\in U$ 
it holds that
\begin{equation*}
d(\ga\tau_1,\ga\tau_2)\geq c\cdot d(\tau_1,\tau_2) .
\end{equation*}
Here, and in what follows the distance $d$ is induced by a fixed 
Riemannian background metric on the flag manifold.

Now we can formulate our second condition:

We say that a subgroup $\Ga <G$ is {\em $\taumod$-CEA} 
(Convergence Expanding Antipodal)
if it is an antipodal $\taumod$-convergence subgroup
which is expanding at the limit set.

We note that the CEA condition 
does not a priori assume word hyperbolicity,
not even finite generation. 

3.\ 
The next condition is motivated by the original definition of Anosov subgroups. 
It is a hybrid of the previous two definitions,
where we weaken 
asymptotic embeddedness (to boundary embeddedness)
and strengthen expansivity. 
We drop the regularity/convergence assumption
and, accordingly, make no use of the limit set in our definition.
Compared to asymptotic embeddedness,
we keep the word hyperbolicity of the subgroup but, 
instead of identifying its Gromov boundary with the limit set 
as in asymptotic embeddedness,
we only require a boundary map embedding the Gromov boundary into the flag manifold.
Compared to CEA,
we require a stronger form of expansivity, now at the image of the boundary map.

We call a subgroup $\Ga<G$ 
{\em $\taumod$-boundary embedded} 
if $\Ga$ is word hyperbolic and 
there exists a $\Ga$-equivariant continuous embedding 
$$\beta:\geo\Ga\lra\Flagt$$
sending distinct visual boundary points 
to antipodal simplices. 
If $\Ga$ is virtually cyclic, we require in addition 
that it is discrete in $G$.
(Otherwise, discreteness is a consequence.)
We will refer to $\beta$ as a {\em boundary embedding}.
In general, boundary embeddings are not unique. 

The {\em infinitesimal expansion factor} 
of an element $g\in G$ at a simplex $\tau\in\Flagt$ is 
\begin{equation*}
\eps(g, \tau)= \min_u |dg(u)|
\end{equation*}
where the minimum is taken over all unit tangent vectors $u\in T_\tau \Flagt$,
again using the Riemannian background metric.

Now we can formulate our version of the Anosov condition:

We say that a subgroup $\Ga<G$ is {\em $\taumod$-Anosov} if 
it is $\taumod$-boundary embedded with boundary embedding $\beta$ 
and satisfies the following expansivity condition:
For every ideal point $\zeta\in \geo \Ga$
and every normalized 
(by $r(0)=e\in \Ga$)
discrete geodesic ray $r: \N\to \Ga$ asymptotic to $\zeta$,
the action $\Ga\acts\Flagt$ satisfies
\begin{equation*}
\eps(r(n)^{-1}, \beta(\zeta))\ge A e^{Cn}
\end{equation*}
for $n\geq 0$ with constants $A, C>0$  independent of $r$. 
(Here, we fix a word metric on $\Ga$.)

The uniformity of expansion in this definition can be significantly weakened: 

We say that a subgroup $\Ga<G$ is {\em non-uniformly $\taumod$-Anosov} if
it is $\taumod$-boundary embedded with boundary embedding $\beta$ and, 
for every ideal point $\zeta\in \geo \Ga$
and every discrete geodesic ray $r: \N\to \Ga$ asymptotic to $\zeta$,
the action $\Ga\acts\Flagt$ satisfies 
\begin{equation*}
\sup_{n\in\N}\, \eps(r(n)^{-1}, \beta(\zeta))=+\infty.
\end{equation*}

The original definition of Anosov subgroups in \cite{Labourie,GW} is rather involved.
It is based on geodesic flows for word hyperbolic groups 
and formulated in terms of expansion/contraction properties 
for lifted flows on associated bundles over the geodesic flow spaces
(see section~\ref{sec:onos}).
Our definition requires only an expansion property for the group action on a suitable flag manifold  
and avoids using the geodesic flow,
whose construction is highly technical for word hyperbolic groups which do not arise as the fundamental group of a closed negatively curved Riemannian manifold.
The geodesic flow is replaced by a simpler coarse geometric object,
the space of quasigeodesics.

\medskip

Now we come to the geometric notions.

4.\ 
The first geometric condition concerns the {\em orbit asymptotics}.
The notion of conicality of limit simplices,
due to Albuquerque \cite[Def.\ 5.2]{Albuquerque}, 
generalizes a well-known condition from the theory of Kleinian groups: 
In the case $\taumod=\simod$,
a limit chamber $\si\in\LasGa$
of a $\simod$-regular subgroup $\Ga<G$
is called {\em conical}
if there exists a sequence $\ga_n\to\infty$ in $\Ga$
such that for a(ny) point $x\in X$ the sequence of orbit points $\ga_nx$
is contained in a tubular neighborhood of the euclidean Weyl chamber $V(x,\si)$
with tip $x$ and asymptotic to $\si$.
For general $\taumod$
and limit simplices $\tau\in\LatGa$ of $\taumod$-regular subgroups $\Ga<G$,
one replaces the euclidean Weyl chamber 
with the Weyl cone $V(x,\st(\tau))$ over the star of $\tau$,
that is, by the union of the euclidean Weyl chambers $V(x,\si)$ 
for all spherical Weyl chambers $\si\supset\tau$.
A $\taumod$-regular subgroup $\Ga<G$ is called {\em conical}
if all limit simplices are conical.
Here is our forth condition:

We say that a subgroup $\Ga<G$
is {\em $\taumod$-RCA} if it is $\taumod$-regular, conical and antipodal.

For nonelementary $\taumod$-regular antipodal subgroups, 
this extrinsic notion of conicality 
is equivalent to an intrinsic one
defined in terms of the dynamics on the $\taumod$-limit set
(Proposition~\ref{prop:concon}),
which enables one to relate it to the dynamical notions above. 

5.\ 
The last set of definitions concerns the {\em coarse extrinsic geometry}.
We recall that a finitely generated subgroup $\Ga<G$ is {\em undistorted} 
if the orbit maps $\Ga\to X$ are quasiisometric embeddings. 
They then send discrete geodesics in $\Ga$ (with respect to a fixed word metric)
to uniform quasigeodesics in $X$.
Undistortion by itself is too weak a restriction, 
compared with the other notions defined previously.
We will strengthen it in two ways.
The first is by adding uniform regularity:

We say that a subgroup $\Ga<G$ is {\em $\taumod$-URU} 
if it is uniformly $\taumod$-regular and undistorted. 

According to the classical {\em Morse Lemma} in negative curvature, 
quasigeodesic segments in rank one symmetric spaces are uniformly Hausdorff close to geodesic segments
with the same endpoints.
This is no longer true in higher rank because it fails already in euclidean plane. 
Another way of strengthening undistortion
is therefore by imposing a ``Morse'' type property on the 
quasigeodesics arising as 
orbit map images of the discrete geodesics in $\Ga$. 

As in the case of conicality above,
where one replaces rays with Weyl cones when passing from rank one to higher rank,
it is natural to replace geodesic segments with ``diamonds'' in a higher rank version of the Morse property.
(This is suggested, for instance, by the geometry of free Anosov subgroups,
see our examples of Morse-Schottky subgroups \cite{morse,manicures}.)
We define diamonds as follows:
If $\taumod=\simod$ and $xy$ is a $\simod$-regular segment,
then the {\em $\simod$-dia\-mond} with tips $x,y$
is the intersection 
$$ \diamo(x,y) = V(x,\si)\cap V(y,\hat\si)$$
of the euclidean Weyl chambers with tips at $x$ and $y$ containing $xy$.
In the case of general $\taumod$, the euclidean Weyl chambers are replaced with $\taumod$-Weyl cones
(see section~\ref{sec:diamo}).

We say that a subgroup $\Ga<G$ is {\em $\taumod$-Morse}
if it is $\taumod$-regular, $\Ga$ is word hyperbolic 
and an(y) orbit map $o_x:\Ga\to\Ga x\subset X$ satisfies the following Morse condition:
The images $o_x\circ s$ of discrete geodesic segments $s:[n_-,n_+]\cap\Z\to\Ga$
are contained in uniform tubular neighborhoods of $\taumod$-diamonds with tips uniformly close 
to the endpoints of $o_x\circ s$
(see Definition~\ref{def:mrs}). 

The definition does not a priori assume undistortion,
but we show in this paper that Morse implies URU.
That, conversely, URU implies Morse may seem unexpected at first 
but follows from our Higher Rank Morse Lemma for regular quasigeodesics \cite{mlem}.

\medskip
We now arrive at our main result on the equivalence of various conditions introduced above.
We state it for nonelementary subgroups because we use this assumption in some of our proofs.

\begin{eqthm}
\label{thm:eqvi}
The following properties for subgroups $\Ga<G$ 
are equivalent in the nonelementary case:

(i)  $\taumod$-asymptotically embedded

(ii) $\taumod$-CEA.

(iii) $\taumod$-Anosov

(iv) non-uniformly $\taumod$-Anosov

(v) $\taumod$-RCA

(vi) $\taumod$-Morse

These properties imply $\taumod$-URU.

Moreover,
the boundary maps for properties (i), (iii) and (iv) coincide.
\end{eqthm}
Here, ``nonelementary'' 
means
$|\geo\Ga|\geq3$
in the Anosov conditions (iii) and (iv),
which assume word hyperbolicity but no $\taumod$-regularity,
and means $|\LatGa|\geq3$ in all other cases.

\begin{rem}
(i)
We prove in \cite{mlem} that, conversely,
$\taumod$-URU implies $\taumod$-Morse (without assuming nonelementary).

(ii) 
All implications between properties (i)-(vi) hold without assuming nonelementary,
with the exception of (ii)$\Rightarrow$(v)$\Rightarrow$(i).
In particular,
the properties (i),(iii),(iv),(vi) and $\taumod$-URU are equivalent in general.

(iii)
The implication Anosov$\Rightarrow$URU had been known before \cite{GW}.

(iv)
Some of the implications in the theorem can be regarded 
as a description of geometric and dynamical properties of Anosov subgroups. 
Different characterizations of Anosov subgroups 
are useful in different contexts.
For example:
Expansivity (ii) is used in \cite{coco13,coco15}
to establish the cocompactness of $\Ga$-actions 
on suitable domains of discontinuity in flag manifolds. 
Asymptotic embeddedness is used in \cite{bordif}
to construct Finsler compactifications 
of locally symmetric spaces for Anosov subgroups.
The Morse property is used in \cite{morse} to prove 
a {\em local-to-global principle} for Anosov subgroups.
The latter in turn leads to new proofs 
of openness and structural stability of Anosov representations,
to a construction of free Anosov subgroups (Morse-Schottky subgroups),
and to the semidecidability of Anosovness,
see \cite{morse}.

(v) 
In our paper \cite{bordif} we establish two more characterizations 
of Anosov subgroups among uniformly regular subgroups,
namely as {\em coarse retracts} and by {\em $S$-cocompactness}.
The former property is a strengthening of undistortion.
The latter means the existence of a certain kind of compactification 
of the corresponding locally symmetric space. 

(vi) 
Other characterizations of Anosov subgroups can be found in \cite{GGKW}.
\end{rem}

\begin{rem}
Boundary embeddedness appears to be a considerable weakening 
of asymptotic embeddedness, even in the regular case. 
Nevertheless two results in this paper establish a close relation 
between the two concepts:

(i) 
For $\simod$-regular subgroups,
boundary embeddedness, conversely, implies asymptotic embeddedness,
while the boundary embedding may have to be modified 
(see Theorem~\ref{thm:bemmdsi}). 

(ii) 
For general $\taumod$-regular subgroups,
there is the following {\em dichotomy} for boundary embeddings
(see Theorem~\ref{thm:bemblim})
which is useful for verifying asymptotic embeddedness:

Either the image of the boundary embedding equals the $\taumod$-limit set
and the subgroup is asymptotically embedded.
Or the image is disjoint from the limit set,
and the limit set is not Zariski dense.
The latter cannot happen for Zariski dense subgroups. 
\end{rem}

While the main results in this paper concern discrete subgroups of Lie groups,
in section~\ref{sec:mqg},
motivated by the Morse property,
we discuss {\em Morse quasigeodesics} and {\em Finsler geodesics}.
We characterize Morse subgroups as word hyperbolic subgroups 
whose intrinsic geodesics are extrinsically uniform Morse quasigeodesics. 
Furthermore, 
we characterize Morse quasigeodesics as bounded perturbations of Finsler geodesics.
Lastly, 
we analyze the $\De$-distance along Finsler geodesics and Morse quasigeodesics.
We show that, via the $\De$-distance function, they project to Finsler geodesics and Morse quasigeodesics
in $\De$.

Most of the results in this paper 
were already contained in chapters 1-6 of the preprint \cite{morse},
however the presentation in this paper is more efficient. 
The further material on the Morse property in \cite[\S 7]{morse}
will appear elsewhere. 

{\bf Acknowledgements.} 
The first author was partly supported by the NSF grants  DMS-12-05312 and DMS-16-04241. 
He is also thankful to KIAS (the Korea Institute for Advanced Study) for its hospitality. 
The third author was partially supported by  grant FEDER/Mineco MTM2015-66165-P.

\section{Geometry of symmetric spaces}
\label{sec:prelim}

In this section, we collect some material from the geometry of symmetric spaces and buildings.
We explain the notions which are most important for the purposes of this paper,
establish notation and give proofs for some of the less standard facts.
No attempt of a complete review is made. For more detailed discussions, 
we refer the reader to \cite{Eberlein}, \cite{BGS}, \cite{qirigid} and \cite{habil}.

We give a brief description of where various parts of this section are used in the paper:

Sections~\ref{sec:spherical}-\ref{sec:stars} are used essentially everywhere.

While the vector valued distance function $d_{\De}$ is used in many places in the paper, 
the rest of the material in sections~\ref{sec:VVD} and~\ref{sec:refined} 
is used primarily in section~\ref{sec:sepnest} on the separation of nested Weyl cones 
and in section~\ref{sec:mqg} where we analyze projections of Morse quasigeodesics to the euclidean model Weyl chamber $\De$. 

The material of section~\ref{sec:ascones-shadows} 
dealing with shadows at infinity is used in section~\ref{sec:reg-con} 
when we prove the equivalence of {\em regularity} and {\em contractivity} for sequences of isometries of $X$. 
The main result of section~\ref{sec:sepnest} on the separation of nested Weyl cones 
is used in section~\ref{sec:asyemb} to prove that Morse subgroups are URU (Theorem~\ref{thm:asymbur}). 

The main results of sections~\ref{sec:horocycles} and~\ref{sec:contraction}  
are Theorem~\ref{thm:expand} and Proposition~\ref{prop:expand} 
establishing estimates for the contraction and expansion of isometries of $X$ acting on flag manifolds. 
(The other results are only used only in  sections~\ref{sec:horocycles} and~\ref{sec:contraction}). 
Theorem~\ref{thm:expand} and Proposition~\ref{prop:expand}  
are used in sections~\ref{sec:expa} and~\ref{sec:anosov} 
while discussing discrete subgroups satisfying expansion properties (CEA and Anosov).

The material of section \ref{sec:fins} is used only in section \ref{sec:mqg}
where it is proven that Morse quasigeodesics are uniformly closed to Finsler geodesics 
and that $\De$-distance projections of Finsler geodesics are again Finsler geodesics.

\subsection{General metric space notation} \label{sec:general}

We will use the notation $B(p,r)$ 
for the open $r$-ball with center $p$ in a metric space, and $\ol B(p,r)$ for the closed $r$-ball.

A {\em geodesic} in a metric space $(Z,d)$ is an isometric embedding $I\to Z$ 
from a (possibly infinite) interval $I\subset\R$. 
In the context of finitely generated groups equipped with word metrics,
we will also work with {\em discrete geodesics};
these are isometric embeddings from intervals $I\cap\Z$ in $\Z$.
The notion of {\em discrete quasigeodesic} will be used similarly.

\subsection{Spherical buildings}
\label{sec:spherical}

Spherical buildings occur in this paper 
as the visual boundaries of symmetric spaces of noncompact type,
equipped with their structures of thick spherical Tits buildings.

\subsubsection{Spherical geometry}
\label{sec:sphgeom}

Let $S$ be a unit sphere in a euclidean space,
and let $\si\subset S$ be a spherical simplex with dihedral angles $\leq\pihalf$.
Then $\diam(\si)\leq\pihalf$.

For a face $\tau\subseteq\si$,
we define the {\em $\tau$-boundary}  $\D_{\tau}\si$ as the union of faces of $\si$ which do not contain $\tau$,
and the {\em $\tau$-interior} $\inte_{\tau}(\si)$ as the union of open faces of $\si$ whose closure contains $\tau$.
We obtain the decomposition 
$$\si=\inte_{\tau}(\si)\sqcup\D_{\tau}\si.$$
If $\tau'\subset\tau$, 
then $\D_{\tau'}\si\subset\D_{\tau}\si$ and $\inte_{\tau'}(\si)\supset\inte_{\tau}(\si)$.
Note that $\D_{\si}\si=\D\si$ and $\inte_{\tau}(\si)=\inte(\si)$.

We need the following fact about projections of spherical simplices to their faces:
\begin{lem}
\label{lem:sphprjfc}
The nearest point projection $\inte_{\tau}(\si)\to\inte(\tau)$ is well-defined.
\end{lem}

In other words,
for every point $x\in\inte_{\tau}(\si)$ there exists a point $p\in\inte(\tau)$ 
such that $px\perp\tau$.
In view of $\diam(\si)\leq\pihalf$, this point is necessarily unique.

\proof
We argue by induction on the dimension of $\si$.

Let $x\in\inte_{\tau}(\si)$.
We apply the induction assumption to the link $\Si_v\si$ at a vertex $v$ of $\tau$.
Note that $\D_{\Si_v\tau}\Si_v\si =\Si_v \D_{\tau}\si$.
Since $\oa{vx}\in\inte_{\Si_v\tau}(\Si_v\si)$,
the nearest point projection $\bar\de$ of this direction to $\Si_v\tau$ is contained in $\inte(\Si_v\tau)$
and has angle $<\pihalf$ with $\oa{vx}$.
It follows that the nearest point projection $p$ of $x$ to $\tau$ 
is different from $v$ and lies on the arc in direction $\bar\de$,
$\oa{vp}=\bar\de$. 
In particular, it is not contained in a face of $\tau$ with vertex $v$.
Letting run $v$ through the vertices of $\tau$,
we conclude that $p\in\inte(\tau)$. 
\qed

\medskip
As a consequence of the lemma, 
the nearest point projection $\inte_{\tau}(\si)\to\tau$ 
agrees with the nearest point projection $\inte_{\tau}(\si)\to s$ 
to the geodesic sphere $s\subset S$ spanned by $\tau$
(i.e.\ containing $\tau$ as a top-dimensional subset),
and its image equals $\inte(\tau)$.

\subsubsection{Spherical Coxeter complexes}

A {\em spherical Coxeter complex} $(\amod,W)$ consists of a unit sphere (in a euclidean space) $\amod$
and a finite reflection group $W$ acting isometrically on $\amod$.
We will refer to $\amod$ as the {\em model apartment}
(because it will serve as the model for apartments in spherical buildings, see below).

A {\em wall} in $\amod$ is the fixed point set of a reflection in $W$. 
A {\em half-apartment} is a closed hemisphere in $\amod$ bounded by a wall. 
A {\em singular sphere} in $\amod$ is an intersection of walls.

A {\em chamber} in $\amod$ is the closure of a connected component of the complement of the union of the walls.
The group $W$ acts transitively on the set of chambers. 
The chambers are simplices with diameter $\leq\pihalf$
iff $W$ fixes no point in $\amod$, equivalently, 
the Coxeter complex does not split off a spherical join factor
(in the category of Coxeter complexes).
In this case,
the collection of chambers defines on $\amod$ the structure of a simplicial complex, 
the simplices being intersections of chambers. 

Every chamber is a fundamental domain for the action $W\acts\amod$.
The {\em spherical model chamber} can be defined as the quotient $\simod=\amod/W$.
We identify it with a chamber in the model apartment, $\simod\subset\amod$, 
which we refer to as the {\em fundamental chamber}.

We call the natural projection
$$ \theta:\amod\to\amod/W\cong\simod $$
the {\em type} map for $\amod$.
It restricts to an isometry on every chamber.
A {\em face type} is a face of $\simod$.
The type of a simplex $\bar\tau\subset\amod$ is then defined as $\theta(\bar\tau)$.
Throughout the paper,
we will use the notation 
$\taumod,\taumod',\numod,\numod',\ldots$ for face types.
Furthermore,
we will denote by $\Wt\leq W$ the stabilizer of the face type $\taumod\subseteq\simod$.

The {\em longest element} of the Weyl group is the unique 
element $w_0\in W$ sending $\simod$ to the opposite chamber $-\simod$. 
The {\em standard involution} (also known as the {\em opposition involution}) 
of the model chamber is given by 
$\iota:=-w_0:\simod\to\simod$.

\subsubsection{Spherical buildings}
\label{sec:sphbld}

A {\em spherical building} modeled on a Coxeter complex $(\amod,W)$ 
is a CAT(1) metric space $\B$ 
equipped with a collection of isometric embeddings $\kappa:\amod\to\B$, called {\em charts}.
The image of a chart is an {\em apartment} in $\B$.
One requires that any two points are contained in an apartment
and that the coordinate changes between charts are induced by isometries in $W$.
(The precise axioms can be found e.g. in \cite{qirigid} and \cite{habil}.)
We will use the notation $\angle$ for the metric on $\B$. 

We assume that $W$ fixes no point, equivalently, that $\simod$ is a simplex with diameter $\leq\pihalf$.

Via the atlas of charts,
the spherical building inherits from the spherical Coxeter complex 
a natural structure of a {\em simplicial complex}
where the simplices are the images of the simplices in the model apartment. 
As already mentioned,
the images of the charts are called {\em apartments}. 
Accordingly, the images of chambers (walls, half-apartments, singular spheres) in $\amod$ 
are called chambers (walls, half-apartments, singular spheres) in the building. 
The codimension one faces are called {\em panels}. 
The {\em interior} $\inte(\tau)$ of a simplex $\tau$ 
is obtained by removing all proper faces; 
the interiors of simplices are called {\em open simplices}. 
The simplex {\em spanned} by a point 
is the smallest simplex containing it, equivalently, the simplex containing the point in its interior.
We will sometimes denote the simplex spanned by $\xi$ by $\tau_{\xi}$.

A spherical building is {\em thick} if every wall is the bounds at least three half-apartments,
equivalently, if every panel is adjacent to (i.e.\ contained in the boundary of) at least three chambers.
One can always pass to a thick spherical building structure by reducing the Weyl group,
thereby coarsifying the simplicial structure.

The {\em space of directions} $\Si_\xi \B$ at a point $\xi\in \B$ is the space of germs $\oa{\xi\eta}$ 
of nondegenerate geodesic segments $\xi\eta\subset \B$, 
equipped with the natural angle metric $\angle_{\xi}$. 
Two segments $\xi\eta$ and $\xi\eta'$ represent the same direction in $\Si_\xi \B$, $\oa{\xi\eta}=\oa{\xi\eta'}$, iff they initially agree.
The space of directions is again a spherical building. 

A subset $C\subset \B$ is called ($\pi$-){\em convex} 
if for any two points $\xi, \eta\in C$ with distance $\angle(\xi,\eta)<\pi$
the (unique) geodesic $\xi\eta$ connecting $\xi$ and $\eta$ in $\B$ is contained in $C$. 

Due to the compatibility of charts,
i.e.\ the property of the building atlas that the coordinate changes are induced by isometries in $W$,
there is a well-defined {\em type} map 
\begin{equation*}
\theta: \B\to\simod .
\end{equation*}
It is 1-Lipschitz and restricts to an isometry on every chamber $\si\subset \B$.
We call the inverse 
$\kappa_{\si}=(\theta|_{\si})^{-1}: \si_{mod}\to\si$
the {\em chart} of the chamber $\si$. 
For a simplex $\tau\subset \B$, 
we call the face $\theta(\tau)\subseteq\simod$ 
the {\em type} of the simplex and 
$\kappa_{\tau}=(\theta|_{\tau})^{-1}:\theta(\tau)\to\tau$
its {\em chart}. 
We define the {\em type} of a point $\xi\in \B$
as its image $\theta(\xi)\in\simod$. 
A point $\xi\in \B$ is called {\em regular} 
if its type is an interior point of $\simod$, $\xi\in\inte(\simod)$, 
and {\em singular} otherwise. 

We will sometimes say that 
a singular sphere has type $\taumod$ if it contains a top-dimensional simplex of type $\taumod$.
(A singular sphere has in general several types.)

For a singular sphere $s\subset \B$,
we define $\B(s)\subset \B$ as the union of all apartments containing $s$.
It is a convex subset and splits off $s$ as a spherical join factor.
Moreover, $\B(s)$ is a subbuilding, 
i.e.\ it inherits from $\B$ a spherical building structure modeled on the same Coxeter complex;
the apartments of $\B(s)$ are precisely the apartments of $\B$ containing $s$.
This building structure is however not thick, except in degenerate cases. 
In order to pass to a thick spherical building structure,
take a maximal atlas of charts $\kappa:\amod\to \B(s)$ 
for which the maps $\kappa^{-1}|_s:s\to\amod$ coincide,
and reduce the Weyl group to the pointwise stabilizer of $s$ in $W$.

Two points $\xi, \hat\xi\in \B$ are {\em antipodal} or {\em opposite}
if $\angle(\xi,\hat\xi)=\pi$,
equivalently, 
if they are antipodal in one (every)  apartment containing them.
We then define the singular sphere $s(\xi,\hat\xi)\subset \B$ {\em spanned} by the points $\xi, \hat\xi$
as the smallest singular sphere containing them.
Moreover,
we define the {\em suspension} $\B(\xi,\hat\xi)\subset \B$ of $\{\xi, \hat\xi\}$ 
as the union of all geodesics connecting $\xi$ and $\hat\xi$,
equivalently, as the union of all apartments containing $\xi$ and $\hat\xi$.
Then $\B(\xi,\hat\xi)=\B(s(\xi,\hat\xi))$.
As above, a thick spherical building structure on $\B(\xi,\hat\xi)$ is obtained
by taking all charts $\kappa:\amod\to \B(\xi,\hat\xi)$ so that $\kappa^{-1}(\xi)=\theta(\xi)\in\simod$,
and reducing the Weyl group to the stabilizer of $\theta(\xi)$ in $W$.

Similarly, one defines {\em antipodal} or {\em opposite faces} $\tau,\hat\tau\subset \B$ 
as faces which are antipodal in the apartments containing them both,
equivalently, whose interiors contain a pair of antipodal points $\xi\in\inte(\tau)$ and $\hat\xi\in\inte(\hat\tau)$. 
We define the singular sphere $s(\tau,\hat\tau)\subset \B$
{\em spanned} by the simplices $\tau,\hat\tau$
again as the smallest singular sphere containing them,
and the suspension $\B(\tau,\hat\tau)$ as the union of all apartments containing $\tau \cup \hat\tau$;
then $s(\tau,\hat\tau)= s(\xi,\hat\xi)$ and $\B(\tau,\hat\tau)= \B(\xi,\hat\xi)$. 

\medskip
We will need some facts about antipodes.

Recall that in a spherical building $\B$
every point $\xi\in \B$ has an antipode in every apartment $a\subset \B$,
and hence for every simplex $\tau\subset \B$ there exists an opposite simplex $\hat\tau\subset a$,
cf.\ e.g.\ the first part of \cite[Lemma 3.10.2]{qirigid}. 
We need the more precise statement
that a point has {\em several} antipodes in an apartment
unless it lies itself in this apartment:

\begin{lem}
\label{lem:oneantip}
Suppose that $\xi\in \B$ has only one antipode in the apartment $a\subset \B$. 
Then $\xi\in a$.
\end{lem}
\proof
Suppose that $\xi\not\in a$ and let $\hat\xi\in a$ be an antipode of $\xi$. 
We choose a``generic'' segment $\xi\hat\xi$ of length $\pi$ tangent to $a$ at $\hat\xi$ as follows.
The suspension $\B(\xi,\hat\xi)$ contains an apartment $a'$ 
with the same unit tangent sphere at $\hat\xi$, $\Si_{\hat\xi}a'=\Si_{\hat\xi}a$. 
Inside $a'$ there exists a segment $\xi\hat\xi$ whose interior does not intersect simplices of codimension $\geq2$. 
Hence $\hat\xi\xi$ leaves $a$ at an interior point $\eta\neq\xi,\hat\xi$ of a panel $\pi\subset a$,
i.e.\ $a\cap\xi\hat\xi=\eta\hat\xi$ and $\pi\cap\xi\hat\xi=\eta$,
and $\eta\xi$ initially lies in a chamber adjacent to $\pi$ but not contained in $a$. 
Let $s\subset a$ be the wall containing $\pi$. 
By reflecting $\hat\xi$ at $s$, one obtains a second antipode for $\xi$ in $a$.
\qed

\medskip
In thick buildings, simplices can be represented as intersections of apartments:
\begin{lem}
\label{lem:intapts}
In a thick spherical building $\B$, 
any simplex $\tau\subset \B$ equals the intersection of the apartments containing it.
\end{lem}
\proof
Since every simplex is an intersection of chambers,
we are reduced to the case when $\tau$ is a chamber. 
Furthermore,
since every chamber is an intersection of half-apartments,
we are reduced to the corresponding assertion for half-apartments.
The latter holds by thickness.
\qed

\subsection{Hadamard manifolds}
\label{sec:had}

In this section only,
$X$ denotes a Hadamard manifold, 
i.e.\ a simply connected complete Riemannian manifold 
with nonpositive sectional curvature.  We will use the notation $\Isom(X)$ for the full isometry group of $X$. 

Any two points in $X$ are connected by a unique geodesic segment. 
We will use the notation $xy$ for the oriented geodesic 
segment connecting $x$ to $y$.

For points $x\neq y,z$ we denote by $\angle_x(y,z)$ 
the angle of the geodesic segments $xy$ and $xz$.
Furthermore,
we denote by $\Si_xX$ the {\em space of directions} of $X$ at $x$
equipped with the angle metric $\angle_x$.
It coincides with the unit tangent sphere at $x$.

A basic feature of Hadamard manifolds is the {\em convexity of the distance function}: 
Given any pair of geodesics $c_1(t), c_2(t)$ in $X$, the function
$t\mapsto d(c_1(t), c_2(t))$ 
is convex. 

Two geodesic rays $\rho_1,\rho_2:[0,+\infty)\to X$ 
are called {\em asymptotic} if the convex function 
$t\mapsto d(\rho_1(t),\rho_2(t))$ on $[0,+\infty)$ is bounded, 
and they are called {\em strongly asymptotic} if 
$d(\rho_1(t),\rho_2(t))\to 0$ as $t\to+\infty$. 

Two geodesic lines $l_1, l_2\subset X$ are {\em parallel} 
if they have finite Hausdorff distance. 
Equivalently,  $l_1\cup l_2$ bounds a flat strip in $X$. 

The {\em ideal} or {\em visual boundary} $\geo X$ of $X$
is the set of asymptote classes of geodesic rays in $X$. 
Points in $\geo X$ are called {\em ideal points}. 
For $x\in X$ and $\xi\in \geo X$ we denote by $x\xi$ the unique geodesic ray 
emanating from $x$ and asymptotic to $\xi$, 
i.e.\ representing the ideal point $\xi$. 
There are natural identifications $\log_x:\geo X\to\Si_xX$
sending the ideal point $\xi$ to the direction $\oa{x\xi}$.

The {\em cone} or {\em visual topology} on $\geo X$ 
is characterized by the property that the maps $\log_x$ are homeomorphisms
with respect to it.
Thus, $\geo X$ is homeomorphic to the sphere of dimension $\dim(X)-1$. 
The visual topology has a natural extension to $\ol{X}=X\sqcup \geo X$ 
which can be described as follows in terms of sequential convergence:
A sequence $(x_n)$ in $\ol X$ converges to an ideal point $\xi\in \geo X$ 
iff, for some (any) base point $x\in X$,
the sequence of geodesic segments or rays $xx_n$ 
converges to the ray $x\xi$ 
(in the pointed Hausdorff topology with base points at $x$).
This topology makes $\ol{X}$ into a closed ball. 
We define the {\em visual boundary} of a subset $A\subset X$ 
as the set $\geo A=\bar A\cap\geo X$ 
of its accumulation points at infinity. 

The visual boundary $\geo X$ carries the natural 
{\em Tits angle metric} $\tangle$ 
defined as 
\begin{equation*}
\tangle(\xi,\eta)= \sup_{x\in X} \angle_{x}(\xi,\eta) ,
\end{equation*}
where $\angle_{x}(\xi,\eta)$ 
is the angle between the geodesic rays $x\xi$ and $x\eta$. 
The {\em Tits boundary} $\tits X$ is the metric space $(\geo X,\tangle)$. 
The Tits metric is lower semicontinuous with respect to the visual topology 
and, accordingly, 
the {\em Tits topology} induced by the Tits metric 
is finer than the visual topology. 
It is discrete if there is an upper negative curvature bound, 
and becomes nondiscrete if $X$ contains nondegenerate flat sectors.
For instance, 
the Tits boundary of flat $r$-space is the unit $(r-1)$-sphere, 
$\tits\R^r\cong S^{r-1}(1)$. 
An isometric embedding $X\to Y$ of Hadamard manifolds
induces an isometric embedding $\tits X\to\tits Y$ of their Tits boundaries.

\medskip
Let $\xi\in \geo X$ be an ideal point. 
For a geodesic ray $\rho:[0,+\infty)\to X$ asymptotic to $\xi$ 
one defines the {\em Busemann function} $b_\xi$ on $X$ as
the uniform monotonic limit
\begin{equation*}
b_\xi(x)=\lim_{t\to+\infty} (d(x, \rho(t))-t).
\end{equation*}
Along the ray, we have  
\begin{equation*}
b_\xi(\rho(t))=-t. 
\end{equation*}
Altering the ray $\rho$ changes $b_\xi$ by an additive constant. 
{The point at infinity $\xi$ thus determines $b_\xi$ up to an additive constant. 
To remove this ambiguity,
given $x\in X$, we define 
$b_{\xi,x}$ to be the Busemann function $b_{\xi,x}$ normalized at the point $x$ by $b_{\xi,x}(x)=0$.}

The Busemann function $b_\xi$ is convex, 1-Lipschitz 
and measures the relative distance from the ideal point $\xi$. 
The sublevel sets 
\begin{equation*}
\Hb_{\xi,x}:=\{b_{\xi}\leq b_{\xi}(x)\}\subset X
\end{equation*}
are called (closed) {\em horoballs} centered at $\xi$. 
As sublevel sets of convex functions, they are convex. 
The visual boundaries of horoballs are $\pihalf$-balls at infinity 
with respect to the Tits metric, 
\begin{equation*}
\geo \Hb_{\xi,x}= \ol B(\xi,\pihalf):=\{\tangle(\xi, \cdot)\le \pi/2\}\subset\geo X .
\end{equation*}
The level sets 
\begin{equation*}
\Hs_{\xi,x}:=\{b_{\xi}= b_{\xi}(x)\} =\D \Hb_{\xi,x}
\end{equation*}
are called {\em horospheres} centered at $\xi$. 

As convex Lipschitz functions,
Busemann functions are {\em asymptotically linear} along rays.
If $\rho:[0,+\infty)\to X$ is a geodesic ray 
asymptotic to $\eta\in\geo X$, $\rho(+\infty)=\eta$,
then 
\begin{equation*}
\lim_{t\to+\infty} \frac{b_\xi(\rho(t))}{t}=-\cos \tangle(\xi, \eta). 
\end{equation*}

\subsection{Symmetric spaces of noncompact type: basic concepts}
\label{sec:symmbas}

In this section, we go through some well known material 
and establish notation.
Standard references are \cite{Eberlein} and \cite{BGS}. 

A symmetric space, denoted by $X$ throughout this paper, 
is said to be of {\em noncompact type} 
if it is nonpositively curved and has no euclidean factor. 
In particular, it is a Hadamard manifold. 
We will write the symmetric space as 
$$ X=G/K $$
where $G$ is a connected\footnote{What is really needed is a weaker property than connectedness,
namely that $G$ has finitely many connected components 
and acts on the Tits building of $X$ by (type preserving) automorphisms.
The latter is equivalent to the triviality of the $G$-action on the model chamber $\simod$, equivalently, on the Dynkin diagram.
Under this assumption, the theory of discrete subgroups presented in this paper goes through unchanged.}
semisimple Lie group with finite center acting isometrically and transitively on $X$, 
and $K<G$ is a maximal compact subgroup.
The natural epimorphism $G\to\Isom(X)_o$ then has compact kernel. 
Every connected semisimple Lie group with finite center occurs in this way.
The Lie group $G$ carries a natural structure of a real algebraic group. 

By the definition of symmetric spaces,
in every point $x\in X$ there is a 
{\em point reflection} or {\em Cartan involution}, that is,
an isometry $\si_x$ which fixes $x$ and has differential $-\id_{T_xX}$ in $x$. 

A {\em transvection} of $X$ is an isometry 
which is the product $\si_{x'}\si_{x}$ of two point reflections; 
it preserves the oriented geodesic through $x$ and $x'$ 
and the parallel vector fields along it. 
The transvections preserving a geodesic line $c(t)$ 
form a one parameter subgroup $(T^c_t)$ of $\Isom(X)_o$ 
where $T^c_t$ denotes the transvection mapping 
$c(s)\mapsto c(s+t)$. 

An isometry $\phi$ of $X$ is called {\em axial} 
if it preserves a geodesic $l$ and does not fix $l$ pointwise. 
Thus, $\phi$ acts as a nontrivial translation on $l$. 
(Note that an axial isometry need not be a transvection.) 
The geodesic $l$ is called an {\em axis} of $\phi$. 
Axes are in general not unique,
but they are parallel to each other. 
For each axial isometry $\phi$, the displacement function 
$x\mapsto d(x, \phi(x))$ on $X$ 
attains its minimum on the convex subset of $X$ which is the union of axes of $\phi$. 
An isometry $\phi$ of $X$ is {\em parabolic} if
$$
\inf_{x\in X} d(x, \phi(x))=0
$$
but $g$ does not fix a point in $X$. 
Isometries fixing points are called {\em elliptic}.

A {\em flat} in $X$ is a complete totally geodesic flat submanifold, 
equivalently, 
a convex subset isometric to a euclidean space. 
A maximal flat in $X$ is a flat 
which is not contained in any larger flat; 
we will use the notation $F$ for maximal flats. 
The group $\Isom(X)_o$ acts transitively on the set of maximal flats; 
the common dimension of maximal  flats is called the {\em rank} of $X$. 
The space $X$ has rank one if and only if it has 
strictly negative sectional curvature. 

A maximal flat $F$ 
is preserved by all transvections along geodesic lines contained in it. 
In general, there exist nontrivial isometries of $X$ fixing $F$ pointwise. 
The subgroup of isometries of $F$ 
which are induced by elements of $G$ 
is isomorphic to a semidirect product $W_{aff}:= \R^r \rtimes W$, 
the {\em affine Weyl group}, 
where $r$ is the rank of $X$. 
The subgroup $\R^r$ acts simply transitively on $F$ by translations. 
The linear part $W$ 
is a finite reflection group, called the {\em Weyl group} of $G$ and $X$. 
Since maximal flats are equivalent modulo $G$, 
the action $W_{aff}\acts F$ is well-defined up to isometric conjugacy. 

We will think of the Weyl group as 
acting on a {\em model flat} $\Fmod\cong \R^r$ fixing the origin $0\in\Fmod$, 
and on its visual boundary sphere at infinity, 
the {\em model apartment} $a_{mod}=\tits\Fmod\cong S^{r-1}$. 
The pair $(\amod,W)$ is the {\em spherical Coxeter complex} associated to $G$.
We identify the {\em euclidean model Weyl chamber} $\De$ 
with the complete cone $V(0,\simod)\subset\Fmod$
with tip in the origin and visual boundary the {\em spherical model Weyl chamber} $\simod\subset\amod$.

For every maximal flat $F\subset X$, 
we have an induced Tits isometric embedding $\geo F\subset\geo X$ 
of its visual boundary sphere. 
The natural identification $F\cong\Fmod$,
unique up to the action of $W_{aff}$, 
induces a natural identification $\geo F\cong a_{mod}$,
unique up to the action of $W$.

The Coxeter complex structure on $a_{mod}$ 
induces simplicial structures on the visual boundary spheres $\geo F$ of the maximal flats $F\subset X$.
The spheres $\geo F$ cover $\geo X$,
and their simplicial structures are compatible
(i.e.\ the intersections are simplicial and the simplicial structures on the intersections agree).
One thus obtains a $G$-invariant 
piecewise spherical {\em simplicial structure} on $\geo X$ 
which makes $\geo X$ into a {\em thick spherical building} and, 
also taking into account the visual topology, 
into a {\em topological} spherical building. 
It is called the {\em spherical} or {\em Tits building} $\tits X$ 
associated to $X$. 
The Tits metric is the path metric 
with respect to the piecewise spherical structure, unless $\rank(X)=1$, in which case $\tits X$ is discrete with distance $\pi$ between distinct points.  
We will sometimes refer to the simplices in $\tits X$ also as {\em faces}.
The visual boundaries of the maximal flats in $X$ 
are precisely the {\em apartments} in $\geo X$,
which in turn are precisely the convex subsets 
isometric, with respect to the Tits metric, to the unit sphere $S^{r-1}$.

We call a flat $f\subset X$ {\em singular} if it is the intersection of maximal flats.
Its visual boundary $\geo f$ is then a {\em singular sphere} in $\geo X$.

\medskip
We define the {\em Weyl sector} $V=V(x,\tau)\subset X$ 
with tip $x$ and asymptotic to a simplex $\tau\subset\geo X$ 
as the union of rays $x\xi$ for the ideal points $\xi\in\tau$. 
Weyl sectors are contained in flats;
they are isometric images of Weyl sectors $V(0, \tau_{mod})\subset\De$ 
under charts $\Fmod\to X$. 
These apartment charts restrict to canonical {\em sector charts} 
$\kappa_{x,\tau}=\kappa_{V(x,\tau)}: V(0,\taumod)\to V(x,\tau)$;
at infinity, they induce simplex charts,
$\geo\kappa_{x,\tau}=\kappa_{\tau}$.

If $\si\subset\geo X$ is a chamber,
the sector $V(x,\si)$ is a {\em euclidean Weyl chamber}.

\medskip
For a flat $f\subset X$,
the {\em parallel set} $P(f)\subset X$ is the union of all flats $f'\subset X$ {\em parallel} to $f$,
equivalently, with the same visual boundary sphere $\geo f'=\geo f$.
The parallel set is a symmetric subspace 
and splits as the metric product 
\begin{equation}
\label{eq:parsplit}
P(f)\cong f\times CS(f)
\end{equation}
of $f$ and a symmetric space $CS(f)$ called the {\em cross section}. 
The latter has no euclidean factor iff $f$ is singular. 
Accordingly, 
the Tits boundary metrically decomposes as the spherical join
\begin{equation}
\label{eq:bdparsersusp}
\tits P(f)\cong\tits f\circ\tits CS(f) .
\end{equation}
It coincides with the subbuilding $(\tits X)(\geo f)\subset\tits X$ 
consisting of the union of all apartments in $\geo X$ containing $\geo f$, 
see section \ref{sec:sphbld}. 

For a singular sphere $s\subset\geo X$,
we define the parallel set $P(s)\subset X$
as the union of the (necessarily singular) flats $f\subset X$ with visual boundary sphere $\geo f=s$,
i.e.\ $P(s)=P(f)$;
we denote its cross section by $CS(s)$.
For a pair of opposite points $\xi,\hat\xi\in\geo X$,
we define 
$P(\xi,\hat\xi)\subset X$ as the parallel set of the singular sphere $s(\xi,\hat\xi)\subset\geo X$ spanned by them, 
$P(\xi,\hat\xi)=P(s(\xi,\hat\xi))$.
Similarly,
for a pair of opposite simplices $\tau,\hat\tau\subset\geo X$,
we define $P(\tau_-, \tau_+)=P(s(\tau_-,\tau_+))$.

\medskip
The action $G\acts\geo X$ on ideal points is not transitive if $\rank(X)\geq2$.
However, every $G$-orbit meets every chamber exactly once. 
The quotient is naturally identified with the spherical model chamber,
and the projection 
\begin{equation*}
\theta:\geo X\to\geo X/G \cong\simod
\end{equation*}
is the {\em type} map, 
cf.\ section~\ref{sec:sphbld}.

A nondegenerate geodesic segment $xy\subset X$ is called {\em regular} 
if the unique geodesic ray $x\xi$ extending $xy$ is asymptotic to a regular ideal point $\xi\in \geo X$.

Two ideal points $\xi,\eta\in\geo X$ are  {\em antipodal}, $\tangle(\xi, \eta)=\pi$, 
iff there exists a geodesic line $l\subset X$ asymptotic to them, $\geo l=\{\xi, \eta\}$. 
Their types are then related by $\theta(\xi_2)=\iota(\theta(\xi_1))$.

We say that two simplices $\tau_1,\tau_2\subset\geo X$ 
are $x$-{\em antipodal} or $x$-{\em opposite} 
if $\tau_2=\si_x\tau_1$, 
using the induced action of the point reflection $\si_x$ on $\geo X$.
Two simplices $\tau_1,\tau_2$ are opposite
iff they are $x$-opposite for some point $x\in X$.
Their types are then related by $\theta(\tau_2)=\iota(\theta(\tau_1))$.
We will frequently use the notation $\tau, \hat\tau$ and $\tau_{\pm}$ for pairs of antipodal simplices. 
A pair of opposite chambers $\si_{\pm}$ is contained in a unique apartment,  
which we will denote by $a(\si_-,\si_+)$.
It is the visual boundary of a unique maximal flat $F(\si_-, \si_+)\subset X$.

We will sometimes say that 
a singular flat $f\subset X$ has {\em type $\taumod$} if its visual boundary $\geo f$ has type $\taumod$,
i.e.\ contains a top-dimensional simplex of type $\taumod$.
(A singular flat has in general several types.)
The set ${\mathcal F}_{\taumod}$ of singular flats of type $\taumod$ is a homogeneous $G$-manifold.
The flats of type $\simod$ are the maximal flats and we denote ${\mathcal F}={\mathcal F}_{\simod}$. 
A family of flats in ${\mathcal F}_{\taumod}$ is {\em bounded}  
if these flats intersect a fixed bounded subset of $X$.

Also, we will sometimes call the parallel set $P(s)$ of a singular sphere $\subset\geo X$ 
of {\em type $\taumod$} or a {\em $\taumod$-parallel set}
if $s$ has type $\taumod$.

\medskip
The stabilizers $P_{\tau}<G$ of the simplices $\tau\subset\geo X$ 
are the {\em parabolic subgroups} of $G$. 
The space $\Flagt$ of simplices of type $\taumod$ is called a (generalized) {\em (partial) flag manifold}. 
The action $G\acts\Flagt$ is transitive
and we can write the flag manifold as a quotient $\Flagt \cong G/P_{\taumod}$, 
where $P_{\taumod}$ stands for a parabolic subgroup 
in the conjugacy class of parabolic subgroups $P_{\tau}$ of type $\theta(\tau)=\taumod$.
Flag manifolds are compact smooth manifolds; 
they admit natural structures of projective real algebraic varieties
(see e.g.\ \cite[p.\ 160]{Jantzen}). 
The topology on flag manifolds induced by the visual topology on $\geo X$
agrees with their manifold topology as homogeneous $G$-spaces.  
For ideal points $\xi\in\geo X$ with type $\theta(\xi)\in\inte(\taumod)$, 
there is a natural $G$-equivariant homeomorphic identification 
of the $G$-orbit $G\xi\subset\geo X$ with $\Flagt$
by assigning to the point $g\xi$ the (unique) simplex of type $\taumod$ containing it. 

The flag manifolds $\Flagt$ and $\Flagit$ are {\em opposite} in the sense 
that the simplices opposite to simplices of type $\taumod$ have type $\iota\taumod$.
To ease notation, we will denote the pair of opposite flag manifolds also by $\Flagpmt$
whenever convenient,
i.e.\ we put $\Flagpt:=\Flagt$ and $\Flagmt:=\Flagit$.
The latter is also reasonable, because 
the simplices $-\taumod, \iota\taumod\subset \amod$ lie in the same $W$-orbit,
i.e.\ $-\taumod$ has type $\iota\taumod$.
(Here we extend the notion of type
to the model apartment, 
defining the type of a simplex in $\amod$ as its image under the natural quotient projection 
$\amod\to\amod/W\cong\simod$.)
Similarly,
we will use the notation $P_{\pm\taumod}$ for a pair of parabolic subgroups 
fixing opposite simplices in $\Flagpmt$.

The stabilizers $B_{\si}<G$ of the chambers $\sigma\subset \geo X$ 
are the {\em minimal parabolic subgroups}\footnote{{When the group $G$ is complex, 
the minimal parabolic subgroups are the {\em Borel subgroups}, 
which is why we use the notation $B$ for these subgroups.}}
of $G$;
they are conjugate.
The space $\DF X:=\Flags$ of chambers is called the (generalized) {\em full flag manifold}
or {\em Furstenberg boundary} of $X$,
and we can write $\DF X=G/B$,
where again $B$ stands for a minimal parabolic subgroup.

For a simplex $\hat\tau\in\Flagit$
we define the {\em open Schubert stratum} $C(\hat\tau)\subset \Flagt$ as the subset of simplices opposite to $\hat\tau$;
it is the open and {\em dense} $P_{\hat\tau}$-orbit.
With respect to the algebraic structure on $\Flagt$,
it is Zariski open,
i.e.\ its complement is a proper subvariety.

We note that, if $\rank(X)=1$,
then there is only one flag manifold, namely $\geo X$,
and the open Schubert strata are the complements of points.

\subsection{Stars, cones and diamonds}\label{sec:stars}

\subsubsection{Stars and suspensions}
\label{sec:star}

We first work inside the spherical model chamber $\simod$.

We recall from section~\ref{sec:sphgeom} that,
for a face type $\taumod\subseteq\simod$,
the {\em $\taumod$-boundary} $\Dt\simod$ of $\simod$ 
is the union of the faces of $\simod$ which do not contain $\taumod$.
The {\em $\taumod$-interior} $\inte_{\taumod}(\simod)$ of $\simod$ 
is the union of the open faces of $\simod$ whose closure contains $\taumod$. 
There is the decomposition
$$\simod=\inte_{\taumod}(\simod)\sqcup\Dt\simod .$$
In particular, $\inte_{\simod}(\simod)=\inte{\simod}$ and 
$\Ds\simod=\D\simod$.

We say that a type in $\simod$ is {\em $\taumod$-regular} if it lies in $\inte_{\taumod}(\simod)$.

\medskip
Now let $\B$ be a spherical building.
As before, we assume that $\diam(\simod)\leq\pihalf$.

A point $\xi\in \B$ is called {\em $\taumod$-regular} if its type is, $\theta(\xi)\in\inte_{\taumod}(\simod)$.
We will {\em quantify} $\taumod$-regularity as follows:
Given a compact subset $\Theta\subset\inte_{\taumod}(\simod)$,
we will say that a $\taumod$-regular point $\xi\in \B$ is {\em $\Theta$-regular} if $\theta(\xi)\in\Theta$.

It will often be natural to impose a convexity property on $\Theta$:
\begin{dfn}[Weyl convex]
A subset $\Theta\subseteq\simod$ is {\em $\taumod$-Weyl convex} 
if its symmetrization $\Wt\Theta\subset\amod$ is convex.
\end{dfn}

Let $\tau\subset \B$ be a simplex of type $\taumod$.
The {\em $\taumod$-star} $\st(\tau)\subset \B$ is the union of all chambers containing $\tau$. 
Its boundary $\D\st(\tau)$ is the union of all simplices in $\st(\tau)$ which do not contain $\tau$;
it consists of the points in $\st(\tau)$ with type in $\Dt\simod$. 
The {\em open $\taumod$-star} $\ost(\tau)$ is the complement $\ost(\tau)=\st(\tau)-\D\st(\tau)$;
it consists of the $\taumod$-regular points in $\st(\tau)$
and is open in $\B$.
For any simplex $\hat\tau$ opposite to $\tau$,
the star $\st(\tau)$ is contained in the suspension $\B(\tau,\hat\tau)$. 

Furthermore,
we define the {\em $\Theta$-star} $\stTh(\tau)\subset\ost(\tau)$
as the subset of points with type $\Theta$,
that is, $\stTh(\tau)=\st(\tau)\cap\theta^{-1}(\Theta)$. 

We will use the following {\em separation} property:
If $\angle(\Theta,\Dt\simod)\geq\eps>0$, 
then $\ost(\tau)$ contains the open $\eps$-neighborhood of $\stTh(\tau)$.

Note that for chambers $\si$ we have $\st(\si)=\si$ and $\ost(\si)=\inte(\si)$.

\medskip
The next result implies that stars are {\em convex}:
\begin{lem}
[Convexity of stars]
\label{lem:starconv}
(i) 
$\st(\tau)$ is an intersection of simplicial $\pihalf$-balls.

(ii)
For any simplex $\hat\tau$ opposite to $\tau$, 
the star $\st(\tau)$ is an intersection 
of the suspension $\B(\tau,\hat\tau)$
with simplicial $\pihalf$-balls 
containing $\st(\tau)$ and centered at points in $\B(\tau,\hat\tau)$.
\end{lem}
\proof
(i) 
Let $\si\not\subset\st(\tau)$ be a chamber,
and let $a$ be an apartment containing $\si$ and $\tau$. 
We can separate $\si$ and $\st(\tau)\cap a$ by a wall in $a$,
i.e.\ there exists a half-apartment $h\subset a$ which contains $\st(\tau)\cap a$ but not $\si$. 
Indeed, choose points $\xi\in\inte(\tau)$ and $\eta\in\inte(\si)$
such that the segment $\xi\eta$ intersects $\D\si$ in a panel,
and take the wall containing this panel.
The simplicial $\pihalf$-ball with the same center as $h$ then contains $\st(\tau)$ but not $\si$. 

(ii) Note first that $\st(\tau)\subset \B(\tau,\hat\tau)$. 
Then we argue as in part (i), observing that 
if $\si\subset \B(\tau,\hat\tau)$ then $a$ can be chosen inside $\B(\tau,\hat\tau)$. 
\qed

\medskip
We extend convexity to $\Theta$-stars:
\begin{lem}[Convexity of $\Theta$-stars]
\label{lem:thetastarconv}
Let $\Theta\subseteq\simod$ be $\taumod$-Weyl convex,
and let $\tau$ be a simplex of type $\taumod$.
Then $\stTh(\tau)$ is an intersection of $\pihalf$-balls.
\end{lem}
\proof
For any apartment $a\supset\tau$, 
the intersection $\stTh(\tau)\cap a$ is convex, 
as a consequence of the Weyl convexity of $\Theta$.

Let $\zeta\in \B$.
Every point in $\stTh(\tau)$ lies in an apartment $a\supset\tau,\zeta$.

For any two apartments $a,a'\supset\tau,\zeta$ 
there exists an isometry $a\to a'$ fixing $\tau$ and $\zeta$. 
(This follows from the compatibility of apartment charts axiom in the definition of spherical buildings.)
It carries $\stTh(\tau)\cap a$ to $\stTh(\tau)\cap a'$.
Hence, $\ol B(\zeta,\pihalf)$ contains the first intersection iff it contains the second.
Letting $a'$ vary,
it follows that $\ol B(\zeta,\pihalf)$ contains $\stTh(\tau)$ iff it contains $\stTh(\tau)\cap a$.

Let $\xi\not\in\stTh(\tau)$.
Then there is an apartment $a\supset\tau,\xi$ 
and, due to the convexity of $\stTh(\tau)\cap a$,
a point $\zeta\in a$ such that $\ol B(\zeta,\pihalf)$ contains $\stTh(\tau)\cap a$ but not $\xi$.
By the above, $\stTh(\tau)\subseteq\ol B(\zeta,\pihalf)$.
\qed

\medskip
In the following, we restrict ourselves to the case $\B=\geo X$
and, besides the metric, also take into account the visual topology on the flag manifolds $\Flagt$. 
The discussion readily generalizes to arbitrary {\em topological} spherical buildings. 

The {\em $\taumod$-regular part} $\geot X$
of the visual boundary equals the union of the open $\taumod$-stars.
The natural projection
\begin{equation}
\label{eq:ostfib}
\geot X=\bigcup_{\tau\in\Flagt}\ost(\tau)\to\Flagt
\end{equation}
is a fiber bundle.

Let $\tau\in\Flagt$ and let $\hat\tau$ be opposite to $\tau$.
Then $\tau$ is the {\em only} simplex in $\B(\tau,\hat\tau)$ which is opposite to $\hat\tau$.
In other words, 
the closed subset 
\begin{equation}
\label{eq:oppinsusp}
\{ \tau'\in\Flagt: \tau'\subset \B(\tau,\hat\tau) \}
\end{equation}
intersects the open Schubert stratum $C(\hat\tau)$ in the single point $\tau$,
which is therefore an {\em isolated} point of this subset.

We know that $\ost(\tau)$ is an open subset of $\B(\tau,\hat\tau)$ with respect to the (Tits) metric.
\begin{lem}[Open stars]
\label{lem:opst}
$\ost(\tau)$ is open in $\B(\tau,\hat\tau)$ also with respect to the visual topology.
\end{lem}
\proof
Consider the fiber bundle (\ref{eq:ostfib}).
The union $U$ of the open $\taumod$-stars over the simplices in $C(\hat\tau)$ is open in $\geo X$.
Since $\tau$ is an isolated point of (\ref{eq:oppinsusp}),
the suspension $\B(\tau,\hat\tau)$ intersects $U$ 
precisely in $\ost(\tau)$, which is therefore open in the suspension. 
\qed

\subsubsection{Cones and parallel sets}
\label{sec:cones}

We transfer notions about stars by coning off. 
Our discussion takes place in $X$ and $\Fmod$.

Consider first the euclidean model chamber $\De=V(0,\simod)$.
Its {\em $\taumod$-boundary} $$\Dt\De:=V(0,\Dt\simod)\subseteq\D\De$$
is the union of the faces which do not contain the face $V(0,\taumod)$.
In particular $\Ds\De=\D\De$.

In the symmetric space $X$,
we define for a point $x\in X$ and a subset $A\subset\geo X$
the {\em cone} $V(x, A)\subset X$ as the union of the rays $x\xi$ for $\xi\in A$. 
We put $V(x, \emptyset):= \{x\}$.  

Let $\tau\subset\geo X$ be a simplex of type $\taumod$.
The {\em Weyl cone} $V(x, \st(\tau))$ with tip at $x\in X$ 
is the union of the euclidean Weyl chambers $V(x,\si)$ 
for all chambers $\si\subseteq\st(\tau)$,
equivalently, $\si\supseteq\tau$.
Its {\em boundary} is given by $\D V(x, \st(\tau))=V(x,\D\st(\tau))$,
and its {\em interior} by $V(x,\ost(\tau))-\{x\}$.
We call the Weyl sector $V(x,\tau)$ the {\em central sector} of the Weyl cone $V(x, \st(\tau))$. 
Similarly, we will refer to $V(0,\taumod)\subseteq\De$ as the {\em central sector} of the cone 
$\Wt \De=V(0,\Wt\simod)\subset\Fmod$. 

For the unique simplex $\hat\tau$ $x$-opposite to $\tau$,
the Weyl cone $V(x, \st(\tau))$ is contained in the parallel set $P(\tau,\hat\tau)$.
We say that the cone {\em spans} the parallel set.

Furthermore, for a compact subset $\Theta\subset\inte_{\taumod}(\simod)$,
we define the {\em $\Theta$-cone} $V(x, \stTh(\tau))$.

Note that for chambers $\si\subset\geo X$ we have $V(x, \st(\si))=V(x,\si)$.

\medskip
We will call two Weyl cones or $\Theta$-cones {\em asymptotic}
if their visual boundary stars coincide.

The {\em Hausdorff distance} of asymptotic Weyl cones $V(y,\st(\tau))$ and $V(y',\st(\tau))$ 
is finite and bounded by the distance $d(y,y')$ of their tips. This follows immediately 
from the corresponding fact for rays.

The distance between boundaries of Weyl cones will be discussed later
in section~\ref{sec:sepnest}.

\medskip
We will need a fact about {\em projections}.
Let 
\begin{equation}
\label{eq:prjwcocntsect}
\pi_{x,\tau}=\pi_{V(x,\tau)}:V(x,\st(\tau))\to V(x,\tau)
\end{equation}
denote the nearest point projection
of the Weyl cone to its central sector.

\begin{lem}
\label{lem:proj}
$\pi_{x,\tau}$ maps the interior of the Weyl cone to the interior of its central sector.
\end{lem}

In other words,
for every point $y$ in the interior of the Weyl cone there exists a point $p$ in the interior of its central sector
such that $py\perp V(x,\tau)$.

\proof
This is a consequence of the general Lemma~\ref{lem:sphprjfc} on projections of spherical simplices to their faces.
It yields at infinity that, for every chamber $\si\supseteq\tau$,
the nearest point projection $\inte_{\tau}(\si)\to\inte(\tau)$ is well-defined.
Equivalently, 
the nearest point projection $\ost(\tau)\to\inte(\tau)$ is well-defined.
The assertion follows by coning off.
\qed

\medskip
As a consequence of the lemma, 
$\pi_{x,\tau}$ agrees with the nearest point projection of the Weyl cone 
to the singular flat spanned by the sector $V(x,\tau)$,
because it does so on the interior.

\medskip
Now we address {\em convexity}. 
We will see that the results on stars carry over to cones.
First of all,
by the definition of Weyl convexity,
the cone $V(0,\Wt\Theta)=\Wt V(0,\Theta)\subset\Fmod$ is convex iff $\Theta$ is $\taumod$-Weyl convex.

\begin{prop}
[Convexity of cones]
\label{prop:thconeconv}
(i) The cones $V(x, \st(\tau))$ are convex.

(ii) If $\Theta$ is $\taumod$-Weyl convex, then also the cones $V(x,\stTh(\tau))$ are convex.
\end{prop}
\proof 
It suffices to verify (ii). 
We show that cones are intersections of horoballs. 

The horoball $\Hb_{\zeta,x}$ contains the cone $V(x,\stTh(\tau))$ iff $\stTh(\tau)\subseteq\ol B(\zeta,\pihalf)$ in $\geo X$.

Let $y\neq x$ be a point and let $x\xi$ be a ray extending $xy$.
Then $y\not\in V(x,\stTh(\tau))$ iff $\xi\not\in\stTh(\tau)$.
Let $F\subset X$ be a maximal flat 
such that $xy\subset F$ and $\tau\subset\geo F$.
According to the proof of Lemma~\ref{lem:thetastarconv}, 
there exists a point $\zeta\in\geo F$ such that 
$\ol B(\zeta,\pihalf)$ contains $\stTh(\tau)$ but not $\xi$.
Since $\Hb_{\zeta,x}\cap F$ is a half-space containing $x$ in its boundary,
it follows that also $y\not\in \Hb_{\zeta,x}$.
\qed

The convexity of cones implies their nestedness:
\begin{cor}
[Nestedness of cones] 
\label{cor:nestcone}
(i) If $y\in V(x,\st(\tau))$, then $V(y,\st(\tau))\subseteq V(x,\st(\tau))$. 

(ii) If $y\in V(x,\stTh(\tau))$, then $V(y,\stTh(\tau))\subseteq V(x,\stTh(\tau))$. 
\end{cor}

\medskip
Next we show an {\em openness} property for Weyl cones in the parallel sets spanned by them:
\begin{lem}[Open cones]
\label{lem:open-cone}
Let $x\in P(\tau,\hat\tau)$.
Then the boundary $\D V(x,\st(\tau))$ of the Weyl cone $V(x, \st(\tau))$
disconnects the parallel set,
and its interior $V(x,\ost(\tau))-\{x\}$ is one of the connected components.
\end{lem}
\proof 
Since parallel sets are cones over their visual boundaries,
$P(\tau,\hat\tau)=V(x,\geo X(\tau,\hat\tau))$,
this follows from the visual openness of stars, cf.\ Lemma~\ref{lem:opst}.
\qed

\subsubsection{Diamonds}
\label{sec:diamo}

We say that a nondegenerate oriented geodesic segment $xy\subset X$ is {\em $\taumod$-regular} 
if the unique geodesic ray $x\xi$ extending $xy$ is asymptotic to a $\taumod$-regular ideal point $\xi\in \geo X$.
In this case, we denote by $\tau(xy)\in\Flagt$ the unique simplex such that $\xi\in\ost(\tau)$. 
Furthermore,
we say that $xy$ is {\em $\Theta$-regular} with $\Theta\in\inte_{\taumod}(\simod)$ if $\theta(\xi)\in\Theta$. 

Note that $xy$ is $\taumod$-regular if and only if $yx$ is $\iota\taumod$-regular,
and $\Theta$-regular iff $yx$ is $\iota\Theta$-regular. 
The types of the simplices $\tau(xy)$ and $\tau(yx)\in\Flagit$ 
are then related by
$$
\theta(\tau(yx))= \iota \theta(\tau(xy)). 
$$
Let  $xy$ be a $\taumod$-regular segment. 
We define its {\em $\taumod$-diamond} as the intersection of Weyl cones
$$\diamot(x, y)=V(x,\st(\tau_+))\cap V(y,\st(\tau_-))\subset P(\tau_-,\tau_+)$$
where $\tau_+=\tau(xy)$ and $\tau_-=\tau(yx)$. 
The points $x, y$ are the {\em tips} of the diamond. 
Furthermore, if $xy$ is $\Theta$-regular, we define its {\em $\Theta$-diamond}
\begin{equation*}
\diamoTh(x,y)=
V(x,\stTh(\tau_+))\cap V(y,\stTh(\tau_-)) \subset \diamot(x, y).
\end{equation*}

The convexity of cones (Proposition~\ref{prop:thconeconv}) implies:
\begin{prop}[Convexity of diamonds]
(i) $\diamot(x,y)$ is convex.

(ii) If $\Theta$ is $\taumod$-Weyl convex, then also $\diamoTh(x,y)$ is convex. 
\end{prop}
And furthermore:
\begin{cor}
[Nestedness of diamonds]
\label{cor:nestdiamo}
Suppose that $xy$ and $x'y'$ are $\taumod$-regular segments such that 
$\tau(x' y')=\tau(xy)$,  $\tau(y' x')=\tau(yx)$ and $x'y'\subset\diamot(x, y)$.
Then:

(i) $\diamot(x', y')\subseteq\diamot(x, y)$. 

(ii) If $xy$ and $x'y'$ are $\Theta$-regular, where $\Theta$ is $\taumod$-Weyl convex,
and if $x'y'\subset\diamoTh(x, y)$, 
then 
$\diamoTh(x', y')\subseteq\diamoTh(x, y)$. 
\end{cor}

\subsection{Vector valued distances}\label{sec:VVD}

The Riemannian distance is not the complete two-point invariant on the symmetric space $X$,
if $\rank(X)\geq2$.
In view of the natural identifications 
$X\times X/G\cong X/K \cong \De$,
the full invariant is given by the quotient map
$$ d_{\De}:X\times X\to \De$$
arising from dividing out the $G$-action, 
which we refer to as the {\em $\De$-distance}. 
We will think of the elements of $\De\subset\Fmod$ as vectors and of $d_\De$ as a {\em vector-valued distance}. 
It relates to the Riemannian distance $d$ on $X$ by
$$
d= \|d_\Delta\|,
$$
where $\|\cdot\|$ is the euclidean norm on $\Fmod$. 

For the model flat, there are corresponding identifications 
$\Fmod\times\Fmod/W_{aff}\cong\Fmod/W\cong\De$
and a $\De$-distance 
$$ d_{\De}:\Fmod\times\Fmod\to \De.$$
It is compatible with the $\De$-distance on $X$
in that the charts $\Fmod\to X$ are $d_{\De}$-isometries. 

Similarly, one defines the $\De$-distance on {\em euclidean buildings} via apartment charts,
see \cite{ccm}.

The distance $d_\Delta$ is not symmetric, but satisfies 
$$
 d_{\De}(y,x)=\iota d_{\De}(x,y) .
$$ 
We refer the reader to \cite{ccm} and \cite{Parreau} for the detailed discussion of 
{\em metric properties} 
(such as ``triangle inequalities'' and ``nonpositive curvature behavior'') 
of $d_\Delta$. 

We note that a geodesic segment $xy\subset X$ is regular iff $d_\Delta(x,y)\in\inte(\De)$. 
Similarly, $xy$ is $\Theta$-regular iff $d_\Delta(x,y)\in V(0, \Theta)$.

We define certain coarsifications of $d_{\De}$ by composing it with linear maps:
For a face type $\taumod$, let 
$$ \pi^{\De}_{\taumod}:\De \to V(0,\taumod) $$
denote the nearest point projection.
The composition 
\begin{equation}
\label{eq:vvdt}
d_{\taumod} := \pi^{\De}_{\taumod} \circ d_{\De} 
\end{equation}
can also be regarded as a vector-valued distance on $X$, with values in the Weyl sector $V(0,\taumod)\subset\De$. 
Note that $d_{\simod}=d_{\De}$.
Obviously,
\begin{equation}
\label{eq:shrk}
\|d_{\taumod}\| \leq d
\end{equation}
because $\pi^{\De}_{\taumod}$ is 1-Lipschitz.

Given a compact subset $\Theta\subset \inte_{\taumod}(\simod)$, for $\Theta$-regular segments $xy\subset X$ it holds that
\begin{equation}
\label{eq:urgtdst}
\|d_{\taumod}(x,y)\| \geq \eps(\Theta)\cdot d(x,y)
\end{equation}
with a constant $\eps(\Theta)>0$, where $\|\cdot\|$ denotes the euclidean norm. 
For the constant $\eps(\Theta)$ one can take the sine of the angular distance $\angle(\Theta,\Dt\simod)$.

\subsection{Refined side lengths of triangles}\label{sec:refined}

In this section, 
we assume more generally that $X$ is a {\em CAT(0) model space},
i.e.\ a nonpositively curved Riemannian symmetric space or a thick euclidean building. 
We denote by 
$${\mathcal P}_3(X)\subset\De^3$$
the set of possible $\De$-side lengths 
$(d_{\De}(x_1,x_2),d_{\De}(x_2,x_3),d_{\De}(x_3,x_1))$
of triangles $\De(x_1,x_2,x_3)$ in $X$.
The following general result reduces the problem of determining ${\mathcal P}_3(X)$
from the symmetric space case to the euclidean building case:

\begin{thm}[{\cite[Thm.\ 1.2]{ccm}}]
\label{thm:slgthsbleqssp}
${\mathcal P}_3(X)$ depends only on the Weyl group $W$,
and not on whether $X$ is a Riemannian symmetric space or a thick euclidean building.
\end{thm}
In the paper \cite{ccm},
a detailed description of the set ${\mathcal P}_3(X)$ is given.

The next result concerns the $\De$-side lengths of triangles $\De(x,y,z)$ in $X$
such that the broken geodesic $xyz$ is a Finsler geodesic
(in the sense of section~\ref{sec:fins} below):
\begin{prop}
\label{prop:prjwcon}
(i) 
If $y\in V(x,\st(\tau))$ and $z\in V(y,\st(\tau))$ with $\tau\in\Flagt$, then
\begin{equation*}
d_{\De}(x,z) \in V(d_{\De}(x,y),\Wt\simod)\cap \De.
\end{equation*}
(ii) If $z\in V(y,\stTh(\tau))$, 
where $\Theta\subset\inte_{\taumod}(\simod)$ is 
$\taumod$-Weyl convex,
then 
\begin{equation*}
d_{\De}(x,z) \in V(d_{\De}(x,y),\Wt\Theta)\cap \De.
\end{equation*}
\end{prop}
Here, the cones $V(d_{\De}(x,y),\cdot)$ are to be understood as subsets of $\Fmod$.
\proof
We prove the stronger claim (ii).

The triangle $\De(x,y,z)$ lies in the parallel set $P=P(\hat\tau,\tau)$ 
for the simplex $\hat\tau\in\Flagit$ $x$-opposite to $\tau$.
The parallel set $P$ is itself a symmetric space (with euclidean factor)
with Weyl group $W'=\Wt\subset W$.
There is a natural inclusion $\simod\subset\simod'\subset\amod$ of spherical Weyl chambers 
such that $\simod'$ equals the convex hull of $\simod$ and the simplex $-\taumod$ opposite to $\taumod$,
and a corresponding  inclusion $\De\subset\De'\subset\Fmod$ of euclidean Weyl chambers 
such that $\De'$ is the convex hull of $\De$ 
and the sector $-V(0,\taumod)$.

Our claim 
is then a consequence of the following assertion on $\De'$-side lengths:
If $d_{\De'}(x,y)\in\De$ and $d_{\De'}(y,z)\in V(0,\Theta)\subset\De$, then 
\begin{equation*}
d_{\De'}(x,z)\in V(d_{\De'}(x,y),\Wt\Theta)\cap\De .
\end{equation*}
Using Theorem~\ref{thm:slgthsbleqssp},
we may pass from symmetric spaces to euclidean buildings:
The assertion is equivalent to the same assertion for any thick euclidean building $\tilde P$ with the same Weyl group $W'$.
(For instance, one can take $\tilde P$ to be the complete euclidean cone over the spherical building $\tits P$,
which is a non-locally compact euclidean building with just one vertex.)
It is easier to verify {the statement} in the building case due to the {\em locally conical} geometry of euclidean buildings.

Suppose therefore that $\De(\tilde x,\tilde y,\tilde z)$ is a triangle in a euclidean building $\tilde P$ with Weyl group $W'$,
satisfying the same assumptions $d_{\De'}(\tilde x,\tilde y)\in\De$ and $d_{\De'}(\tilde y,\tilde z)\in V(0,\Theta)$.
Taking advantage of the local conicality of buildings, we will do ``induction along $\tilde y\tilde z$''
and show that 
\begin{equation}
\label{eq:triineqfinspsind}
d_{\De'}(\tilde x,\tilde z')\in V(d_{\De'}(\tilde x,\tilde y),\Wt\Theta)\cap\De 
\end{equation}
for all points $\tilde z'\in \tilde y\tilde z$.
Since this is a closed condition on $\tilde z'$,
it suffices to show that the subset of points satisfying it is half-open to the right. 
Moreover, 
since the points $\tilde z'\in \tilde y\tilde z$ satisfying (\ref{eq:triineqfinspsind})
also satisfy, like $\tilde y$,  the assumptions 
that $d_{\De'}(\tilde x,\tilde z')\in\De$ and $d_{\De'}(\tilde z',\tilde z)\in V(0,\Theta)$,
it suffices to verify (\ref{eq:triineqfinspsind}) for all points $\tilde z'\in \tilde y\tilde z$ sufficiently close to $\tilde y$.

This however reduces our claim to the {\em flat} case,
because there exists a maximal flat $\tilde F\subset\tilde P$ 
which contains $\tilde x\tilde y$ 
along with a nondegenerate initial portion of the segment $\tilde y\tilde z$.\footnote{This is clear for discrete euclidean buildings.
(In particular, for buildings with only one vertex, like the complete euclidean cone over $\tits P$.)
For the general case, see e.g.\ \cite[\S 4.1.3]{qirigid}.}
We may therefore assume that the triangle $\De(\tilde x,\tilde y,\tilde z)$ lies entirely in $\tilde F$.
Identifying $\tilde F\cong\Fmod$,
we can once more reformulate our claim:
If $\de\in\De$ and $v\in V(0,\Wt\Theta)$, then 
\begin{equation}
\label{eq:triineqfinsflind}
d_{\De'}(0,\de+tv)\in V(\de,\Wt\Theta)\cap\De 
\end{equation}
for all sufficiently small $t\geq0$.

The stabilizer of $\de$ in $W'=\Wt$ 
is a subgroup $\Wn\leq\Wt$ for a face type $\numod$ with $\taumod\subseteq\numod\subseteq\simod$
(namely, for the minimal face type $\numod\supseteq\taumod$ such that $\de\in V(0,\numod)$).
We observe that the cone $\de+V(0,\Wt\Theta)$ is $\Wn$-invariant and can be represented 
{\em locally} near $\de$ as
$$\de+V(0,\Wt\Theta) = \Wn\bigl((\de+V(0,\Wt\Theta))\cap\De\bigr) .$$
The $\Wt$-invariance of $d_{\De'}(0,\cdot)$ yields the assertion. 
\qed

\subsection{Strong asymptote classes}
\label{sec:strongasy}

Let $\rho_1(t)$ and $\rho_2(t)$ be asymptotic geodesic rays in $X$, 
i.e.\ with the same ideal endpoint 
$\rho_1(+\infty)=\rho_2(+\infty)=\xi$.
Equivalently, 
the convex function 
$t\mapsto d(\rho_1(t),\rho_2(t))$ on $[0,+\infty)$ is bounded.
The rays are called {\em strongly asymptotic} if 
$d(\rho_1(t),\rho_2(t))\to 0$ as $t\to+\infty$. 
One sees then using Jacobi fields 
that $d(\rho_1(t),\rho_2(t))$ decays exponentially with rate depending on the type of $\xi$ (see \cite{Eberlein}). 

\medskip
Strong asymptote classes are represented by rays in a parallel set:
\begin{lem}
\label{lem:exstrasygeo}
Let $\xi, \hat\xi\in\geo X$ be antipodal. 
Then every geodesic ray asymptotic to $\xi$ is strongly asymptotic to a geodesic ray 
in the parallel set $P=P(\xi,\hat\xi)$.
\end{lem}
\proof 
Let $c_1(t)$ be a geodesic line forward asymptotic to $\xi$
(extending the given ray).
Then the function $t\mapsto d(c_1(t),P)$ 
is convex and bounded on $[0,+\infty)$, 
and hence non-increasing. We claim that the limit 
$$
D:=\lim_{t\to+\infty} d(c_1(t),P)
$$
equals zero. To see this, 
we choose a geodesic line $c_2(t)$ 
in $P$ forward asymptotic to $\xi$ 
and use the transvections $T^{c_2}_t$ along $c_2$ to ``pull back'' $c_1$: 
The geodesics $c_1^s:=T^{c_2}_{-s}c_1(\cdot+s)$ 
form a bounded family as $s\to+\infty$ 
and subconverge to a geodesic $c_1^{+\infty}$. 
Since the transvections $T^{c_2}_s$ preserve $P$, 
the distance functions 
$d(c_1^s(\cdot),P)=d(c_1(\cdot+s),P)$ 
converge locally uniformly on $\R$ and uniformly on $[0,+\infty)$ 
to the constant $D$. 
It follows that the limit geodesic $c_1^{+\infty}$ 
has distance $\equiv D$ from $P$. 
The same argument, applied to $c_2$ instead of the parallel set, 
implies that $c_1^{+\infty}$ is parallel to $c_2$. 
Thus, $c_1^{+\infty}\subset P(c_2)= P$ and, hence, 
$D=0$. 

Now we find a geodesic in $P$ strongly asymptotic to $c_1$ 
as follows. 
Let $t_n\to+\infty$. 
We choose geodesics $c'_n(t)$ in $P$ forward asymptotic to $\xi$ 
by requiring that $c'_n(t_n)\in P$ is the nearest point projection of $c_1(t_n)$.
Then $d(c_1(t_n),c'_n(t_n))=d(c_1(t_n),P)\to0$. 
The geodesics $c'_n\subset P$ are parallel,
and their mutual Hausdorff distances $d_{mn}$
are bounded above by the distances $d(c'_m(t),c'_n(t))$ independent of $t$.
To estimate the Hausdorff distances, we observe that 
$$ d_{mn}\leq d(c'_m(t),c'_n(t))\leq d(c'_m(t),c_1(t))+d(c_1(t),c'_n(t))\leq d(c'_m(t_m),c_1(t_m))+d(c_1(t_n),c'_n(t_n))$$
for $t\geq t_m,t_n$.
The right-hand side converges $\to0$ as $m,n\to+\infty$,
and hence also $d_{mn}$.
Thus, the geodesics $c'_n$ form a Cauchy sequence   
and therefore converge to a geodesic in $P$.
The limit geodesic is strongly asymptotic to $c_1$. 
\qed

\medskip
We now derive a {\em criterion} for the {\em strong asymptoticity of rays}.

Consider a geodesic line $c(t)$ asymptotic to $\xi\in\geo X$.
We observe that for every $\eta\in\geo P(c)$ the restriction $b_{\eta}\circ c$ is {\em linear},
because there exists a flat $f$ containing $c$ with $\eta\in\geo f$. 

As a consequence, 
for any two strongly asymptotic geodesic lines $c_1(t)$ and $c_2(t)$ 
asymptotic to $\xi$,
the restricted Busemann functions $b_{\eta}\circ c_i$ {\em coincide}
for every $\eta\in\st(\tau_{\xi})\subset\geo P(c_1)\cap\geo P(c_2)$,
where $\tau_{\xi}$ denotes the simplex spanned by $\xi$.

There is the following useful criterion for strong asymptoticity:

\begin{lem}
\label{lem:strasycrit}
For geodesic lines $c_1(t)$ and $c_2(t)$ asymptotic to $\xi$
the following are equivalent:

(i) $c_1(t)$ and $c_2(t)$ are strongly asymptotic. 

(ii) $b_{\eta}\circ c_1=b_{\eta}\circ c_2$ for every $\eta\in\st(\tau_{\xi})$. 

(ii') $b_{\eta}\circ c_1=b_{\eta}\circ c_2$ for every $\eta\in B(\xi,\eps)$ for some $\eps>0$.
\end{lem}
\proof (i)$\Ra$(ii) follows from the above discussion and (ii)$\Ra$(ii') is immediate. 

In order to prove (ii')$\Ra$(i), 
we replace the geodesics $c_i$ by a pair of parallel ones without changing their strong asymptote classes,
applying Lemma~\ref{lem:exstrasygeo}.
Using the implication (i)$\Ra$(ii), which we already proved,
we see that the $c_i$ keep satisfying hypothesis (ii').
Since they now lie in a common flat,
(ii') immediately implies that they coincide, i.e.\ (i) follows.
\qed

\medskip
We generalize the discussion of strong asymptoticity to {\em sectors}. 

Two Weyl sectors in $X$ are {\em asymptotic} iff their visual boundary simplices coincide,
equivalently, iff they have finite Hausdorff distance.

Fix a simplex $\tau\in\Flagt$ and consider two asymptotic sectors $V(x_1,\tau)$ and $V(x_2,\tau)$. 
The function $V(0,\taumod)\to[0,+\infty)$ given by 
\begin{equation}
\label{eq:distsect} 
y\mapsto d(\kappa_{x_1,\tau}(y),\kappa_{x_2,\tau}(y)) ,
\end{equation}
where $\kappa_{x_i,\tau}$ are the sector charts, 
is convex and bounded.
We denote its infimum by $d_{\tau}(x_1,x_2)$.
This defines a pseudo-metric $d_{\tau}$ on $X$, 
viewed as the set of (tips of) sectors asymptotic to $\tau$.\footnote{Observe that $d_{\tau}(x_1,x_2)$
depends only on the strong asymptote classes of the sectors $V(x_i,\tau)$,
and hence $d_{\tau}$ descends to $X_{\tau}^{par}\times X_{\tau}^{par}$.
The triangle inequality is a consequence of Proposition~\ref{prop:strasspariso} below.
One can also verify the triangle inequality for $d_{\tau}$ directly,
using the fact that,
for bounded convex functions $\phi,\psi:V(0,\taumod)\to[0,+\infty)$,
it holds that $\inf\phi+\inf\psi=\inf(\phi+\psi)$.}

We say that the sectors $V(x_1,\tau)$ and $V(x_2,\tau)$ are {\em strongly asymptotic} 
if $d_{\tau}(x_1,x_2)=0$. 
For any ideal point $\xi\in \inte(\tau)$ 
this is equivalent to the rays $x_1\xi$ and $x_2\xi$ being strongly asymptotic. 
We denote by 
\begin{equation*}
X_{\tau}^{par}=X/\sim_{d_{\tau}}
\end{equation*}
the {\em space of strong asymptote classes}
of Weyl sectors asymptotic to $\tau$. 

We show now that, also in the case of sectors, parallel sets represent strong asymptote classes.
For a simplex $\hat\tau$ opposite to $\tau$ 
we consider the restriction 
\begin{equation}
\label{eq:parasy}
P(\tau,\hat\tau)\to X_{\tau}^{par}
\end{equation}
of the natural projection $X\to X_{\tau}^{par}$. 

\begin{prop}
\label{prop:strasspariso}
The map (\ref{eq:parasy}) is an isometry. 
\end{prop}
\proof
For points $x_1,x_2\in P(\tau,\hat\tau)$
the function (\ref{eq:distsect}) is constant
$\equiv d(x_1,x_2)$. 
Hence (\ref{eq:parasy}) is an isometric embedding. 
To see that it is also surjective, 
we need to verify that 
every sector $V(x,\tau)$ is strongly asymptotic to a sector 
$V(x',\tau)\subset P(\tau,\hat\tau)$. 
This follows from the corresponding fact for geodesic rays,
see Lemma~\ref{lem:exstrasygeo}.
\qed

\subsection{Asymptotic Weyl cones}\label{sec:ascones-shadows}

\subsubsection{Separation of nested Weyl cones}
\label{sec:sepnest}

Suppose that $y\in V(x,\st(\tau))$ with $\tau\in\Flagt$.
By nestedness (Coroillary \ref{cor:nestcone}), 
we have the inclusion of Weyl cones $V(y,\st(\tau))\subseteq V(x,\st(\tau))$.
We now determine the separation of their boundaries:

\begin{prop}[Separation]
\label{prop:hdstwcones}
The nearest point distance of the boundaries $\D V(x,\st(\tau))$ and $\D V(y,\st(\tau))$ 
equals 
$d(\de,\Dt\De)=d(y,\D V(x,\st(\tau)))$,
where $\de=d_{\De}(x,y)$.
\end{prop}
\proof
The natural submersion 
\begin{equation*}
d_{\De}(x,\cdot) : X\to\De
\end{equation*}
is 1-Lipschitz and 
restricts to an isometry on every euclidean Weyl chamber with tip at $x$.
By restricting it to the Weyl cone $V(x,\st(\tau))$,
one sees that 
\begin{equation*}
d(\cdot,\D V(x,\st(\tau))) = d\bigl(d_{\De}(x,\cdot) , \Dt\De \bigr) 
\end{equation*}
on $V(x,\st(\tau))$.
According to Proposition~\ref{prop:prjwcon}(i),
the values of $d_{\De}(x,\cdot)$ on $V(y,\st(\tau))$
are contained in 
\begin{equation*}
V(\de,\Wt\simod) \cap\De ,
\end{equation*}
and clearly all these values are attained
(on a euclidean Weyl chamber with tip at $x$ and containing $y$).
It follows that the nearest point distance of $V(y,\st(\tau))$ and $\D V(x,\st(\tau))$
equals the nearest point distance of 
$V(\de,\Wt\simod) \cap\De$ and $\Dt\De$.

In order to see that 
the latter is given by $d(\de,\Dt\De)$,
note that $d(\cdot,\Dt\De)$ is the minimum of finitely many root functionals on $\De$,
namely of those corresponding to the walls of $\De$ not containing the sector $V(0,\taumod)$,
equivalently, of those which are nonnegative on $\Wt\De$.
Each of these functionals attains its minimum on the cone $V(\de,\Wt\simod)$ at its tip $\de$.
\qed

\subsubsection{Shadows at infinity and strong asymptoticity of Weyl cones}
\label{sec:shadw}

For a simplex $\tau_-\in\Flagit$ and a point $x\in X$,
we consider the function 
\begin{equation}
\label{eq:distfrpar}
\tau\mapsto d(x,P(\tau_-,\tau))
\end{equation}
on the open Schubert stratum $C(\tau_-)\subset\Flagt$.
We denote by $\tau_+\in C(\tau_-)$ the simplex $x$-opposite to $\tau_-$.
\begin{lem}
\label{lem:contpr}
The function (\ref{eq:distfrpar}) is continuous and proper.
\end{lem}
\proof
This follows from the fact that 
$C(\tau_-)$ and $X$ are homogeneous spaces for the parabolic subgroup $P_{\tau_-}$.
Indeed, 
continuity follows from the continuity of the function 
$$g\mapsto d(x,P(\tau_-,g\tau_+))=d(g^{-1}x,P(\tau_-,\tau_+))$$
on $P_{\tau_-}$ which factors through the orbit map $P_{\tau_-}\to C(\tau_-),g\mapsto g\tau_+$.

Regarding properness,
note that a simplex $\tau\in C(\tau_-)$ is determined by any point $y$ contained in the parallel set $P(\tau_-,\tau)$,
namely as the simplex $y$-opposite to $\tau_-$.
Thus,
if $P(\tau_-,\tau)\cap B(x,R)\neq\emptyset$ for some fixed $R>0$,
then there exists $g\in P_{\tau_-}$ such that $\tau=g\tau_+$ and $d(x,gx)<R$. 
In particular, $g$ lies in a compact subset.
This implies properness.
\qed

\medskip
Moreover, 
the function (\ref{eq:distfrpar}) has a unique minimum zero in $\tau_+$. 

We define the following open subsets of $C(\tau_-)$ 
which can be regarded as {\em shadows} of balls in $X$ with respect to $\tau_-$. 
For $x\in X$ and $r>0$, we put
\begin{equation}
\label{eq:shadwdf}
U_{\tau_-,x,r}:=\{\tau\in C(\tau_-) | d(x,P(\tau_-,\tau))<r\} .
\end{equation}
The next fact expresses the {\em strong asymptoticity}
of asymptotic Weyl cones: 
\begin{lem}
\label{lem:expconvsect}
For $r,R>0$ there exists $d=d(r,R)>0$ such that:

If $y\in V(x,\st(\tau_-))$ with $d(y,\D V(x,\st(\tau_-)))\geq d(r,R)$, 
then 
$U_{\tau_-,x,R}\subset U_{\tau_-,y,r}$.
\end{lem}
\proof
If $U_{\tau_-,x,R}\not\subset U_{\tau_-,y,r}$ 
then there exists $x'\in B(x,R)$ such that $d(y,V(x',\st(\tau_-)))\geq r$. 
Thus, if the assertion is wrong, 
there exist a sequence $x_n\to x_{\infty}$ in $\ol B(x,R)$ 
and an $\iota\taumod$-regular sequence 
$(y_n)$ in $V(x,\st(\tau_-))$ such that 
$d(y_n,V(x_n,\st(\tau_-)))\geq r$.  

Let $\rho:[0,+\infty)\to V(x,\tau_-)$ be a geodesic ray with initial point $x$ and asymptotic to an interior point of $\tau_-$.
By $\iota\taumod$-regularity,
the sequence $(y_n)$ eventually enters every Weyl cone $V(\rho(t),\st(\tau_-))$. 
Since the distance function $d(\cdot,V(x_n,\st(\tau_-)))$ is convex and bounded, and hence non-increasing 
along rays asymptotic to $\st(\tau_-)$, 
we have that 
\begin{equation*}
R\geq d(x,V(x_n,\st(\tau_-))) 
\geq d(\rho(t),V(x_n,\st(\tau_-)))\geq d(y_n,V(x_n,\st(\tau_-)))\geq r
\end{equation*}
for $n\geq n(t)$. 
It follows that 
\begin{equation*}
R\geq d(\rho(t),V(x_{\infty},\st(\tau_-)))\geq r
\end{equation*}
for all $t\geq0$. 
However, the ray $\rho$ is strongly asymptotic to $V(x_{\infty},\st(\tau_-))$, 
cf.\ Proposition~\ref{prop:strasspariso},
a contradiction. 
\qed

\subsection{Horocycles}
\label{sec:horocycles}

We discuss various foliations of $X$ naturally associated to a simplex $\tau\subset\geo X$.

We begin with foliations by flats and parallel sets:
First, 
we denote by ${\mathcal F}_{\tau}$ the partition of $X$ 
into the singular flats $f\subset X$ such that $\tau\subset\geo f$ is a top-dimensional simplex.
Second, 
we consider the partition ${\mathcal P}_{\tau}$ of $X$ 
into the parallel sets $P(\tau,\hat\tau)$ for the simplices $\hat\tau$ opposite to $\tau$.
Note that ${\mathcal P}_{\tau}$ is a coarsening of ${\mathcal F}_{\tau}$,
and coincides with it iff $\tau$ is a chamber.
The parabolic subgroup $P_{\tau}$ preserves both partitions 
and acts transitively on their leaves.
This implies that these partitions are {\em smooth foliations}. 

We will now show that there exist complementary orthogonal foliations.
To do so, we describe preferred mutual identifications between the leaves of ${\mathcal F}_{\tau}$ as well as of ${\mathcal P}_{\tau}$
by the actions of certain subgroups of $P_{\tau}$. 
Their orbits will be submanifolds orthogonal and complementary to the foliations,
i.e.\ the integral submanifolds of the distributions normal to them.

The tuple $(b_{\xi})_{\xi\in\Vert(\tau)}$ 
of Busemann functions for the vertices $\xi$ of $\tau$ 
(well-defined up to additive constants) 
provides affine coordinates simultaneously for each flat $f\in{\mathcal F}_{\tau}$. 
The Busemann functions at the other ideal points in $\tau$ are linear combinations of these.
The group $P_{\tau}$ preserves the family of horospheres at every $\xi\in\tau$, 
and the action on it yields a natural ``shift'' homomorphism $\phi_{\xi}:P_{\tau}\to\R$.
The intersection of their kernels forms the normal subgroup 
\begin{equation}
\label{eq:largehorocgp}
\bigcap_{\xi\in\Vert(\tau)}\Stab(b_{\xi})
=\bigcap_{\xi\in\tau}\Stab(b_{\xi})
\triangleleft P_{\tau} .
\end{equation}
It acts transitively on the set ${\mathcal F}_{\tau}$ of flats and preserves the coordinates; 
it thus provides consistent identifications between these flats.
The level sets of $(b_{\xi})_{\xi\in\Vert(\tau)}$ 
are submanifolds orthogonal and complementary to these flats,
because the gradient directions of the Busemann functions $b_{\xi}$
at a point $x\in f\in{\mathcal F}_{\tau}$ 
constitute a basis of the tangent space $T_xf$. 
These level sets form a smooth foliation ${\mathcal F}_{\tau}^{\perp}$ 
and are the orbits of the subgroup (\ref{eq:largehorocgp}). 

In order to describe the foliation normal to ${\mathcal P}_{\tau}$,
we define the {\em horocyclic subgroup} at $\tau$ 
as the (smaller) normal subgroup $N_\tau \triangleleft P_{\tau}$ given by
\begin{equation*}
N_{\tau}=\bigcap_{\xi\in\st(\tau)}\Stab(b_{\xi}) \triangleleft \Fix(\st(\tau))\triangleleft P_{\tau} .
\end{equation*}
It is the kernel of the $P_\tau$-action 
on the set of all (unnormalized) Busemann functions centered at ideal points in $\st(\tau)$.

Note that 
as a consequence of Lemma~\ref{lem:strasycrit},
$N_{\tau}$ preserves the strong asymptote classes of geodesic rays at all ideal points $\xi\in\ost(\tau)$.

\medskip
We now give a method for constructing isometries in $N_{\tau}$.

Let 
$\xi\in\inte(\tau)$, 
and let $c(t)$ be a geodesic line forward asymptotic to it,
$c(+\infty)=\xi$.
Consider the one parameter group
$(T_t^c)_{t\in \R}$ 
of transvections along $c$. 
The transvections $T_t^c$ fix $\geo P(c)$ pointwise 
and shift the Busemann functions $b_{\eta}$ centered at ideal points $\eta\in\geo P(c)$
by additive constants:
\begin{equation*}
b_{\eta}\circ T_t^c-b_{\eta} \equiv -t\cdot\cos\tangle(\eta,\xi)
\end{equation*}
Note that $\st(\tau)\subset\geo P(c)$.

\begin{lem}
\label{lem:constrhiso}
Let $c_1(t)$ and $c_2(t)$ be geodesic lines forward asymptotic to $\xi\in\inte(\tau)$,
which are strongly asymptotic.
Then there exists an isometry\footnote{This isometry is unipotent but we will not need this fact.}  
$n\in G$ with the properties:

(i) $n\circ c_1=c_2$.

(ii) $n$ fixes $\geo P(c_1)\cap\geo P(c_2)$ pointwise.

(iii) $b_{\eta}\circ n\equiv b_{\eta}$ for all $\eta\in\geo P(c_1)\cap\geo P(c_2)$. 

\no
In particular, $n\in N_{\tau}$.
\end{lem}
\proof
By our observation above, 
the isometries $T_{-t}^{c_2}\circ T_t^{c_1}$ 
fix $\geo P(c_1)\cap\geo P(c_2)\supseteq\st(\tau)$ pointwise
and preserve the Busemann functions $b_{\eta}$ for all $\eta\in\geo P(c_1)\cap\geo P(c_2)$.
Thus, they belong to $N_{\tau}$.
Moreover, 
they form a bounded family.
Therefore,
as $t\to+\infty$, 
they subconverge to an isometry $n\in N_{\tau}$ 
which maps $c_1$ to $c_2$ while preserving parameterizations.
\qed

\begin{cor}
\label{cor:horotranspar}
$N_{\tau}$ acts transitively on 

(i) every strong asymptote class of geodesic rays at every ideal point $\xi\in\interior(\tau)$;

(ii) the set of leaves of ${\mathcal P}_{\tau}$.
\end{cor}
\proof
Part (i) is a direct consequence of the lemma. 

Also (ii) follows because
every parallel set in ${\mathcal P}_{\tau}$ contains a (in fact, exactly one) geodesic ray of every strong asymptote class 
at any point $\xi\in\inte(\tau)$, 
cf.\ Proposition~\ref{prop:strasspariso}.
\qed

\begin{rem}
One also obtains 
that every geodesic asymptotic to an ideal point $\xi\in\D\tau$
can be carried by an isometry in $N_{\tau}$ to any other strongly asymptotic geodesic.
However, 
$N_{\tau}$ does not preserve strong asymptote classes at $\xi$ in that case. 
\end{rem}

\begin{lem}
If $n\in N_{\tau}$ 
preserves a parallel set $P(\tau,\hat\tau)$, 
$n\hat\tau=\hat\tau$, 
then it acts trivially on it. 
\end{lem}
\proof 
The hypothesis implies that 
$n$ fixes $\st(\tau)$ and $\hat\tau$ pointwise, 
and hence also their convex hull $\geo P(\tau,\hat\tau)$ in $\tits X$. 
Thus $n$ preserves every maximal flat $F\subset P(\tau,\hat\tau)$. 
Moreover it preserves all Busemann functions $b_{\xi}$ 
centered at points $\xi\in\geo F\cap\st(\tau)$, 
and therefore must fix $F$ pointwise,
compare Lemma~\ref{lem:strasycrit}.
\qed

\begin{cor}
The stabilizer of $P(\tau,\hat\tau)$ in $N_{\tau}$ is its pointwise fixator $K_{\tau,\hat\tau}<G$.
\end{cor}
\proof 
The claim follows from the obvious inclusion 
$K_{\tau,\hat\tau}\subset N_{\tau}$ together with the lemma.
\qed

\begin{rem}
The subgroup $N_\tau$ decomposes as the semidirect product $U_\tau\rtimes K_{\tau,\hat\tau}$, 
where $U_\tau \triangleleft P_\tau$ is the {\em unipotent radical} of $P_\tau$.
\end{rem}

By the above, $N_{\tau}$ provides consistent 
identifications between 
the parallel sets $P(\tau,\hat\tau)$. 
The $N_{\tau}$-orbits are submanifolds orthogonal to the parallel sets
and must have complementary dimension. 
They form a smooth foliation 
\begin{equation}
\label{eq:horfol}
{\mathcal H}_{\tau}={\mathcal P}_{\tau}^{\perp}
\end{equation}
refining ${\mathcal F}_{\tau}^{\perp}$, 
which we call the {\em horocyclic foliation} 
and its leaves the {\em horocycles} at $\tau$. 
We denote the horocycle at $\tau$ through the point $x$ by $\Hc_{\tau, x}$, i.e.\ $\Hc_{\tau,x}=N_{\tau}x$.

For incident faces, 
the associated subgroups and foliations 
are contained in each other:
If $\ups\subset\tau$, 
then $\st(\ups)\supset\st(\tau)$ and $N_{\ups}< N_{\tau}$. 
Therefore, e.g.\ 
${\mathcal H}_{\ups}$ refines ${\mathcal H}_{\tau}$.

Note that in rank one, horocycles are horospheres.

\medskip
We also see how horocycles and strong asymptote classes relate;
by Corollary~\ref{cor:horotranspar}(i):
\begin{cor}[Strong asymptote classes are horocycles]
The sectors $V(x_1,\tau)$ and $V(x_2,\tau)$ are strongly asymptotic 
if and only if $x_1$ and $x_2$ lie in the same horocycle at $\tau$. 
\end{cor}

Moreover,
the discussion shows that 
for the stabilizer $P_{\tau}\cap P_{\hat\tau}$ of $P(\tau,\hat\tau)$ in $P_{\tau}$
it holds that 
$N_{\tau}(P_{\tau}\cap P_{\hat\tau})=P_{\tau}$
and 
$P_{\tau}\cap P_{\hat\tau}\cap N_{\tau}=K_{\tau,\hat\tau}$, 
and so the sequence 
\begin{equation*}
1\to N_{\tau}\to P_{\tau}\to\Isom(X_{\tau}^{par}) 
\end{equation*}
is exact.

\begin{rem}
Note that the homomorphism $P_{\tau}\to\Isom(X_{\tau}^{par})$ is in general not surjective.
Namely, let $X_{\tau}^{par}=:f_{\tau}\times CS(\tau)$ 
denote the decomposition (\ref{eq:parsplit}) of $X_{\tau}^{par}\cong P(\tau,\hat\tau)$.
Then $P_{\tau}$ acts on the flat factor $f_{\tau}$ only by the group $A_{\tau}$ of translations. 
On the cross section, it acts by a subgroup $M_\tau\leq\Isom(CS(\tau))$ containing the identity component.
The above exact sequence is then a part of the  {\em Langlands' decomposition} of $P_\tau$, 
$$
1\to N_\tau \to P_\tau \to A_\tau\times M_\tau \to 1, 
$$
which, on the level of Lie algebras, is a split exact sequence.  
\end{rem}

We return now to Lemma~\ref{lem:constrhiso}.
For later use, 
we elaborate on the special case 
when the geodesics $c_i$ are contained in the parallel set of a 
singular flat of dimension rank minus one.

Consider a half-apartment  $h\subset\geo X$; 
it is a simplicial $\pihalf$-ball in $\geo X$.
We call its center $\zeta$ the {\em pole} of $h$. 
We define the {\em star} $\st(h)$ as the union of the stars $\st(\tau)$
where $\tau$ runs through all simplices with $\inte(\tau)\subset\inte(h)$,
equivalently, which are spanned by interior points of $h$.
Similarly, we define the {\em open star} $\ost(h)$ as the union of the corresponding open stars $\ost(\tau)$.
Note that $\inte(h)\subset\ost(h)$.
Furthermore, 
we define the subgroup $N_h< G$ as the intersection of the horocyclic subgroups $N_{\tau}$
at these simplices $\tau$,
$$ N_h = \bigcap_{\inte(\tau)\subset\inte(h)} N_{\tau}.$$

We observe that $N_h$ 
preserves the strong asymptote classes of geodesic rays at all ideal points $\xi\in\ost(h)$,
and it preserves the family of maximal flats $F$ with $\geo F\supset h$.
The action on this set of flats is transitive.
Indeed, parallel to Lemma~\ref{lem:constrhiso}, we have:

\begin{lem}
\label{lem:constrhisospec}
Let $F_1,F_2\subset P(\D h)$ be maximal flats with $\geo F_i\supset h$.
Then there exists an isometry $n\in N_h$ with the properties:

(i) $nF_1=F_2$.

(ii) $n$ fixes $\st(h)$ pointwise.

(iii) $b_{\eta}\circ n\equiv b_{\eta}$ for all $\eta\in\st(h)$. 
\end{lem}
\proof The parallel set  $P(\D h)$ splits as the product $f\times  CS(\D h)$,
see \eqref{eq:parsplit}, 
where $f\subset X$ is a singular flat with $\geo f=\D h$,
and the cross section $CS(\D h)$ is a rank one symmetric space. 
Accordingly, 
the maximal flats $F_i$ split as products $f\times \bar c_i$ with geodesics $\bar c_i\subset CS(\D h)$
asymptotic to the pole $\zeta\in CS(\D h)$ of $h$.

Let $\xi\in\inte(h)$.
We choose geodesics $c_1(t),c_2(t)$ in $F_1,F_2$ asymptotic to $\xi$.
Their $f$-com\-po\-nents are parallel geodesics in $f$,
and their $CS(\D h)$-components are geodesics in $CS(\D h)$ asymptotic to $\zeta$, 
equal to $\bar{c}_1, \bar{c}_2$ up to reparametrization. 
The geodesics $c_1,c_2$ are strongly asymptotic 
iff they have the same $f$-component 
and their $CS(\D h)$-components are strongly asymptotic. 
We choose them in this way,
using the fact that any two asymptotic geodesics in a rank one symmetric space 
become strongly asymptotic after suitable reparameterization. 

We then can apply the limiting argument (in the proof of Lemma~\ref{lem:constrhiso}) 
to the compositions 
$T_{-t}^{c_2}\circ T_t^{c_1}$
and obtain an isometry $n\in N_{\tau_{\xi}}$
where $\tau_{\xi}\subset h$ denotes the simplex spanned by $\xi$.
The isometry $n$ carries $F_1$ to $F_2$,
fixes $\st(\tau_{\xi})$ pointwise 
and satisfies (iii) for all $\eta\in\st(\tau_{\xi})$.

We observe that the isometries $T_{-t}^{c_2}\circ T_t^{c_1}$ act trivially on $f$
and the limiting isometry $n$ depends only on the $CS(\D h)$-components of the geodesics $c_i$.
Thus, by replacing the $f$-component of the $c_i$,
we are not affecting $n$, but we can change the ideal endpoint $\xi$ of the $c_i$ 
to any other ideal point $\xi'\in\inte(h)$.
(We work here with constant speed parametrizations $c_i(t)$.)
It follows that $n$ fixes also $\st(\tau_{\xi'})$ pointwise 
and satisfies (iii) also for all $\eta\in\st(\tau_{\xi'})$.
Varying $\xi'$, 
we let $\tau_{\xi'}$ run through all simplices with $\interior(\tau)\subset\interior(h)$ 
and conclude also parts (ii)+(iii) of the assertion. 
\qed

\medskip
We obtain an analogue of Corollary~\ref{cor:horotranspar}:

\begin{cor}
\label{cor:horotransparspec}
$N_h$ acts transitively on 

(i) every strong asymptote class of geodesic rays at every ideal point $\xi\in\interior(h)$;

(ii) the set of maximal flats $F$ with $\geo F\supset h$.
\end{cor}

We describe a consequence of our discussion for the horocyclic foliations. 

The maximal flats $F$ with $\geo F\supset h$
are contained in the parallel set $P(\D h)\cong f\times  CS(\D h)$ 
and form the leaves of a smooth foliation ${\mathcal P}_h$ of $P(\D h)$.
This foliation is the pullback (via the natural projection $P(\D h)\to CS(\D h)$) 
of the one-dimensional foliation of the rank one symmetric space $CS(\D h)$ by the geodesics 
asymptotic to the ideal point $\zeta\in\geo CS(\D h)$, the center of $h$. 
There exists a foliation 
${\mathcal H}_h$ of $P(\D h)$ whose leaves are normal (orthogonal and complementary) 
to those of ${\mathcal P}_h$. 
The leaves of ${\mathcal H}_h$ have the form $\{y\}\times \Hs_{\zeta,z}$, where $y\in f$ and 
$\Hs_{\zeta,z}\subset CS(\D h)$ is the horosphere centered at $\zeta$ and passing through $z\in CS(\D h)$. 
We call the leaves of ${\mathcal H}_h$ the {\em horocycles at $h$} 
and the foliation ${\mathcal H}_h$ the {\em horocyclic} foliation. 
The leaf of ${\mathcal H}_h$ passing through $x\in P(\D h)$ will be denoted $\Hc_{h,x}$. 
Corollary~\ref{cor:horotransparspec} implies that $\Hc_{h,x}=N_hx$.

Let $\tau$ be a simplex so that $\interior(\tau)\subset\interior(h)$.
Then the foliation ${\mathcal P}_{\tau}$ of $X$ by parallel sets
restricts on $P(\D h)$ to the foliation ${\mathcal P}_h$ by maximal flats,
and the horocyclic foliation ${\mathcal H}_{\tau}$
restricts to the horocyclic foliation ${\mathcal H}_h$.
(This follows from the fact that the foliations ${\mathcal P}_{\tau}$ and ${\mathcal H}_{\tau}$
are normal to each other, cf.\ (\ref{eq:horfol}).)
In other words, the horocyclic foliations ${\mathcal H}_{\tau}$ 
for the various simplices $\tau$ with $\interior(\tau)\subset \interior(h)$ 
{\em coincide} on the parallel set $P(\D h)$.

\subsection{Contraction at infinity}
\label{sec:contraction}

\subsubsection{Identifications of horocycles}

We fix a simplex $\tau\subset\geo X$.
Since every horocycle at $\tau$ intersects every parallel set $P(\tau,\hat\tau)$, $\hat\tau\in C(\tau)$, exactly once, 
there are $N_{\tau}$-equivariant diffeomorphisms
\begin{equation}
\label{eq:hc-schubert}
\Hc_{\tau,x} \buildrel\cong\over\to C(\tau)
\end{equation}
sending a point $y\in \Hc_{\tau,x}$ to the unique simplex $\hat\tau\in C(\tau)$ such that 
$\Hc_{\tau,x} \cap P(\tau, \hat\tau)=\{y\}$. 
(The smoothness of these identifications follows from their $N_{\tau}$-equivariance.)
Composing the maps \eqref{eq:hc-schubert} 
and their inverses, we obtain $N_\tau$-equivariant diffeomorphisms
\begin{equation}
\label{eq:horo-diffeo}
\pi^\tau_{x'x}: \Hc_{\tau,x}\to \Hc_{\tau,x'},
\end{equation}
sending the intersection point $\Hc_{\tau,x} \cap P(\tau, \hat\tau)$ to the intersection $\Hc_{\tau,x'} \cap P(\tau, \hat\tau)$
for $\hat\tau\in C(\tau)$.

Let $h\subset\geo X$ be a half-apartment such that $\interior(\tau)\subset\interior(h)$. 
Then, as discussed in the end of the previous section, the horocycles at $\tau$ intersect the parallel set 
$P(\D h)$ in the horocycles at $h$.
The latter are homogeneous spaces for the subgroup $N_h<N_{\tau}$.
Thus, 
for $x,x'\in P(\D h)$, 
the diffeomorphisms \eqref{eq:horo-diffeo} 
restrict to $N_h$-equivariant diffeomorphisms 
\begin{equation*}
\pi^h_{x'x}: \Hc_{h,x}\buildrel\cong\over\to \Hc_{h,x'}
\end{equation*}
between the horocycles at $h$, while the diffeomorphisms \eqref{eq:hc-schubert} restrict to 
$N_h$-equivariant diffeomorphisms  
\begin{equation*}
\Hc_{h,x} \buildrel\cong\over\to C(h)
\end{equation*}
between the horocycles at $h$ and the $N_h$-orbit 
$C(h)\subset C(\tau)$ consisting of the simplices which are contained in $\geo P(\D h)$.  

\medskip
We estimate now the {\em contraction-expansion} of the identifications $\pi^h_{x'x}$. 

We build on the discussion at the end of the previous section. 
As we saw, the horocycles $\Hc_{h,x}$ in $P(\D h)\cong f\times CS(\D h)$ 
are horospheres in the cross sections $pt\times CS(\D h)$.
They therefore project isometrically 
onto the horospheres $\Hs_{\zeta,\bar x}$ in $CS(\D h)$, 
where $\bar x$ denotes the projection of $x$. 
Under these projections,
the identifications $\pi^h_{x'x}$ 
correspond to the 
identifications 
\begin{equation}
\label{eq:hrsphidfcrsc}
\pi^{\zeta}_{\bar x'\bar x}: \Hs_{\zeta,\bar x}\buildrel\cong\over\to\Hs_{\zeta,\bar x'}
\end{equation}
of horospheres,
i.e.\ for $x,x'\in P(\D h)$,
we have the commutative diagram:
$$
\begin{array}{ccc}
\Hc_{h,x} & \buildrel{\pi^h_{x'x}}\over{\longrightarrow} &			 \Hc_{h,x'}\\
\downarrow & ~                    &			\downarrow\\
\Hs_{\zeta,\bar x} & \buildrel{\pi^{\zeta}_{\bar x'\bar x}}\over{\longrightarrow} & \Hs_{\zeta,\bar x'}
\end{array}
$$
Estimating the contraction rate of $\pi^h_{x'x}$ therefore reduces to estimating it for $\pi^{\zeta}_{\bar x'\bar x}$
in the rank one symmetric space $CS(\D h)$. 

We estimate the {\em infinitesimal} contraction.
We assume that $\bar x'$ is closer to $\zeta$ than $\bar x$,
$b_{\zeta}(\bar x)\geq b_{\zeta}(\bar x')$.
Then there is actual contraction, 
at a uniform rate in terms of the distance between the horospheres. 
For the differential 
$d\pi^{\zeta}_{\bar x'\bar x}$ of $\pi^{\zeta}_{\bar x'\bar x}$,
one has the estimate
\begin{equation*}
e^{-c_1(b_{\zeta}(\bar x)-b_{\zeta}(\bar x'))} \|\bar v\| \leq \|(d\pi^{\zeta}_{\bar x'\bar x}) \bar v\|   
\leq e^{-c_2(b_{\zeta}(\bar x)-b_{\zeta}(\bar x'))} \|\bar v\|
\end{equation*}
for all tangent vectors $\bar v\in T\Hs_{\zeta,\bar x}$, 
with constants $c_1\geq c_2>0$ 
depending only on the rank one symmetric space $CS(\D h)$, 
in fact, depending only on $X$,
because there are only finitely many isometry types of rank one symmetric spaces 
occurring as cross sections of parallel sets in $X$. 
The estimate follows from the standard fact that the exponential decay rate of decaying 
Jacobi fields along geodesic rays in $CS(\D h)$ 
is bounded below and above 
(in terms of the eigenvalues of the curvature tensor).

In view of 
$b_{\zeta}(x)-b_{\zeta}(x')=b_{\zeta}(\bar x)-b_{\zeta}(\bar x')$,
we obtain for $\pi^h_{x'x}$:
\begin{lem}[Infinitesimal contraction of horocycle identifications]
\label{lem:contrac}
If $b_{\zeta}(x)\geq b_{\zeta}(x')$, then 
\begin{equation}
\label{ineq:contrestrk1}
e^{-c_1(b_{\zeta}(x)-b_{\zeta}(x'))} \|v\| \leq \|(d\pi^h_{x'x}) v\|   \leq e^{-c_2(b_{\zeta}(x)-b_{\zeta}(x'))} \|v\|   
\end{equation}
for all tangent vectors $v$ to $\Hc_{h,x}$, 
with constants $c_1,c_2>0$ depending only on $X$. 
\end{lem}

\subsubsection{Infinitesimal contraction of transvections}

We now focus on transvections and their action at infinity.

Suppose that $x,x'\in P(\tau,\hat\tau)$ are distinct points. 
Let $\vartheta_{xx'}$ denote the transvection with axis $l=l_{xx'}$ through $x$ and $x'$ 
mapping $x'\mapsto x$;
we orient the geodesic $l_{xx'}$ from $x'$ to $x$,
i.e.\ so that $\vartheta_{xx'}$ translates along it in the positive direction. 
The transvection $\vartheta_{xx'}$  preserves the parallel set $P(\tau,\hat\tau)$ and 
fixes the simplices $\tau,\hat\tau$ at infinity.

We consider the action of $\vartheta_{xx'}$ on $C(\tau)$ 
and its differential at the fixed point $\hat\tau$. 
Modulo the identifications 
\eqref{eq:hc-schubert} and \eqref{eq:horo-diffeo}, 
the action of $\vartheta_{xx'}$ on $C(\tau)$ 
corresponds to the action of $\vartheta_{xx'}\circ \pi^{\tau}_{x'x}$ on $Hc^{\tau}_x$, 
and the differential $(d\vartheta_{xx'})_{\hat\tau}$ of $\vartheta_{xx'}$ at $\hat\tau$ 
to the differential of $\vartheta_{xx'}\circ \pi^{\tau}_{x'x}$ at $x$. 

We first consider the case when $\vartheta_{xx'}$ when $\xi:=l_{xx'}(-\infty)\in\ost(\tau)$,
equivalently, when $x'$ lies in the interior of the Weyl cone $V(x,\st(\tau))$. 
Then $(d\vartheta_{xx'})_{\hat\tau}$ strictly contracts:
\begin{lem}
\label{lem:diffposev}
If $\xi\in\ost(\tau)$,
then $(d\vartheta_{xx'})_{\hat\tau}$ 
is diagonalizable with eigenvalues in $(0,1)$.
\end{lem}
\proof
Since $\xi\in\ost(\tau)$,
the group $N_{\tau}$ preserves the strong asymptote classes of geodesic rays at $\xi$,\footnote{However, 
$N_{\tau}$ does not act transitively on it, unless $\xi\in\inte(\tau)$.}
cf.\ section~\ref{sec:horocycles}, 
i.e.\ the geodesics $nl_{xx'}$ for $n\in N_{\tau}$ are strongly backward asymptotic to $l_{xx'}$.
Thus,
by assigning to $n\hat\tau\in C(\tau)$ the geodesic $nl_{xx'}$,
which is the unique geodesic in the parallel set $P(\tau,n\hat\tau)$
strongly backward asymptotic to $l_{xx'}$,
we obtain a smooth family of geodesics in the strong backward asymptote class of $l_{xx'}$,
parametrized by the manifold $C(\tau)$.

By differentiating this family,
we obtain a linear embedding of the tangent space $T_{\hat\tau} C(\tau)$ 
into the vector space $\Jac_{l_{xx'},\xi}$ of Jacobi fields along $l_{xx'}$ 
which decay to zero at $\xi$. 
The effect of the differential $(d\vartheta_{xx'})_{\hat\tau}$ on $C(\tau)$ is given,
in terms of these Jacobi fields, by the push-forward
\begin{equation*}
J\mapsto (\vartheta_{xx'})_*(J) = d\vartheta_{xx'}\circ J\circ \vartheta_{x'x}
\end{equation*}
The Jacobi fields in $\Jac_{l,\xi}$, 
which are of the form of a decaying exponential function times a parallel vector field along $l_{xx'}$, 
correspond to the eigenvectors of $(d\vartheta_{xx'})_{\hat\tau}$ with eigenvalues in $(0,1)$. 
It is a standard fact from the Riemannian geometry of symmetric spaces
that the vector space $\Jac_{l_{xx'},\xi}$ has a basis consisting of such special Jacobi fields.\footnote{A transvection 
along a geodesic acts
on the space of Jacobi fields along this geodesic as a diagonalizable transformation,
see \cite{Eberlein, Helgason}.}
The same then follows for the linear subspace $L\subseteq\Jac_{l_{xx'},\xi}$ corresponding to $T_{\hat\tau} C(\tau)$.
Thus the eigenvectors of $(d\vartheta_{xx'})_{\hat\tau}$ for positive eigenvalues span $T_{\hat\tau}C(\tau)$.
\qed

\medskip
We now give a uniform estimate for the contraction of $(d\vartheta_{xx'})_{\hat\tau}$:

\begin{lem}
\label{lem:evestinst}
If $\xi\in\ost(\tau)$,
then the eigenvalues $\la$ of $(d\vartheta_{xx'})_{\hat\tau}$ satisfy an estimate
\begin{equation}
\label{eq:evest}
-\log\la\geq c\cdot d(x',\D V(x,\st(\tau)))
\end{equation}
with a constant $c>0$ depending only on $X$.
\end{lem}
\proof
We continue the argument in the previous proof.

Let $F\supset l_{xx'}$ be a maximal flat. 
Then $F\subset P(\tau,\hat\tau)$. 
The smooth family $n\hat\tau\mapsto nl_{xx'}$ of geodesics parametrized by $C(\tau)$
embeds into the smooth family of maximal flats $n\hat\tau\mapsto nF$.
They are all asymptotic to $\st(\tau)\cap\geo F$,
i.e.\ $\geo(nF)\supset\st(\tau)\cap\geo F$.
Accordingly, each Jacobi field $J\in L\subseteq\Jac_{l_{xx'},\xi}$ 
extends to a Jacobi field $\hat J$ along $F$ which decays to zero at all ideal points in $\ost(\tau)\cap\geo F$.
(Here we use again that 
$N_{\tau}$ preserves the strong asymptote classes of geodesic rays at all points in $\ost(\tau)$.)
Thus, we obtain a natural identification of $T_{\hat\tau} C(\tau)$ and $L$ 
with a linear subspace $\hat L$ 
of the vector space $\Jac_{F,\ost(\tau)\cap\geo F}$
of Jacobi fields along $F$ 
which decay to zero at all ideal points in $\ost(\tau)\cap\geo F$.

The decomposition of Jacobi fields mentioned in the previous proof
works in the same way along flats.\footnote{As in the case of geodesics,
a transvection along a flat acts
on the space of Jacobi fields along this flat as a diagonalizable transformation,
see \cite{Eberlein, Helgason}.}
The vector space $\Jac_{F,\ost(\tau)\cap\geo F}$
has a basis consisting of Jacobi fields of the form $e^{-\al}V$ 
with an affine linear form $\al$ on $F$ 
and a parallel vector field $V$ along $F$. 
Furthermore,
since $G$ acts transitively on maximal flats,
only {\em finitely} many affine linear forms $\al$ occur for these basis elements,
independently of $F$.
(The possible forms are determined by the root system of $G$, but we do not need this fact here.)

The decay condition on the forms $\al$ occurring in our decomposition
is equivalent to the property that 
$\al\geq\al(x)$ on $V(x,\st(\tau)\cap\geo F)\subset F$ 
and $\al>\al(x)$ on the interior of this cone.
It implies an estimate 
\begin{equation*}
\al(x')-\al(x)\geq c\cdot \underbrace{d(x',\D V(x,\st(\tau)\cap\geo F))}_{=d(x',\D V(x,\st(\tau)))}
\end{equation*}
with a constant $c=c(\al)>0$.
(The equality of distances follows from Proposition \ref{prop:hdstwcones}.) 
Since there are only finitely many forms $\al$ involved,
the constant $c$ can be taken 
{\em independent} of $\al$.

Notice that 
the eigenvalues $\la$ of $(d\vartheta_{xx'})_{\hat\tau}$ are of the form 
$$ e^{-(\al(x')-\al(x))} .$$
The claimed upper bound for the eigenvalues follows.
\qed

\medskip
By continuity, the result extends to the case 
when $x'$ lies in the boundary of the Weyl cone $V(x,\st(\tau))$. 
We obtain:
\begin{cor}
\label{cor:wcontrincone}
If $x'\in V(x,\st(\tau))$, 
then $(d\vartheta_{xx'})_{\hat\tau}$ is diagonalizable with eigenvalues in $(0,1]$
satisfying an estimate \eqref{eq:evest}.

In particular, 
the eigenvalues lie in $(0,1)$,
if $x'$ lies in the interior of $V(x,\st(\tau))$.
\end{cor}

If $x'$ lies outside the Weyl cone $V(x,\st(\tau))$, 
then $d(\vartheta_{xx'})_{\hat\tau}$ has expanding directions.
In order to see this, 
we consider the action of $\vartheta_{xx'}$
on certain invariant submanifolds of $C(\tau)$
corresponding to parallel sets of singular hyperplanes.

Again, there exists a maximal flat $F$ with $l_{xx'}\subset F\subset  P(\tau,\hat\tau)$.
Let $h\subset\geo F$ be a half-apartment such that $\inte(\tau)\subset\inte(h)$. 
Then $l_{xx'}\subset F\subset P(\D h)$.
The transvection $\vartheta_{xx'}$ fixes $\geo F$ pointwise.
Hence it preserves the parallel set $P(\D h)$ and the submanifold 
$C(h)=N_h\hat\tau\subset C(\tau)$. 

If $l_{xx'}$ is parallel to the euclidean factor of $P(\D h)$,
equivalently, if $\geo l_{xx'}\subset\D h$,
then $\vartheta_{xx'}$ acts trivially on $\geo P(\D h)$.
Hence, 
$\vartheta_{xx'}$ acts also trivially on $C(h)$,
because the latter consists of simplices contained in $\geo P(\D h)$.

In the general case,
the action of $\vartheta_{xx'}$ on $C(h)$ corresponds to the restriction of the action of 
$\vartheta_{xx'}\circ \pi^{\tau}_{x'x}$ to $\Hc_{h,x}=\Hc_{\tau,x}\cap P(\D h)$.
When projecting to $CS(\D h)$, 
the latter action in turn corresponds to the action of 
$\vartheta_{\bar x\bar x'}\circ \pi^{\zeta}_{\bar x'\bar x}$
on the horosphere $\Hs_{\zeta,\bar x}$. 
Here, $\vartheta_{\bar x\bar x'}$ denotes the transvection on $CS(\D h)$ 
with axis $l_{\bar x\bar x'}$ through $\bar x$ and $\bar x'$ 
mapping $\bar x'\mapsto\bar x$, 
and $ \pi^{\zeta}_{\bar x'\bar x}$ is the natural identification \eqref{eq:hrsphidfcrsc}.
The axis $l_{\bar x\bar x'}$ is the image of $F$ under the projection
(if $\bar x=\bar x'$, we define it in this way).
It is asymptotic to $\zeta$ and another ideal point $\hat\zeta\in C(\zeta)=\geo CS(\D h)-\{\zeta\}$.
The simplex $\hat\tau$ corresponds to $\hat\zeta$ 
under the natural $N_h$-equivariant identification
$C(h)\cong C(\zeta)$,
and the action of $\vartheta_{xx'}$ on $C(h)$ corresponds to the action of $\vartheta_{\bar x\bar x'}$ on $C(\zeta)$. 

We now obtain analogues of Lemmas~\ref{lem:diffposev} and~\ref{lem:evestinst}.
Recall that $\xi=l_{xx'}(-\infty)$.
\begin{lem}
\label{lem:evestinstrk1}
If $\xi\in\inte(h)$, 
then $(d\vartheta_{xx'})_{\hat\tau}|_{T_{\hat\tau}C(h)}$ 
is diagonalizable with eigenvalues $\la\in (0,1)$
satisfying an estimate
\begin{equation}
\label{ineq:eigvsubsp}
c_2 \leq\frac{-\log\la}{b_{\zeta}(x)-b_{\zeta}(x')}\leq c_1
\end{equation}
with constants $c_1,c_2>0$ depending only on $X$.
\end{lem}
\proof
The diagonalizablility follows by applying Lemma~\ref{lem:diffposev} 
to $CS(\D h)$ and $(d\vartheta_{\bar x\bar x'})_{\hat\zeta}$. 

Since $\xi\in\interior(h)$, we have that 
$b_{\zeta}(x)-b_{\zeta}(x')=b_{\zeta}(\bar x)-b_{\zeta}(\bar x')>0$,
and the eigenvalue estimate follows from the contraction estimate (\ref{ineq:contrestrk1})
\qed

\begin{cor}
\label{cor:diffnotcontr}
If $x'\in P(\tau,\hat\tau)-V(x,\st(\tau))$, 
then $(d\vartheta_{xx'})_{\hat\tau}$ 
has some eigenvalues 
in $(1,+\infty)$.
\end{cor}
\proof
By our assumption, we have that $\xi\not\in\st(\tau)$.
Therefore, the half-apartment $h\subset\geo F$ can be chosen so that its interior contains, besides $\interior(\tau)$, 
also $l_{xx'}(+\infty)$.
(Recall that the convex subcomplex $\st(\tau)\cap\geo F$ is an intersection of half-apartments in $\geo F$, cf. Lemma \ref{lem:intapts}.)
Then the estimate (\ref{ineq:eigvsubsp}) applied to $\vartheta_{x'x}=\vartheta_{xx'}^{-1}$ 
yields that $(d\vartheta_{xx'})_{\hat\tau}^{-1}$ has some eigenvalues in $(0,1)$.
\qed

\medskip
Complementing Corollary~\ref{cor:wcontrincone},
we bound the contraction rate from above,
if $x'\in V(x,\st(\tau))$:

\begin{lem}
\label{lem:contrestinst}
If $\xi\in\st(\tau)$,
then $(d\vartheta_{xx'})_{\hat\tau}$ has some eigenvalue $\la\in(0,1]$ 
satisfying an estimate
\begin{equation*}
-\log\la\leq c_1\cdot d(x',\D V(x,\st(\tau)))
\end{equation*}
with a constant $c_1>0$ depending only on $X$.
\end{lem}
\proof Since $xx'\subset F$, 
some nearest point $y'$ to $x'$ on $\D V(x,\st(\tau))$ lies in $F$, 
cf.\ Proposition \ref{prop:hdstwcones}.
Hence we can choose the half-apartment $h\subset\geo F$ 
such that $b_{\zeta}(y')=b_{\zeta}(x)$ and 
\begin{equation*}
d(x',\D V(x,\st(\tau))) =
b_{\zeta}(x)-b_{\zeta}(x').
\end{equation*}
Now let $\la$ be an eigenvalue of $(d\vartheta_{xx'})_{\hat\tau}|_{T_{\hat\tau}C(h)}$
and apply the upper bound in (\ref{ineq:eigvsubsp}).
\qed

\medskip
Putting the information 
(Corollaries~\ref{cor:wcontrincone}, \ref{cor:diffnotcontr} 
and Lemmas~\ref{lem:evestinstrk1}, \ref{lem:contrestinst})
together, we obtain:
\begin{prop}[Infinitesimal contraction of transvections at infinity]
\label{prop:infcontrtrans}
Let $\tau,\hat\tau\subset\geo X$ be opposite simplices,
and let $\vartheta$ be a nontrivial transvection with an axis $l\subset P(\tau,\hat\tau)$
through the point $x$. 
Then the following hold for the differential $d\vartheta_{\hat\tau}$ of $\vartheta$ on $C(\tau)$
at the fixed point $\hat\tau$:

(i) 
$d\vartheta_{\hat\tau}$ is 
diagonalizable with eigenvalues in $(0,1]$ 
iff $\vartheta^{-1}x\in V(x,\st(\tau))$,
and 
diagonalizable with eigenvalues in $(0,1)$ 
iff $\vartheta^{-1}x\in V(x,\ost(\tau))$.

(ii)
If $\vartheta^{-1}x\in V(x,\st(\tau))$,
then the eigenvalues $\la$ of $d\vartheta_{\hat\tau}$
satisfy an estimate 
\begin{equation*}
c_2\cdot d(\vartheta^{-1}x,\D V(x,\st(\tau))) \leq
-\log\la
\leq c_1\cdot d(\vartheta^{-1}x,\D V(x,\st(\tau)))
\end{equation*}
with constants $c_1,c_2>0$ depending only on $X$.
\end{prop}

We deduce a consequence for the action of general isometries in $G$.
For later use, 
we will formulate it in terms of expansion (of their inverses) rather than contraction.

We need the following notion:
For a diffeomorphism $\Phi$ of a Riemannian manifold $M$,
we define the {\em expansion factor} at $x\in M$ as
\begin{equation}
\label{eq:expfacearl}
\eps(\Phi,x) = \inf_{v\in T_xM-\{0\}} \frac{\|d\Phi(v)\|}{\|v\|} = \| (d\Phi_x)^{-1} \|^{-1} ,
\end{equation}
compare \eqref{eq:expfac} in section~\ref{sec:trexpand} below.

We equip the flag manifolds $\Flagt$ with auxiliary Riemannian metrics.

\begin{thm}[Infinitesimal expansion of isometries at infinity]
\label{thm:expand}
Let $\tau\in\Flagt$,
$x\in X$, and $g\in G$ such that 
$d(gx,V(x,\st(\tau)))\leq r$.
Then for the action of $g^{-1}$ on $\Flagt$ we have the estimate 
\begin{equation*}
C^{-1}\cdot d(gx,\D V(x,\st(\tau))) - A\leq
\log \eps(g^{-1},\tau)
\leq C\cdot d(gx,\D V(x,\st(\tau))) + A
\end{equation*}
with constants $C,A>0$ depending only on $x$, $r$ 
and the chosen Riemannian metric on $\Flagt$.\footnote{The estimate depends also on the point $x$
because the choice of the auxiliary metric on $\Flagt$ reduces the symmetry:
The action of a compact subgroup of $G$ on $\Flagt$ is uniformly bilipschitz,
but not the $G$-action.}
\end{thm}
\proof 
We write $g$ as a product $g=t b$ 
of a transvection $t$ along a geodesic $l$ through $x$ with $l(+\infty)\in\st(\tau)$ 
and an isometry $b\in G$ such that $d(x,bx)\leq r$.
Then $t$ fixes $\tau$ on $\Flagt$, 
and the expansion factor $\eps(g^{-1},\tau)$ equals 
$\eps(t^{-1},\tau)$
up to a multiplicative constant depending on $r$ and the chosen Riemannian metric on $\Flagt$.

When replacing the metric, $\eps(t^{-1},\tau)$ changes at most by another multiplicative constant,
and we may therefore assume that the Riemannian metric is invariant under the maximal compact subgroup $K_x<G$ 
fixing $x$.
Now the eigenspace decomposition of $dt_{\tau}$ on $T_{\tau}\Flagt$ is orthogonal.
Consequently,
$$\eps(t^{-1},\tau)=\la_{max}^{-1}$$
where $\la_{max}$ denotes the maximal eigenvalue of $dt_{\tau}$.

Let $\hat\tau$ denote the simplex $x$-opposite to $\tau$.
Applying Proposition~\ref{prop:infcontrtrans}(ii) to 
$\vartheta=t$ 
while exchanging the roles of $\tau$ and $\hat\tau$,
we obtain the estimate 
\begin{equation*}
c_2\cdot d(t^{-1}x,\D V(x,\st(\hat\tau))) \leq
-\log\la
\leq c_1\cdot \underbrace{d(t^{-1}x,\D V(x,\st(\hat\tau)))}_{=d(tx,\D V(x,\st(\tau)))}
\end{equation*}
for the eigenvalues $\la$ of $dt_{\tau}$,
and so
$$ c_2\cdot d(tx,\D V(x,\st(\tau))) \leq \log\eps(t^{-1},\tau) \leq c_1\cdot d(tx,\D V(x,\st(\tau))) ,$$
which is the desired estimate.
\qed

\medskip
Let us now consider sequences $(g_n)$ in $G$.
The theorem can be used to draw conclusions 
from the expansion behavior at infinity of the sequence of inverses $(g_n^{-1})$
on the geometry of an orbit sequence $(g_nx)$ in $X$:
If $(g_nx)$ lies in a tubular neighborhood of the Weyl cone $V(x,\st(\tau))$,
then the expansion factors $\eps(g_n^{-1},\tau)$ on $\Flagt$ 
are bounded below,
and their logarithms measure the distance of $(g_nx)$ to the boundary of the Weyl cone. 
In particular, if the expansion factors diverge, $\eps(g_n^{-1},\tau)\to+\infty$,
then (the projection of) $(g_nx)$ enters deep into the cone $V(x,\st(\tau))$.

The next result shows how to recognize from expansion
whether the orbit sequence $(g_nx)$ remains in a tubular neighborhood of the Weyl cone $V(x,\st(\tau))$,
once it stays close to the parallel set spanned by it:

\begin{prop}
\label{prop:expand}
Let $\tau,\hat\tau\subset\geo X$ be opposite simplices.
Suppose that $(g_n)$ is a sequence in $G$ such that, 
for some point $x\in X$, 
the sequence $(g_nx)$ is contained in a tubular neighborhood of the parallel set $P(\tau,\hat\tau)$,
but drifts away from the Weyl cone $V(x,\st(\tau))$,
\begin{equation*}
d(g_nx,V(x,\st(\tau)))\to+\infty
\end{equation*}
as $n\to+\infty$. 
Then $\eps(g_n^{-1}, \tau)\to0$.
\end{prop}
\proof 
We may assume that $x\in P=P(\tau,\hat\tau)$. 
As in the proof of Theorem~\ref{thm:expand}, 
we can reduce to the case that the $g_n$ are transvections along geodesics $l_n$ in $P$
through the point $x$.
We need to show that the differentials $(dg_n^{-1})_{\tau}$ on $\Flagt$ have (some) small eigenvalues,
i.e.\ that their minimal eigenvalue goes $\to0$.

We proceed as in the proof of Corollary~\ref{cor:diffnotcontr}.
Let $F_n\subset P$ be a maximal flat containing $l_n$.
Then also 
\begin{equation*}
d(g_nx,V(x,\st(\tau))\cap F_n)\to+\infty ,
\end{equation*}
cf.\ Proposition~\ref{prop:hdstwcones}.
There exist half-apartments $h_n\subset\geo F_n$ with centers $\zeta_n$,
so that $b_{\zeta_n}\leq b_{\zeta_n}(x)$ on $V(x,\st(\tau))\cap F_n$ (and hence also on $V(x,\st(\tau))$)
and $b_{\zeta_n}(g_nx)-b_{\zeta_n}(x)\to+\infty$.
Let $\hat h_n\subset\geo F_n$ denote the complementary half-apartments, $\D\hat h_n=\D h_n$,
and $\hat\zeta_n$ their centers.
Then $b_{\zeta_n}+b_{\hat\zeta_n}\equiv const$ on $F_n$.
It suffices to show that the differentials $(dg_n^{-1})_{\tau}$ are contracting on the invariant subspaces 
$T_{\tau}C(\hat h_n)\subseteq T_{\tau}C(\hat\tau)$
with norms going $\to0$.
According to Lemma~\ref{lem:evestinstrk1},
the eigenvalues of $(dg_n^{-1})_{\tau}|_{T_{\tau}C(\hat h_n)}$
are positive and bounded above by 
$$ e^{-c_2(b_{\hat\zeta_n}(x)-b_{\hat\zeta_n}(g_nx))} = e^{-c_2(b_{\zeta_n}(g_nx)-b_{\zeta_n}(x))}\to0.$$
This finishes the proof.
\qed

\subsection{Finsler geodesics}
\label{sec:fins}

We will work with the following notion of Finsler geodesic:
\begin{dfn}[Finsler geodesics]
\label{def:finsgeo}
A continuous path $c:I\to X$ is a {\em $\taumod$-Finsler geodesic} if
it is contained in a parallel set $P(\tau_-,\tau_+)$ with $\tau_{\pm}\in\Flagpmt$
such that 
\begin{equation}
\label{eq:finsgeo}
c(t_+)\in V(c(t_-),\st(\tau_+))
\end{equation}
for all subintervals $[t_-,t_+]\subseteq I$.
It is {\em $\Theta$-regular}
if, moreover,  
\begin{equation}
\label{eq:finsgeoth}
c(t_+)\in V(c(t_-),\stTh(\tau_+))
\end{equation}
We call a $\taumod$-Finsler geodesic {\em uniformly $\taumod$-regular}
if it is $\Theta$-regular for some $W_{\taumod}$-convex compact subset $\Theta\subset \inte_{\taumod}(\simod)$.
\end{dfn}
Note that we do not require the parameterization of Finsler geodesics to be by arc length. 
The terminology is justified by the fact that $\taumod$-Finsler geodesics 
are (up to parameterization) the geodesics for certain $G$-invariant ``polyhedral'' Finsler metrics,
see \cite[\S 5.1.3]{bordif}.

The condition (\ref{eq:finsgeo}) is equivalent to $c(t_-)\in V(c(t_+),\st(\tau_-))$,
and it follows that 
the subpaths $c|_{[t_-,t_+]}$ are contained in the diamonds $\diamot(c(t_-),c(t_+))$.
Similarly,
(\ref{eq:finsgeoth}) is equivalent to $c(t_-)\in V(c(t_+),\stTh(\tau_-))$,
because $\Theta$ is assumed $\iota$-invariant,
and in the $\Theta$-regular case $c|_{[t_-,t_+]}$ is contained in $\diamoTh(c(t_-),c(t_+))$.

It is worth mentioning the following Finsler geometric interpretation of {\em diamonds}: 
They are Finsler versions of Riemannian geodesic segments 
in the sense that the union of all $\taumod$-Finsler geodesic segments with endpoints $x_{\pm}$
fills out $\diamot(x_-,x_+)$,
see also \cite[\S 5.1.3]{bordif}.

\medskip
We now discuss the ``drift'' component of $\taumod$ Finsler geodesics.

We work with the vector valued distance $d_{\taumod} =\pi^{\De}_{\taumod} \circ d_{\De}$.
introduced in (\ref{eq:vvdt}).
We first consider the case of broken geodesics $xyz$ which are $\taumod$-Finsler geodesics:
\begin{lem}[Additivity]
Let $\tau\in\Flagt$.
If $y\in V(x,\st(\tau))$ and $z\in V(y,\st(\tau))$, then 
$$ d_{\taumod}(x,y)+d_{\taumod}(y,z) =d_{\taumod}(x,z) .$$
\end{lem}
\proof
The $\taumod$-distance can be expressed 
in terms of the projections of Weyl cones to their central sectors.
Consider the nearest point projection
$$ \pi_{x,\tau}:V(x,\st(\tau))\to V(x,\tau),$$
cf.\ (\ref{eq:prjwcocntsect}). 
Note that it coincides with the nearest point projection from $V(x,\st(\tau))$ to the singular flat 
spanned by the sector $V(x,\tau)$,
compare Lemma~\ref{lem:proj} and the comment thereafter.
Then
$$ d_{\taumod}(x,\cdot) =d_{\De}(x,\pi_{x,\tau}(\cdot)) $$
on $V(x,\st(\tau))$. 

In order to relate $d_{\taumod}(y,z)$ to $d_{\taumod}(x,y)$ and $d_{\taumod}(x,z)$,
we observe that 
the sectors $V(y,\tau)$ and $V(\pi_{x,\tau}(y),\tau)\subseteq V(x,\tau)$ are parallel
and isometrically identified by $\pi_{x,\tau}$.
Moreover, 
$$\pi_{x,\tau}|_{V(y,\st(\tau))}=(\pi_{x,\tau}|_{V(y,\tau)})\circ\pi_{y,\tau} .$$
Therefore,
$$ d_{\taumod}(y,z) = d_{\De}(y,\pi_{y,\tau}(z)) = d_{\De}(\pi_{x,\tau}(y),\pi_{x,\tau}(z))  .$$
The additivity formula follows in view of the nestedness $\pi_{x,\tau}(z)\in V(\pi_{x,\tau}(y),\tau)$.
\qed

\medskip
Applying the lemma to $\taumod$-Finsler geodesics yields:

\begin{prop}[Additivity of $\taumod$-distance along Finsler geodesics]
If $c:I\to X$ is a $\taumod$-Finsler geodesic,
then
$$ d_{\taumod}(c(t_0),c(t_1))+d_{\taumod}(c(t_1),c(t_2))=d_{\taumod}(c(t_0),c(t_2))$$
for all $t_0\leq t_1\leq t_2$ in $I$.
\end{prop}
We reformulate this as:
\begin{prop}
[$\taumod$-projection of Finsler geodesics]
\label{prop:tmdstalfins}
If $c:[0,T]\to X$ is a $\taumod$-Finsler geodesic,
then so is
$$ \bar c_{\taumod}:= d_{\taumod}(c(0),c) : [0,T]\to V(0,\taumod),$$
and 
$$ \bar c_{\taumod}(t_2) = \bar c_{\taumod}(t_1)+d_{\taumod}(c(t_1),c(t_2)) $$
for all $0\leq t_1\leq t_2\leq T$.
\end{prop}
Note that the equality in the last proposition implies:
\begin{equation} 
\label{eq:tmantdl}
d(\bar c_{\taumod}(t_1),\bar c_{\taumod}(t_2)) = \| d_{\taumod}(c(t_1),c(t_2)) \|
\end{equation}

\medskip
We now study the {\em $\De$-distance along Finsler geodesics}.

This is based on Proposition~\ref{prop:prjwcon} 
which concerns the $\De$-side lengths of triangles $\De(x,y,z)$ in $X$
such that the broken geodesic $xyz$ is a Finsler geodesic.
Applying this proposition to Finsler geodesics, we obtain our main result concerning their geometry:
\begin{thm}[$\De$-projection of Finsler geodesics]
\label{thm:dedstalfins}
(i) If $c:[0,T]\to X$ is a $\taumod$-Fins\-ler geodesic, 
then so is 
$$ \bar c_{\De}:=d_{\De}(c(0),c):[0,T]\to \De .$$

(ii) If $c$ is also $\Theta$-regular, with $\Theta\subset\inte_{\taumod}(\simod)$ compact and $\taumod$-Weyl convex,
then so is $\bar c_{\De}$.
Moreover, the distances between points on $c$ and $\bar c_{\De}$ are comparable:
$$ d(\bar c_{\De}(t_1),\bar c_{\De}(t_2)) \geq \eps(\Theta)\cdot d(c(t_1),c(t_2)) $$
for $0\leq t_1\leq t_2\leq T$ with a constant $\eps(\Theta)>0$.
\end{thm}
We note that $d(\bar c_{\De}(t_1),\bar c_{\De}(t_2))\leq d(c(t_1),c(t_2))$,
because $d_{\De}(c(0),\cdot)$ is 1-Lipschitz.
\proof
(i)
Applying Proposition~\ref{prop:prjwcon} to the triangles 
$\De(c(0),c(t_1),c(t_2))$, $0\leq t_1\leq t_2\leq T$,
yields
\begin{equation*}
\bar c_{\De}(t_2) \in V(\bar c_{\De}(t_1),\Wt\De),
\end{equation*}
the cone being understood as a subset of $\Fmod$,
which means that $\bar c_{\De}$ is a $\taumod$-Finsler geodesic.

(ii) That $\bar c_{\De}$ is now $\Theta$-regular, follows similarly.
The comparability of distances we deduce using our earlier discussion of $\taumod$-distances along Finsler geodesics. 
We estimate:
$$ d(\bar c_{\De}(t_1),\bar c_{\De}(t_2)) 
\geq d(\bar c_{\taumod}(t_1),\bar c_{\taumod}(t_2))
= \|d_{\taumod}(c(t_1),c(t_2))\| 
\geq \eps(\Theta)\cdot d(c(t_1),c(t_2)) $$
The first inequality holds, because $\bar c_{\taumod} = \pi^{\De}_{\taumod} \circ \bar c_{\De}$ 
and $\pi^{\De}_{\taumod}$ is 1-Lipschitz. 
The equality follows from (\ref{eq:tmantdl}).
The last inequality comes from the lower bound for the length of the $\taumod$-component of $\Theta$-regular segments,
cf.\ (\ref{eq:urgtdst}).
\qed

\section{Topological dynamics}

\subsection{Expansion}
\label{sec:trexpand}

Let first $Z$ be a metric space
and let $\Ga\acts Z$ be a continuous action by a discrete group.
We will use the following notions of {\em metric expansion}, compare \cite[\S 9]{Sullivan}:

\begin{dfn}[Metric expansion]
\label{def:metexpan}
(i) 
A homeomorphism $h$ of $Z$ is {\em expanding} at a point $z\in Z$
if there exists a neighborhood $U$ of $z$ and a constant $c>1$
such that $h|_U$ is $c$-expanding in the sense that 
\begin{equation*}
d(h z_1,h z_2)\geq c\cdot d(z_1,z_2) .
\end{equation*}
for all points $z_1,z_2\in U$.

(ii)
A sequence of homeomorphisms $h_n$ of $Z$ has {\em diverging expansion} at the point $z\in Z$
if there exists a sequence of neighborhoods $U_n$ of $z$
and numbers $c_n\to+\infty$ such that $h_n|_{U_n}$ is $c_n$-expanding.

(iii) 
The action $\Ga\acts Z$ is {\em expanding} 
at 
$z\in Z$ 
if there exists an element $\ga\in\Ga$ 
which is expanding at $z$.
The action has {\em diverging expansion} at $z\in Z$
if $\Ga$ contains a sequence which has diverging expansion at $z$.

(iv)
The action $\Ga\acts Z$ is {\em expanding} 
at a compact $\Ga$-invariant {\em subset} $E\subset Z$ 
if it is expanding at all points $z\in E$. 
\end{dfn}

We observe that the properties of diverging expansion depend only on the bilipschitz class of the metric. 
Furthermore,
if an action is expanding at an invariant compact subset
then, due to iteration, it has diverging expansion at every point of the subset.

\medskip
Now let $M$ be a Riemannian manifold and let $\Ga\acts M$ be a smooth action.
There are infinitesimal analogs of the above expansion conditions. 

We recall from \eqref{eq:expfacearl} that,
for a diffeomorphism $\Phi$ of $M$, 
the {\em expansion factor} $\eps(\Phi,x)$ at a point $x\in M$ is defined as:
\begin{equation}
\label{eq:expfac}
\eps(\Phi,x) = \inf_{v\in T_xM-\{0\}} \frac{\|d\Phi(v)\|}{\|v\|} = \| (d\Phi_x)^{-1} \|^{-1} 
\end{equation}

\begin{definition}[Infinitesimal expansion]
\label{def:infexpan}
(i) 
A diffeomorphism $\Phi$ of $M$ is {\em infinitesimally expanding} at a point $x\in M$ 
if $\eps(\Phi,x)>1$.

(ii)
A sequence of diffeomorphisms $\Phi_n$ of $M$ has {\em diverging infinitesimal expansion} at $x$
if $\eps(\Phi_n,x)\to+\infty$ as $n\to+\infty$.

(iii) 
The action $\Ga\acts M$ is {\em infinitesimally expanding} 
at $x$
if there exists an element $\ga\in\Ga$ 
which is infinitesimally expanding at $x$.
The action has {\em diverging infinitesimal expansion} at $x$
if $\Ga$ contains a sequence which has diverging infinitesimal expansion at $x$.

(iv)
The action $\Ga\acts M$ is {\em infinitesimally expanding} 
at a compact $\Ga$-invariant {\em subset} $E\subset M$ 
if it is infinitesimally expanding at all points $x\in M$. 
\end{definition}

If the manifold $M$ is compact,
the properties of diverging infinitesimal expansion are independent of the Riemannian metric.
In the general case,
if an action is infinitesimally expanding at an invariant compact subset
then it has diverging infinitesimal expansion at every point of the subset.

We note that for smooth actions on Riemannian manifolds infinitesimal and metric expansion are equivalent.

\subsection{Discontinuity and dynamical relation}
\label{sec:discdynrl}

Let $Z$ be a compact metrizable space,
and let $\Ga<\Homeo(Z)$ be a countably infinite subgroup (although in the definition of a proper action below we allow for subsemigroups).
We consider the action $\Ga\acts Z$.
\begin{dfn}[Discontinuous]
\label{dfn:wand}
A point $z\in Z$ is called {\em wandering} 
with respect to the $\Ga$-action 
if the action is {\em discontinuous} at $z$,
i.e.\ if $z$ has a neighborhood $U$ such that 
$U\cap\ga U\neq\emptyset$ for at most finitely many $\ga\in\Ga$.
\end{dfn}
Nonwandering points are called {\em recurrent}. 

\begin{dfn}[Domain of discontinuity]
\label{dfn:domdisc}
We call the set 
$$\Om_{disc}\subset Z$$ 
of wandering points 
the {\em wandering set} or {\em domain of discontinuity} for the action $\Ga\acts Z$.
\end{dfn}
Note that $\Om_{disc}$ is open and $\Ga$-invariant.

\begin{dfn}[Proper]
The action of a {\em subsemigroup} $\Ga<  Homeo(X)$ on an open 
subset $U\subset Z$
is called {\em proper} if
for every compact subset $K\subset U$ 
$K\cap\ga K\neq\emptyset$ for at most finitely many $\ga\in\Ga$.
\end{dfn}

If $\Ga$ is a subgroup of $Homeo(X)$ acting properly discontinuously on $U\subset X$ 
then the action of $\Ga$ on $U$ is then discontinuous, $U\subseteq\Om_{disc}$,
and therefore is  called {\em properly discontinuous}. 

\begin{dfn}[Domain of proper discontinuity]
\label{dfn:dompropdisc}
If $\Ga< Homeo(X)$ is a subgroup, 
we call a $\Ga$-invariant open subset $\Om\subseteq\Om_{disc}$ on which $\Ga$ acts properly
a {\em domain of proper discontinuity} for $\Ga$.
\end{dfn}

The orbit space $\Om/\Ga$ is then Hausdorff. 
Note that in general there is {\em no unique} maximal {\em proper} domain of discontinuity. 

Discontinuity and proper discontinuity can be nicely expressed 
using the notion of dynamical relation.
The following definition is due to Frances \cite[Def.\ 1]{Frances}:
\begin{dfn}[Dynamically related]
Two points $z,z'\in Z$ 
are called 
{\em dynamically related} with respect to a sequence $(h_n)$ in $\Homeo(Z)$,
$$ z\stackrel{(h_n)}{\sim}z' $$
if there exists a sequence $z_n\to z$ in $Z$ such that $h_nz_n\to z'$. 

The points $z,z'$ are called 
{\em dynamically related} with respect to the $\Ga$-action, 
$$ z\stackrel{\Ga}{\sim}z' $$
if there exists a sequence $\ga_n\to\infty$ in $\Ga$ such that $z\stackrel{(\ga_n)}{\sim}z'$.
\end{dfn}
Here, for a sequence $(\ga_n)$ in $\Ga$ we write $\ga_n\to\infty$
if every element of $\Ga$ occurs at most finitely many times in the sequence. 

One verifies (see e.g. \cite{manicures}):

(i) Dynamical relation is a closed relation in $Z\times Z$. 

(ii) Points in different $\Ga$-orbits are dynamically related 
if and only if their orbits cannot be separated 
by disjoint $\Ga$-invariant open subsets. 

The concept of dynamical relation is useful for our discussion of discontinuity, because:

(i) A point is nonwandering if and only if it is dynamically related to itself. 

(ii) The action is proper on an open subset $U\subset Z$ 
if and only if no two points in $U$ are dynamically related.

\subsection{Convergence groups}
\label{sec:convdy}

Let $Z$ be a compact metrizable space with at least three points.

A sequence $(h_n)$ in $\Homeo(Z)$ is {\em contracting}
if there exist points $z_{\pm}\in Z$ such that 
\begin{equation}
\label{eq:contrtaucv}
h_n|_{Z-\{z_-\}}\to z_+
\end{equation} 
uniformly on compacts as $n\to+\infty$.
Equivalently, there is no dynamical relation 
$z\stackrel{(h_n)}{\sim}z'$ between points $z\neq z_-$ and $z'\neq z_+$.
This condition is clearly symmetric, i.e.\ (\ref{eq:contrtaucv})
is equivalent to the dual condition 
that
\begin{equation}
\label{eq:contrtaucvdual}
h^{-1}_n|_{Z-\{z_+\}}\to z_-
\end{equation} 
uniformly on compacts as $n\to+\infty$.
The points $z_{\pm}$ are uniquely determined, since $|Z|\geq3$.

A sequence $(h_n)$ in $\Homeo(Z)$
is said to {\em converge} to a point $z\in Z$,
\begin{equation}
\label{eq:convseqhm}
h_n\to z
\end{equation} 
if every subsequence contains a contracting subsequence 
which, outside its exceptional point, converges to the constant map $\equiv z$.

One considers the following stronger form of convergence:

\begin{dfn}[Conical convergence]
A converging sequence $h_n\to z$
converges {\em conically},
\begin{equation}
\label{eq:convseqhmcon}
h_n\stackrel{con}{\to} z
\end{equation} 
if for some relatively compact sequence $(\hat z_n)$ in $Z-\{z\}$,
the sequence of pairs of distinct points
$h_n^{-1}(\hat z_n,z)$ is relatively compact in $(Z\times Z)^{dist}$.
\end{dfn}

Here, $(Z\times Z)^{dist}\subset Z\times Z$ denotes the complement of the diagonal.
\begin{lem}
If $h_n\stackrel{con}{\to} z$,
then
the condition in the definition 
holds 
for all relatively compact sequences $(\hat z_n)$ in $Z-\{z\}$. 
\end{lem}
\proof
Let $(\hat z_n)$ be a relatively compact sequence in $Z-\{z\}$. 
For every contracting subsequence $(h_{n_k})$
there exists a point $\hat z\in Z$ such that 
$$h^{-1}_{n_k}|_{Z-\{z\}}\to\hat z$$ uniformly on compacts.
In particular,
$h^{-1}_{n_k}\hat z_{n_k}\to\hat z$ 
and the relative compactness of $(h_{n_k}^{-1}(\hat z_{n_k},z))$ in $(Z\times Z)^{dist}$
becomes equivalent to the condition 
that the sequence $(h^{-1}_{n_k}z)$ does not accumulate at $\hat z$.
The latter condition is independent of the sequence $(\hat z_n)$.
\qed

\medskip
The following criterion for being a conical limit point of a subsequence is immediate:\footnote{Here
it suffices that $|Z|\geq2$.}
\begin{lem}
\label{lem:recogncnlimcv}
A sequence $(h_n)$ in $\Homeo(Z)$
has a subsequence conically converging to $z\in Z$
iff there exists a subsequence $(h_{n_k})$ and a point $z_-\in Z$
such that the following conditions are satisfied:

(i) 
$h_{n_k}^{-1}|_{Z-\{z\}}\to z_-$ uniformly on compacts.

(ii)
$(h_{n_k}^{-1}z)$ converges to a point different from $z_-$.
\end{lem}

\medskip
Now we pass to group actions. 

A continuous action $\Ga\acts Z$ of a 
discrete group $\Ga$ is a {\em convergence action}
if every sequence $(\ga_n)$ of pairwise distinct elements in $\Ga$ 
contains a subsequence converging to a point, equivalently,
a contracting subsequence. 
The kernel of a convergence action is finite,
and we will identify $\Ga$ with its image in $\Homeo(Z)$
which we will call a {\em convergence group}.

The {\em limit set} $\La\subset Z$ of a convergence group $\Ga<\Homeo(Z)$ is the subset of all points 
which occur as limits $z_+$ as in (\ref{eq:contrtaucv}),
equivalently,
as limits $z$ as in (\ref{eq:convseqhm})
for sequences $\ga_n\to\infty$ in $\Ga$.
The limit set is $\Ga$-invariant and compact.
A limit point $\la\in\La$ is {\em conical} if it occurs as the limit of a conically converging sequence.
A convergence group is said to have {\em conical limit set} if all limits points are conical,
and to be {\em non-elementary} if $|\La|\geq3$.
Tukia \cite[Thm.\ 2S]{Tukia_convgps} has shown that in the non-ele\-men\-ta\-ry case the limit set is perfect 
and the $\Ga$-action on it is minimal. 

If the limit set is conical, then $\Ga$ and its action on $\La$ are very special:
\begin{thm}[Bowditch {\cite{Bowditch_charhyp}}]
\label{thm:bowdchr}
Suppose that $\Ga<\Homeo(Z)$ is a non-elementary convergence group with conical limit set $\La$.
Then $\Ga$ is word hyperbolic and $\La\cong\geo\Ga$ equivariantly.
\end{thm}
The converse is easier to see:
\begin{thm}[{\cite{Gromov_hypgps,Tukia_convgps,Freden}}]
\label{thm:hypgpbdac}
The natural action of a non-virtually cyclic word hyperbolic group on its Gromov boundary 
is a minimal conical convergence action.
\end{thm} 

\subsection{Expanding convergence groups}
\label{sec:expconvgps}

The following result connects expansion with convergence dynamics.
\begin{lem}
\label{lem:conical}
If $\Gamma\acts Z$ is an expanding convergence action 
on a perfect compact metric space, 
then all points in $Z$ are conical limit points. 
\end{lem}
\proof 
We start with a general remark concerning expanding actions. 
For every point $z\in Z$ there exist an element $\ga\in\Ga$
and constants $r>0$ and $c>1$ such that 
$\ga$ is a $c$-expansion on the ball $B(z,r)$ and 
$\ga(B(z,r'))\supset B(\ga z,cr')$ for all radii $r'\leq r$. 
To see this, 
suppose that $c$ is a local expansion factor for $\ga$ at $z$ and, 
by contradiction, 
that there exist sequences of radii $r_n\to0$ 
and points $z_n\not\in B(z,r_n)$ 
such that $\ga z_n\in B(\ga z,cr_n)$. 
Then $z_n\to z$ due to the continuity of $\ga^{-1}$ and, 
for large $n$, 
we obtain a contradiction to the local $c$-expansion of $\ga$. 
Since $Z$ is compact, 
the constants $r$ and $c$ can be chosen uniformly. 
It follows by iterating expanding maps 
that for every point $z$ and every neighborhood $V$ of $z$ 
there exists $\ga\in\Ga$ such that 
$\ga(V)\supset B(\ga z,r)$, 
equivalently, 
$\ga(Z-V)\subset Z-B(\ga z,r)$. 

To verify that a point $z$ is conical, 
let $V_n$ be a shrinking sequence of neighborhoods of $z$, 
$$
\bigcap_{n} V_n= \{z\},
$$ 
and let $\ga_n\in\Ga$ be elements such that 
$\ga_n^{-1}(Z-V_n)\subset Z-B(\ga_n^{-1}z,r)$. 
Since $V_n$ is shrinking 
and $\ga_n^{-1}(V_n)\supset B(\ga_n^{-1}z,r)$ 
contains balls of uniform radius $r$, 
it follows that 
the $\ga_n^{-1}$ do not subconverge uniformly 
on any neighborhood of $z$; 
here we use that $Z$ is perfect. 
In particular, $\ga_n\to\infty$. 
The convergence action property implies that, 
after passing to a subsequence, 
the $\ga_n^{-1}$ must converge locally uniformly on $Z-\{z\}$. 
Moreover, we can assume that the sequence of points $\ga_n^{-1}z$ converges. 
By construction, 
its limit will be different 
(by distance $\geq r$) from the limit of the sequence of maps 
$\ga_n^{-1}|_{Z-\{z\}}$. 
Hence the point $z$ is conical. 
\qed

\medskip
Combining this with Bowditch's dynamical characterization of hyperbolic groups, we obtain:
\begin{cor}
If $\Gamma\acts Z$ is an expanding convergence action 
on a perfect compact metric space, 
then $\Ga$ is word hyperbolic and $Z\cong\geo\Ga$ equivariantly.
\end{cor}

Note that, conversely,
the natural action $\Ga\acts\geo\Ga$ of a word hyperbolic group $\Ga$ on its Gromov boundary
is expanding with respect to a {\em visual} metric, 
see e.g.\ \cite{CP}.

\section{Regularity and contraction}
\label{sec:regcontr}

In this section, we discuss a class of discrete subgroups of semisimple Lie groups 
which will be the framework for most of our investigations in this paper.
In particular, it contains Anosov subgroups. 
The class of subgroups will be distinguished by an 
asymptotic {\em regularity} condition
which in rank one just amounts to discreteness, but in higher rank is strictly stronger. 
The condition will be formulated in two equivalent ways.
First dynamically in terms of the action on a flag manifold, 
then geometrically in terms of the orbits in the symmetric space.

\subsection{Contraction}

Consider the action $$G\acts\Flagt$$
on the flag manifold of type $\taumod$.
Recall that for a simplex $\tau_-$ of type $\iota\taumod$ 
we denote by $C(\tau_-)\subset\Flagt$ the open dense $P_{\tau_-}$-orbit;
it consists of the simplices opposite to $\tau_-$.

We introduce the following dynamical conditions for sequences and subgroups in $G$:
\begin{dfn}[Contracting sequence]
\label{def:contracting_sequence}
A sequence $(g_n)$ in $G$ is {\em $\taumod$-con\-trac\-ting} 
if there exist simplices $\tau_{+}\in \Flagt, \tau_-\in \Flagit$ such that 
\begin{equation}
\label{eq:contrtau}
g_n|_{C(\tau_-)}\to\tau_+
\end{equation} 
uniformly on compacts as $n\to+\infty$.
\end{dfn}

\begin{dfn}[Convergence type dynamics]
\label{def:conv}
A subgroup $\Ga<G$ is a {\em $\taumod$-convergence subgroup}
if every sequence $(\ga_n)$ of distinct elements in $\Ga$ contains a $\taumod$-contrac\-ting subsequence.
\end{dfn}

Note that $\taumod$-contracting sequences diverge to infinity
and therefore $\taumod$-convergence subgroups are necessarily {\em discrete}.

A notion for sequences in $G$ equivalent to $\taumod$-contraction had been introduced by Benoist in \cite{Benoist},
see in particular part (5) of his Lemma 3.5.

The contraction property exhibits a symmetry:
\begin{lem}[Symmetry]
\label{lem:contrsym}
Property (\ref{eq:contrtau}) is equivalent to the dual property that 
\begin{equation}
\label{eq:contrtaudual}
g^{-1}_n|_{C(\tau_+)}\to\tau_-
\end{equation} 
uniformly on compacts as $n\to+\infty$.
\end{lem}
\proof
Suppose that (\ref{eq:contrtau}) holds but (\ref{eq:contrtaudual}) fails.
Equivalently, after extraction
there exists a sequence 
$\xi_n\to\xi\neq\tau_-$ in $\Flagit$ 
such that $g_n\xi_n\to\xi'\in C(\tau_+)$.
Since $\xi\neq\tau_-$,
there exists $\hat\tau_-\in C(\tau_-)$ not opposite to $\xi$.
(For instance, take an apartment in $\geo X$ containing $\tau_-$ and $\xi$, 
and let $\hat\tau_-$ be the simplex opposite to $\tau_-$ in this apartment.)
Hence there is a sequence $\tau_n\to\hat\tau_-$ in $\Flagt$ such that $\tau_n$ is not opposite to $\xi_n$ for all $n$.
(It can be obtained e.g.\ by taking a sequence $h_n\to e$ in $G$ such that $\xi_n=h_n\xi$
and putting $\tau_n=h_n\hat\tau_-$.)
Since $\hat\tau_-\in C(\tau_-)$,
condition (\ref{eq:contrtau}) implies that $g_n\tau_n\to\tau_+$.
It follows that $\tau_+$ is not opposite to $\xi'$, 
because $g_n\tau_n$ is not opposite to $g_n\xi_n$ 
and being opposite is an open condition.
This contradicts $\xi'\in C(\tau_+)$.
Therefore,
condition (\ref{eq:contrtau}) implies (\ref{eq:contrtaudual}). 
The converse implication follows by replacing the sequence $(g_n)$ with
$(g_n^{-1})$. 
\qed

\begin{lem}[Uniqueness]
\label{lem:contruniq}
The simplices $\tau_{\pm}$ in (\ref{eq:contrtau}) are uniquely determined. 
\end{lem}
\proof
Suppose that besides (\ref{eq:contrtau}) we also have $g_n|_{C(\tau'_-)}\to\tau'_+$
with simplices $\tau'_{\pm}\in\Flagpmt$.
Since the subsets $C(\tau_-)$ and $C(\tau'_-)$ are open dense in $\Flagt$,
their intersection is nonempty 
and hence $\tau'_+=\tau_+$.
Using the equivalent dual conditions (\ref{eq:contrtaudual})
we similarly obtain that $\tau'_-=\tau_-$.
\qed

\subsection{Regularity}
\label{sec:reg}

The second set of asymptotic properties concerns the geometry of the orbits in $X$.

We first consider sequences in the euclidean model Weyl chamber $\De$.
Recall that $\Dt\De=V(0,\Dt\simod)\subset\De$ 
is the union of faces of $\De$ which do not contain the sector $V(0,\taumod)$.
Note that $\Dt\De\cap V(0,\taumod)=\D V(0,\taumod)=V(0,\D\taumod)$. 

\begin{dfn}
A sequence $(\de_n)$ in $\De$ is 

(i) {\em $\taumod$-regular} if it drifts away from $\Dt\De$,
$$ d(\de_n,\Dt\De) \to+\infty .$$

(ii) {\em $\taumod$-pure} if it is contained in a tubular neighborhood of the sector $V(0,\taumod)$
and drifts away from its boundary,
$$ d(\de_n,\D V(0,\taumod)) \to+\infty .$$
\end{dfn}
Note that $(\de_n)$ is $\taumod$-regular/pure iff $(\iota\de_n)$ is $\iota\taumod$-regular/pure.

We extend these notions to sequences in $X$ and $G$:
\begin{dfn}[Regular and pure]
\label{def:pureg}
(i) A sequence $(x_n)$ in $X$ is {\em $\taumod$-regular}, respectively, {\em $\taumod$-pure}
if for some (any) base point $o\in X$ the sequence of $\De$-distances $d_{\De}(o,x_n)$ in $\De$
has this property.

(ii) A sequence $(g_n)$ in $G$ is {\em $\taumod$-regular}, respectively, {\em $\taumod$-pure}
if for some (any) point $x\in X$ the orbit sequence $(g_nx)$ in $X$ has this property.

(iii) A subgroup $\Ga<G$ is {\em $\taumod$-regular}
if all sequences of distinct elements in $\Ga$ have this property.
\end{dfn}

That these properties 
are independent of the base point and stable under bounded perturbation of the sequences,
is due to the triangle inequality
$|d_{\De}(x,y)-d_{\De}(x',y')|\leq d(x,x')+d(y,y')$.

Subsequences of $\taumod$-regular/pure sequences are again $\taumod$-regular/pure.

Clearly, $\taumod$-pureness is a strengthening of $\taumod$-regularity;
a sequence in $\De$ is $\taumod$-pure iff it is $\taumod$-regular and contained in a tubular neighborhood of $V(0,\taumod)$.

The face type of a pure sequence is uniquely determined.
Moreover,
a $\taumod$-regular sequence is $\taumod'$-regular for every face type $\taumod'\subset\taumod$,
because $\Dtp\De\subset\Dt\De$.

A sequence $(g_n)$ is $\taumod$-regular/pure iff 
the inverse sequence $(g_n^{-1})$ is $\iota\taumod$-regular/pure,
because $d_{\De}(x,g_n^{-1}x)=d_{\De}(g_nx,x)=\iota d_{\De}(x,g_nx)$.

Note that $\taumod$-regular subgroups are in particular {\em discrete}.
If $\rank(X)=1$, then discreteness is equivalent to ($\simod$-)regularity.
In higher rank, {\em regularity} can be considered as a {\em strengthening of discreteness}:
A discrete subgroup $\Ga<G$ may not be $\taumod$-regular for any face type $\taumod$;
this can happen e.g.\ for free abelian subgroups of transvections of rank $\geq2$.

A property for sequences in $G$ equivalent to regularity 
had appeared in \cite[Lemma 3.5(1)]{Benoist}.

\begin{lem}[Pure subsequences]
\label{lem:obspureg}
Every sequence, which diverges to infinity, 
contains a $\taumod$-pure subsequence for some face type $\taumod\subseteq\simod$.
\end{lem}
\proof
In the case of sequences in $\De$,
take $\taumod$ to be a minimal face type 
so that a subsequence is contained in a tubular neighborhood of $V(0,\taumod)$.
\qed

\medskip
Note also that a sequence, which diverges to infinity, 
is $\taumod$-regular iff it contains $\numod$-pure subsequences only for face types $\numod\supseteq\taumod$. 

\medskip
The lemma implies in particular,
that every sequence $\ga_n\to\infty$ in a discrete subgroup $\Ga<G$ contains a subsequence 
which is $\taumod$-regular, even $\taumod$-pure, for some face type $\taumod$.

\begin{rem}
\label{rem:regfins}
Regularity has a natural Finsler geometric interpretation, cf.\ \cite{bordif}:
A sequence in $X$ is $\taumod$-regular iff, 
in the Finsler compactification $\ol X^{Fins}=X\sqcup\geo^{Fins}X$ of $X$, 
it accumulates at the closure of the stratum $S_{\taumod}\subset\geo^{Fins}X$ at infinity.
\end{rem}

\subsection{Contraction implies regularity}
\label{sec:contrimplrg}

In this section and the next,
we relate contractivity and regularity for sequences and, as a consequence, 
establish the equivalence between $\taumod$-regularity and the $\taumod$-convergence property for discrete subgroups.

To relate contraction and regularity,
it is useful to consider the $G$-action on flats.
We recall that ${\mathcal F}_{\taumod}$ denotes the space of flats $f\subset X$ of type $\taumod$
(see section~\ref{sec:symmbas}).
Two flats $f_{\pm}\in{\mathcal F}_{\taumod}$ are {\em dynamically related} 
with respect to a sequence $(g_n)$ in $G$,
$$f_-\stackrel{(g_n)}{\sim}f_+ ,$$
if there exists a sequence of flats $f_n\to f_-$ in ${\mathcal F}_{\taumod}$
such that $g_nf_n\to f_+$.
The action of $(g_n)$ on ${\mathcal F}_{\taumod}$ is {\em proper} 
iff there are no dynamical relations with respect to subsequences,
cf.\ section~\ref{sec:discdynrl}.

Dynamical relations between singular flats yield dynamical relations between maximal ones:
\begin{lem}
\label{lem:dynrsingmx}
If $f_{\pm}\in{\mathcal F}_{\taumod}$ are flats such that $f_-\stackrel{(g_n)}{\sim}f_+$,
then for every maximal flat $F_+\supseteq f_+$ there exist a maximal flat $F_-\supseteq f_-$
and a subsequence $(g_{n_k})$ 
such that $F_-\stackrel{(g_{n_k})}{\sim}F_+$.
\end{lem}
\proof
Let $f_n\to f_-$ be a sequence in ${\mathcal F}_{\taumod}$
such that $g_nf_n\to f_+$.
Then there exists a sequence of maximal flats $F_n\supseteq f_n$ such that $g_nF_n\to F_+$.
The sequence $(F_n)$ is bounded 
because the sequence $(f_n)$ is,
and hence $(F_n)$ subconverges to a maximal flat $F_-\supseteq f_-$.
\qed

\medskip
For pure sequences there are dynamical relations between singular flats of the corresponding type 
with respect to suitable subsequences:
\begin{lem}
\label{lem:improp}
If $(g_n)$ is $\taumod$-pure,
then the action of $(g_n)$ on ${\mathcal F}_{\taumod}$ is not proper.

More precisely,
there exist simplices $\tau_{\pm}\in\Flagt$ such that 
for every flat $f_+\in{\mathcal F}_{\taumod}$ asymptotic to $\tau_+$
there exist a flat $f_-\in{\mathcal F}_{\taumod}$ asymptotic to $\tau_-$
and a subsequence $(g_{n_k})$ such that 
$$f_-\stackrel{(g_{n_k})}{\sim}f_+.$$

\end{lem}
\proof
By pureness, there exists a sequence $(\tau_n)$ in $\Flagt$
such that 
\begin{equation}
\label{eq:clsct}
\sup_n d(g_nx,V(x,\tau_n))  <+\infty 
\end{equation}
for any point $x\in X$. 
There exists a subsequence $(g_{n_k})$ such that 
$\tau_{n_k}\to\tau_+$ and $g_{n_k}^{-1}\tau_{n_k}\to\tau_-$.

Let $f_+\in{\mathcal F}_{\taumod}$ be asymptotic to $\tau_+$.
We choose $x\in f_+$ and consider the sequence of flats $f_k\in{\mathcal F}_{\taumod}$ 
through $x$ asymptotic to $\tau_{n_k}$.
Then $f_k\to f_+$.
The sequence of flats $(g_{n_k}^{-1}f_k)$ is bounded
as a consequence of (\ref{eq:clsct}).
Therefore, after further extraction, we obtain convergence $g_{n_k}^{-1}f_k\to f_-$.
The limit flat $f_-$ is asymptotic to $\tau_-$
because the $f_k$ are asymptotic to $g_{n_k}^{-1}\tau_{n_k}$.
\qed

\medskip
By a diagonal argument one can also show that the subsequences $(g_{n_k})$ 
in the two previous lemmas 
can be made independent of the flats $F_+$ respectively $f_+$.

For contracting sequences, 
the possible dynamical relations between maximal flats are restricted as follows:
\begin{lem}
\label{lem:contrdynrel}
Suppose that $(g_n)$ is {\em $\taumod$-contracting} with (\ref{eq:contrtau}),
and that 
$F_-\stackrel{(g_n)}{\sim}F_+$
for maximal flats $F_{\pm}\in{\mathcal F}$. 
Then $\tau_{\pm}\subset\geo F_{\pm}$.
\end{lem}
\proof
Suppose that $\tau_-\not\subset\geo F_-$.
Then the visual boundary sphere $\geo F_-$ contains 
at least two different simplices $\hat\tau_-,\hat\tau'_-$ opposite to $\tau_-$,
cf.\ Lemma~\ref{lem:oneantip}.

Let $F_n\to F_-$ be a sequence in ${\mathcal F}$ such that $g_nF_n\to F_+$.
Due to $F_n\to F_-$,
there exist sequences of simplices $\tau_n,\tau'_n\subset\geo F_n$ 
such that $\tau_n\to\hat\tau_-$ and $\tau'_n\to\hat\tau'_-$.
In particular, $\tau_n\neq\tau'_n$ for large $n$.
After extraction,
we also obtain convergence 
$g_n\tau_n\to\hat\tau_+$ and $g_n\tau'_n\to\hat\tau'_+$.
Moreover, 
since $g_nF_n\to F_+$,
it follows that the limits $\hat\tau_+,\hat\tau'_+$
are {\em different} simplices in $\geo F_+$. 

This is however in conflict with the contraction property (\ref{eq:contrtau}).
In view of $\hat\tau_-,\hat\tau'_-\in C(\tau_-)$, 
the latter implies that 
$g_n\tau_n\to\tau_+$ and $g_n\tau'_n\to\tau_+$,
convergence to the same simplex, 
a contradiction.
Thus, $\tau_-\subset\geo F_-$.

Considering the inverse sequence $(g_n^{-1})$
yields that also $\tau_+\subset\geo F_+$,
cf.\ Lemma~\ref{lem:contrsym}.
\qed

\medskip
Combining the previous lemmas, we obtain:
\begin{lem}
\label{lem:contrpuresq}
If a sequence in $G$ is $\taumod$-contracting and $\numod$-pure, then $\taumod\subseteq\numod$.
\end{lem}
\proof
We denote the sequence by $(g_n)$ and assume (\ref{eq:contrtau}).
According to Lemmas~\ref{lem:improp} and~\ref{lem:dynrsingmx},
by $\numod$-purity, there exist simplices $\nu_{\pm}\in\Flagn$
such that for {\em every} maximal flat $F_+$ 
with $\geo F_+\supset\nu_+$ 
there exist a maximal flat $F_-$ 
with $\geo F_-\supset\nu_-$ 
and a subsequence $(g_{n_k})$ such that 
$$F_-\stackrel{(g_{n_k})}{\sim}F_+.$$
By Lemma~\ref{lem:contrdynrel},
always $\tau_+\subset\geo F_+$.
Varying $F_+$, it follows that $\tau_+\subseteq\nu_+$, 
cf.\ Lemma~\ref{lem:intapts}.
\qed

\medskip
From these observations, we conclude:
\begin{prop}[Contracting implies regular]
\label{prop:contrimpreg}
If a sequence in $G$ is $\taumod$-contracting, then it is $\taumod$-regular.
\end{prop}
\proof
Consider a sequence in $G$ which is not $\taumod$-regular.
Then a subsequence is $\numod$-pure for some face type 
$\numod\subseteq\Dt\simod$, 
compare Lemma~\ref{lem:obspureg}.
The condition on the face type is equivalent to $\numod\not\supseteq\taumod$.
By the last lemma, the subsequence cannot be $\taumod$-contracting.
\qed

\subsection{Regularity implies contraction}\label{sec:reg-con}

We now prove a converse to Proposition~\ref{prop:contrimpreg}.
Since contractivity involves a convergence condition,
we can expect regular sequences to be contracting only after extraction. 

Consider a $\taumod$-regular sequence $(g_n)$ in $G$.
After fixing a point $x\in X$,
there exist simplices $\tau_n^{\pm}\in\Flagpmt$
(unique for large $n$)
such that 
\begin{equation}
\label{eq:shad}
g_n^{\pm1}x \in V(x,\st(\tau_n^{\pm})) .
\end{equation}
Note that the sequence $(g_n^{-1})$ is $\iota\taumod$-regular,
compare the comment after Definition~\ref{def:pureg}.
\begin{lem}
\label{lem:flconvcontr}
If $\tau_n^{\pm}\to\tau_{\pm}$ in $\Flagpmt$,
then $(g_n)$ is $\taumod$-contracting with (\ref{eq:contrtau}).
\end{lem}
\proof
Since $x\in g_nV(x,\st(\tau_n^-))=V(g_nx,\st(g_n\tau_n^-))$,
it follows together with $g_nx \in V(x,\st(\tau_n^+))$
that the Weyl cones $V(g_nx,\st(g_n\tau_n^-))$ and $V(x,\st(\tau_n^+))$ 
lie in the same parallel set, namely in $P(g_n\tau_n^-,\tau_n^+)$, and face in opposite directions. 
In particular, the simplices $g_n\tau_n^-$ and $\tau_n^+$ are $x$-opposite,
and thus $g_n\tau_n^-$ converges to the simplex $\hat\tau_+$ which is $x$-opposite to $\tau_+$,
$$ g_n\tau_n^- \to\hat\tau_+ .$$
Since the sequence $(g_n^{-1}x)$ is $\iota\taumod$-regular,
it holds that 
$$ d(g_n^{-1}x,\D V(x,\st(\tau_n^-))) \to+\infty  .$$
By Lemma~\ref{lem:expconvsect},
for any $r,R>0$,
one has for $n\geq n(r,R)$ 
the inclusion of shadows 
(cf.\ (\ref{eq:shadwdf}))
$$ U_{\tau_n^-,x,R} \subset U_{\tau_n^-,g_n^{-1}x,r} .$$
Consequently, there exist sequences of positive numbers $R_n\to+\infty$ and $r_n\to0$ such that 
$$ U_{\tau_n^-,x,R_n} \subset U_{\tau_n^-,g_n^{-1}x,r_n} $$
for large $n$,
equivalently
\begin{equation}
\label{eq:mpfshd}
g_nU_{\tau_n^-,x,R_n} \subset U_{g_n\tau_n^-,x,r_n} .
\end{equation}
Since $\tau_n^-\to\tau_-$ and $R_n\to+\infty$,
the shadows $U_{\tau_n^-,x,R_n}\subset C(\tau_n^-)\subset\Flagt$ {\em exhaust} $C(\tau_-)$
in the sense that every compact in $C(\tau_-)$ is contained in $U_{\tau_n^-,x,R_n}$ for large $n$.\footnote{Indeed, 
for fixed $R>0$ we have Hausdorff convergence 
$U_{\tau_n^-,x,R}\to U_{\tau_-,x,R}$ in $\Flagt$,
which follows e.g.\ from the transitivity of the action $K_x\acts\Flagit$
of the maximal compact subgroup $K_x<G$ fixing $x$.
Furthermore,
the shadows $U_{\tau_-,x,R}$ exhaust $C(\tau_-)$ as $R\to+\infty$,
cf.\ the continuity part of Lemma~\ref{lem:contpr}.}
On the other hand,
since $g_n\tau_n^-\to\hat\tau_+$ and $r_n\to0$,
the 
$U_{g_n\tau_n^-,x,r_n}$ {\em shrink}, i.e.\ Hausdorff converge to the point $\tau_+$.\footnote{Indeed, 
$U_{g_n\tau_n^-,x,r}\to U_{\hat\tau_+,x,r}$ in $\Flagt$ for fixed $r>0$,
and $U_{\hat\tau_+,x,r}\to\tau_+$ as $r\to0$,
using again the continuity part of Lemma~\ref{lem:contpr}
and the fact that the function (\ref{eq:distfrpar}) assumes the value zero only in $\tau_+$.}
Therefore,
(\ref{eq:mpfshd}) implies that 
$$ g_n|_{C(\tau_-)} \to \tau_+ $$
uniformly on compacts,
i.e.\ $(g_n)$ is $\taumod$-contracting. 
\qed

\medskip
With the lemma, we can add the desired converse to Proposition~\ref{prop:contrimpreg}
and obtain a characterization of regularity in terms of contraction:

\begin{prop}
\label{prop:regequivcontr}
The following properties are equivalent for sequences in $G$:

(i) Every subsequence contains a $\taumod$-contracting subsequence. 

(ii) The sequence is $\taumod$-regular.
\end{prop}

\proof
This is a direct consequence of the lemma. 
For the implication (ii)$\Ra$(i) one uses the compactness of flag manifolds.
The implication (i)$\Ra$(ii) is obtained as follows, 
compare the proof of Proposition~\ref{prop:contrimpreg}:
If a sequence is not $\taumod$-regular, then it contains a $\numod$-pure subsequence
for some face type $\numod\not\supseteq\taumod$.
Every subsequence of this subsequence is again $\numod$-pure
and hence not $\taumod$-contracting by Lemma~\ref{lem:contrpuresq}.
\qed

\medskip
A version of Propostition \ref{prop:regequivcontr}  had already been proven by Benoist  in \cite[Lemma 3.5]{Benoist}.

\medskip
We conclude for subgroups:
\begin{thm}
\label{thm:regimplcontrgp}
A subgroup $\Ga<G$ is $\taumod$-regular iff it is a $\taumod$-convergence subgroup.
\end{thm}
\proof
By definition, 
$\Ga$ is $\taumod$-regular iff every sequence $(\ga_n)$ of distinct elements in $\Ga$ is $\taumod$-regular,
and $\taumod$-convergence iff every such sequence $(\ga_n)$ has a $\taumod$-contracting subsequence.
According to the proposition, both conditions are equivalent.
\qed

\subsection{Convergence at infinity and limit sets}
\label{sec:conv}

The discussion in the preceding two sections 
leads to a natural notion of convergence at infinity 
for regular sequences in $X$ and $G$. 
As regularity, 
it can be expressed both in terms of orbit geometry in $X$ and dynamics on flag manifolds.

We first consider a $\taumod$-regular sequence $(g_n)$ in $G$.
Flexibilizing condition (\ref{eq:shad}),
we choose points $x,x'\in X$ and consider a sequence $(\tau_n)$ in $\Flagt$ 
such that 
\begin{equation}
\label{eq:bdddstfrcona}
\sup_nd\bigl(g_nx,V(x',\st(\tau_n))\bigr) < +\infty.
\end{equation}
Note that the condition is independent of the choice of the points $x$ and $x'$.\footnote{Recall that 
the Hausdorff distance of asymptotic Weyl cones $V(y,\st(\tau))$ and $V(y',\st(\tau))$
is bounded by the distance $d(y,y')$ of their tips.}

\begin{lem}
\label{lem:regseqindepch}
The accumulation set of $(\tau_n)$ in $\Flagt$ depends only on $(g_n)$.
\end{lem}
\proof
Let $(\tau'_n)$ be another sequence in $\Flagt$                                                                                                               
such that $d(g_nx,V(x',\st(\tau'_n)))$ is uniformly bounded.
Assume that after extraction $\tau_n\to\tau$ and $\tau'_n\to\tau'$.
We must show that $\tau=\tau'$.

We may suppose that $x'=x$.
There exist bounded sequences $(b_n)$ and $(b'_n)$ in $G$
such that
$$g_nb_nx\in V(x,\st(\tau_n))                                                                                              
\qquad\hbox{ and }\qquad                                                                                                   
g_nb'_nx\in V(x,\st(\tau'_n))$$
for all $n$.
Note that the sequences $(g_nb_n)$ and $(g_nb'_n)$ in $G$ are again $\taumod$-regular.
By Lemma~\ref{lem:flconvcontr},
after further extraction,
they are $\taumod$-contracting with
$$g_nb_n|_{C(\tau_-)}\to\tau                                                                                               
\qquad\hbox{ and }\qquad                                                                                                   
g_nb'_n|_{C(\tau'_-)}\to\tau'$$
uniformly on compacts for some $\tau_-,\tau'_-\in\Flagit$.
Moreover, we may assume convergence $b_n\to b$ and $b'_n\to b'$.
Then
$$g_n|_{C(b\tau_-)}\to\tau
\qquad\hbox{ and }\qquad
g_n|_{C(b'\tau'_-)}\to\tau'$$
uniformly on compacts.
With Lemma~\ref{lem:contruniq} it follows that $\tau=\tau'$.
\qed

\medskip
In view of the lemma, we can define the following notion of convergence:

\begin{dfn}[Flag convergence of sequences in $G$]
\label{def:flgcnvg}
A $\taumod$-regular sequence $(g_n)$ in $G$ 
{\em $\taumod$-flag converges} 
to a simplex $\tau\in\Flagt$,
$$ g_n\to\tau ,$$
if $\tau_n\to\tau$ in $\Flagt$
for some sequence $(\tau_n)$ in $\Flagt$ 
satisfying (\ref{eq:bdddstfrcona}).
\end{dfn}

We can now characterize contraction in terms of flag convergence.
We rephrase Lemma~\ref{lem:flconvcontr} and show that its converse holds as well:

\begin{lem}
\label{lem:flconvcontrconvreph}
For a sequence $(g_n)$ in $G$ and simplices $\tau_{\pm}\in\Flagpmt$,
the following are equivalent:

(i) $(g_n)$ is $\taumod$-contracting 
with 
$g_n|_{C(\tau_-)}\to\tau_+$ uniformly on compacts.

(ii) $(g_n)$ is $\taumod$-regular and $g_n^{\pm1}\to\tau_{\pm}$.
\end{lem}
In part (ii), the sequence $(g_n^{-1})$ is $\iota\taumod$-regular 
and $g_n^{-1}\to\tau_-$ means $\iota\taumod$-flag convergence. 
\proof
The implication (ii)$\Ra$(i) is Lemma~\ref{lem:flconvcontr}.

Conversely,
suppose that (i) holds.
Since the sequence $(g_n)$ is $\taumod$-contracting,
it is $\taumod$-regular by Proposition~\ref{prop:contrimpreg}.
Let $(\tau_n^{\pm})$ be sequences satisfying (\ref{eq:shad}).
We must show that $\tau_n^{\pm}\to\tau_{\pm}$.
Otherwise, after extraction we obtain that
$\tau_n^{\pm}\to\tau'_{\pm}$
with $\tau'_+\neq\tau_+$ or $\tau'_-\neq\tau_-$.
Then also
$g_n|_{C(\tau'_-)}\to\tau'_+$
by Lemma~\ref{lem:flconvcontr},
and Lemma~\ref{lem:contruniq} implies that $\tau'_{\pm}=\tau_{\pm}$,
a contradiction.
\qed

\medskip
Vice versa,
we can characterize flag convergence in terms of contraction
and thus give an alternative {\em dynamical definition} of it:

\begin{lem}
\label{lem:flconvitcontr}
For a sequence $(g_n)$ in $G$, 
the following are equivalent:

(i) $(g_n)$ is $\taumod$-regular and $g_n\to\tau$.

(ii) There exists a bounded sequence $(b_n)$ in $G$ and $\tau_-\in\Flagit$
such that 
$g_nb_n|_{C(\tau_-)}\to\tau$ uniformly on compacts.

(iii) There exists a bounded sequence $(b'_n)$ in $G$ 
such that 
$b'_ng_n^{-1}|_{C(\tau)}$ converges to a constant map uniformly on compacts.
\end{lem}
\proof
(ii)$\Ra$(i):
According to the previous lemma
the sequence
$(g_nb_n)$ is $\taumod$-regular and $\taumod$-flag converges,
$g_nb_n\to\tau$.
Since $d(g_nx,g_nb_nx)$ is uniformly bounded,
this is equivalent to $(g_n)$ being $\taumod$-regular and $g_n\to\tau$.

(i)$\Ra$(ii):
The sequence $(g^{-1}_n)$ is $\iota\taumod$-regular.
There exists a bounded sequence $(b'_n)$ in $G$ such that $(b'_ng^{-1}_n)$
$\iota\taumod$-flag converges,
$b'_ng^{-1}_n\to\tau_-\in\Flagit$.
We put $b_n={b'_n}^{-1}$.
Since also $(g_nb_n)$ is $\taumod$-regular and $g_nb_n\to\tau$,
it follows
from the previous lemma
that $g_nb_n|_{C(\tau_-)}\to\tau$ uniformly on compacts.

The equivalence (ii)$\Leftrightarrow$(iii) with $b'_n=b_n^{-1}$ follows from Lemma~\ref{lem:contrsym}.
\qed

\medskip
We 
carry over 
the notion of flag convergence to sequences in $X$.

Consider now a $\taumod$-regular sequence $(x_n)$ in $X$. 
We choose again a base point $x\in X$ and consider a sequence $(\tau_n)$ in $\Flagt$ such that 
\begin{equation}
\label{eq:bdddstfrcon}
\sup_nd\bigl(x_n,V(x,\st(\tau_n))\bigr) < +\infty ,
\end{equation}
analogous to (\ref{eq:bdddstfrcona}).
As before, the condition is independent of the choice of the point $x$, 
and we obtain a version of Lemma~\ref{lem:regseqindepch}:

\begin{lem}
\label{lem:regseqindepchx}
The accumulation set of $(\tau_n)$ in $\Flagt$ depends only on $(x_n)$.
\end{lem}
\proof
Let $(g_n)$ be a sequence in $G$ such that the sequence $(g_n^{-1}x_n)$ in $X$ is bounded.
Then $(g_n)$ is $\taumod$-regular 
and (\ref{eq:bdddstfrcon}) becomes equivalent to (\ref{eq:bdddstfrcona}).
This reduces the claim to Lemma~\ref{lem:regseqindepch}.
\qed

\medskip
We therefore can define, analogous to Definition~\ref{def:flgcnvg} above:
\begin{dfn}[Flag convergence of sequences in $X$]
A $\taumod$-regular sequence $(x_n)$ in $X$ 
{\em $\taumod$-flag converges} 
to a simplex $\tau\in\Flagt$,
$$ x_n\to\tau ,$$
if $\tau_n\to\tau$ in $\Flagt$ for some sequence $(\tau_n)$ in $\Flagt$ satisfying (\ref{eq:bdddstfrcon}).
\end{dfn}

For any $\taumod$-regular sequence $(g_n)$ in $G$ and any point $x\in X$,
we have 
$g_n\to\tau$ iff $g_nx\to\tau$.

Flag convergence and flag limits are stable under bounded perturbations of sequences:
\begin{lem}
(i)
For any $\taumod$-regular sequence $(g_n)$ and any bounded sequence $(b_n)$ in $G$,
the sequences $(g_n)$ and $(g_nb_n)$ have the same $\taumod$-flag accumulation sets in $\Flagt$.

(ii)
If $(x_n)$ and $(x'_n)$ are $\taumod$-regular sequences in $X$ such that $d(x_n,x'_n)$ is uniformly bounded,
then both sequences have the same $\taumod$-flag accumulation set in $\Flagt$.
\end{lem}
\proof
(i)
The sequence $(g_nb_n)$ is also $\taumod$-regular and satisfies condition (\ref{eq:bdddstfrcona}) iff $(g_n)$ does.

(ii)
The sequence $(x'_n)$ satisfies condition (\ref{eq:bdddstfrcon}) iff $(x'_n)$ does. 
\qed

\begin{rem}
\label{rem:flconvfins}
There is a natural topology on the bordification $X\sqcup\Flagt$
which induces $\taumod$-flag convergence.
Moreover,
the bordification embeds into a natural Finsler compactification of $X$,
compare Remark~\ref{rem:regfins}.
\end{rem}

Flag convergence leads to a notion of limit sets in flag manifolds
for subgroups:

\begin{dfn}[Flag limit set]
\label{def:flagl}
For a subgroup $\Ga<G$, the {\em $\taumod$-limit set} 
$$\LatGa\subset\Flagt$$
is the set of possible limit simplices of $\taumod$-flag converging $\taumod$-regular sequences in $\Ga$,
equivalently,
the set of simplices $\tau_+$ as in (\ref{eq:contrtau})
for all $\taumod$-contracting sequences in $\Ga$.
\end{dfn}

The limit set is $\Ga$-invariant and closed,
as a diagonal argument shows.

\begin{rem}
Benoist introduced in \cite[\S 3.6]{Benoist}
a notion of limit set $\La_{\Ga}$ for Zariski dense subgroups $\Ga$ of reductive algebraic groups over local fields
which in the case of real semisimple Lie groups
is equivalent to (the dynamical version of) our concept of $\simod$-limit set $\Las$.\footnote{Benoist's limit set $\La_{\Ga}$ 
is contained in the flag manifold $Y_{\Ga}$
which in the case of real Lie groups is the full flag manifold $G/B$, see the beginning of \S 3 of his paper.
It consists of the limit points of sequences contracting on $G/B$, cf.\ his Definitions 3.5 and 3.6.}
What we call the $\taumod$-limit set $\Lat$ for other face types $\taumod\subsetneq\simod$
is mentioned in his Remark 3.6(3),
and his work implies that, in the Zariski dense case,
$\Lat$ is the image of $\Las$ under the natural projection $\Flags\to\Flagt$ of flag manifolds.
\end{rem}

\subsection{Uniform regularity}

In this section we introduce stronger forms of the regularity conditions discussed in section~\ref{sec:reg}.

We first consider sequences in the euclidean model Weyl chamber $\De$.
\begin{dfn}
\label{def:unifreg}
A sequence $\de_n\to\infty$ in $\De$ is 
{\em uniformly $\taumod$-regular} if it drifts away from $\Dt\De$
at a linear rate with respect to 
its norm,
$$ \liminf_{n\to+\infty} \frac{d(\de_n,\Dt\De)}{\|\de_n\|} > 0.$$
\end{dfn}

We extend these notions to sequences in $X$ and $G$,
compare Definition~\ref{def:pureg}:
\begin{dfn}[Uniformly regular]
(i) A sequence $(x_n)$ in $X$ is {\em uniformly $\taumod$-regular}
if for some (any) base point $o\in X$ the sequence of $\De$-distances $d_{\De}(o,x_n)$ in $\De$
has this property.

(ii) A sequence $(g_n)$ in $G$ is {\em uniformly $\taumod$-regular}
if for some (any) point $x\in X$ the orbit sequence $(g_nx)$ in $X$ has this property.

(iii) A subgroup $\Ga<G$ is {\em uniformly $\taumod$-regular}
if all sequences of distinct elements in $\Ga$ have this property.
\end{dfn}

For a subgroup $\Ga<G$, uniform $\taumod$-regularity is equivalent to 
the visual limit set $\LaGa\subset\geo X$ being contained in the union of the open $\taumod$-stars.

\section{Asymptotic and coarse properties of discrete subgroups}
\label{sec:prp}

This chapter is the core of the paper.
In section~\ref{sec:bdemb},
motivated by the boundary map part of the original Anosov notion,
we study equivariant embeddings of the Gromov boundaries of word hyperbolic subgroups into flag manifolds.
We show how these boundary embeddings can be used, especially for regular subgroups,
to control the geometry of the orbits in the symmetric space:
Intrinsic geodesic lines in the group are uniformly close to parallel sets in the symmetric space.
Moreover, in the generic case, for instance for Zariski dense subgroups,
intrinsic rays in the group are close to Weyl cones.
This conicality property implies in particular
that the boundary map continuously extends the orbit maps to infinity and identifies the Gromov boundary with the limit set.
This leads us to notion of asymptotically embedded subgroups
discussed in section~\ref{sec:asyemb}.
We find that asymptotic embeddedness has strong implications for the coarse extrinsic geometry of subgroups:
They are undistorted,
and moreover their intrinsic geodesics satisfy a higher rank version of the ``Morse property'';
they are uniformly close to diamonds. 
This motivates the notion of Morse subgroups studied in section~\ref{sec:morse}. 
The higher rank Morse property immediately implies that the limit set is conical and antipodal.
We call regular subgroups with the latter properties RCA and study them in section~\ref{sec:conic}. 
Using Bowditch's dynamical characterization of hyperbolic groups,
we show that RCA subgroups are asymptotically embedded,
closing part of the circle.
In section~\ref{sec:expa}, we observe that conicality implies expansive dynamics at the limit set,
which yields another equivalent property for subgroups,
this time formulated purely in terms of the dynamics on flag manifolds.
In sections~\ref{sec:anosov} and~\ref{sec:onos}, 
we discuss different (uniform and non-uniform) versions of our Anosov condition
and show that it is equivalent to the previous conditions as well as to the original definition of Anosov subgroups.
In section~\ref{sec:mqg} we take up the discussion of the Morse property.
Leaving the context of discrete subgroups, 
we study the geometry of Morse quasigeodesics in symmetric spaces.
We characterize them as bounded perturbations of Finsler quasigeodesics
and study the behavior of the $\De$-distance along them:
we prove that via the $\De$-distance they project to Morse quasigeodesics in $\De$. 
We also obtain another characterization of Morse subgroups by the quasiconvexity property
that their intrinsic geodesics are extrinsically Morse quasigeodesics,
equivalently, are uniformly close to Finsler geodesics.

\subsection{Antipodality}
\label{sec:antip}

If $X$ has rank one,
then $G$ acts transitively on pairs of distinct points in $\geo X$.
Thus there are only two possibilities for the {\em relative position} of two points in the visual boundary:
They can coincide or be different.
In higher rank,
the $G$-actions on the associated flag manifolds are in general not two point transitive
and there are more possibilities for the relative position.

We recall (see section~\ref{sec:symmbas}) that two simplices $\tau,\tau'\subset\geo X$ 
are called {\em opposite} or {\em antipodal}
if they are opposite simplices in the apartments $a\subset\geo X$ containing them both. 
Their types are then related by $\theta(\tau')=\iota\theta(\tau)$.
In particular, if three simplices are pairwise opposite,
their types must be equal and $\iota$-invariant. 

\begin{dfn}[Antipodal]
Suppose that $\taumod$ is $\iota$-invariant.

(i) 
A subset of $\Flagt$ is {\em antipodal}
if it consists of pairwise opposite simplices.

(ii)
A map into $\Flagt$ is {\em antipodal} if it sends different elements to opposite simplices.

(iii)
A subgroup $\Ga<G$ is {\em $\taumod$-antipodal}
if $\LatGa$ is antipodal.
\end{dfn}

Being antipodal is an {\em open} condition for pairs of points in flag manifolds.
It is the {\em generic} relative position.
Antipodal maps are in particular {\em injective}.

We note that for a $\taumod$-antipodal $\taumod$-convergence subgroup $\Ga<G$ 
the action $$\Ga\acts\LatGa$$ has {\em convergence dynamics} in the usual sense,
see section~\ref{sec:convdy}:
If $(\ga_n)$ is a sequence in $\Ga$ such that 
$\ga_n|_{C(\tau_-)}\to\tau_+$,
then $\tau_{\pm}\in\LatGa$.
Due to antipodality, $\LatGa-\{\tau_-\}\subset C(\tau_-)$ and we obtain the intrinsic convergence property.

\subsection{Boundary embeddings and limit sets}
\label{sec:bdemb}

In this section,
we study embeddings of word hyperbolic groups into semisimple Lie groups
which admit a certain kind of continuous boundary map.
We will assume that $\taumod$ is $\iota$-invariant.

\begin{dfn}[Boundary embedded]
\label{def:bdemb}
A subgroup $\Ga<G$ is {\em $\taumod$-boundary embedded}
if it is intrinsically word hyperbolic 
and there exists an antipodal $\Ga$-equivariant continuous embedding 
\begin{equation}
\label{eq:bdemb}
\beta:\geo\Ga\to\Flagt 
\end{equation}
of the Gromov boundary $\geo\Ga$ of $\Ga$. The map $\beta$ is 
called a {\em boundary embedding}.
If $|\geo\Ga|\leq2$,
we require in addition that $\Ga$ is discrete in $G$.
\end{dfn}

Thus, $\taumod$-boundary embedded subgroups are necessarily {\em discrete}, since 
$\Ga$ acts on $\beta(\geo\Ga)$ as a discrete convergence group if $|\geo\Ga|\geq3$.\footnote{Note that 
boundary embedded subgroups are not required to be regular, although they frequently are, see Theorem 3.11 in \cite{manicures}.} 

Boundary embeddings are in general not unique.
This is so by trivial reasons if $|\geo\Ga|=2$, cf.\ below,
but it also happens if $|\geo\Ga|\geq3$, see \cite[Example 6.20]{morse}.

\medskip
In order to understand the implications 
of a boundary embedding,
we will first use it 
to obtain
control on the geometry of the $\Ga$-orbits in $X$.

We fix a word metric on $\Ga$. 
Via the antipodal boundary embedding $\beta$
one can assign to every discrete geodesic {\em line}\footnote{Recall that by a {discrete geodesic line},
we mean an isometric embedding of $\Z$, cf, section \ref{sec:general}.} 
$l:\Z\to\Ga$
a parallel set in $X$. 
Namely, let $\zeta_{\pm}:=l(\pm\infty)\in\geo\Ga$ denote the ideal endpoints of the line.
Their image simplices $\beta(\zeta_{\pm})\in\Flagt$ are opposite
and determine the parallel set 
$$P(\beta(\zeta_-),\beta(\zeta_+))\subset X .$$
We consider the images of the discrete geodesic lines $l$ in $\Ga$ under the orbit map 
$o_x:\Ga\to\Ga x\subset X$ for a point $x\in X$ (fixed throughout the discussion)
and claim that the discrete paths $lx:\Z\to X$ are uniformly close to the corresponding parallel sets:\footnote{
For a map $\phi:N\to\Ga$ and a point $x\in X$ we denote by $\phi x:N\to X$ 
the map sending $n\in N$ to $\phi(n)x\in X$.}
\begin{lem}[Lines go close to parallel sets]
\label{lem:lnclpar}
The discrete path 
$lx$ is contained in a tubular neighborhood of the parallel set $P(\beta(\zeta_-),\beta(\zeta_+))$
with uniform radius 
$\rho=\rho(\Ga,x)$.
\end{lem}
Here and below,
we mean by the dependence of a constant on $\Ga$ that it depends on $\Ga$ as a subgroup of $G$ 
and also on the chosen word metric on $\Ga$.
\proof
This can be seen by a simple compactness argument: 
Let 
\begin{equation}
\label{eq:oppfl}
(\Flagt\times\Flagt)^{opp}
\subset\Flagt\times\Flagt
\end{equation}
denote the subspace of pairs of opposite simplices.
It is the open and dense $G$-orbit 
and in particular a homogeneous $G$-space.
The latter implies that the function 
on $(\Flagt\times\Flagt)^{opp}\times X$
assigning 
\begin{equation}
\label{eq:dstcnt}
(\tau_-,\tau_+,x') \mapsto d\bigl(x',P(\tau_-,\tau_+)\bigr)
\end{equation}
is continuous, because 
$d(gx',P(h\tau_-,h\tau_+))=d(h^{-1}gx',P(\tau_-,\tau_+))$ for $g,h\in G$.
Also the map
$$ {\mathcal L} \to (\Flagt\times\Flagt)^{opp}\times X$$
from the space ${\mathcal L}$ of discrete geodesic lines $l:\Z\to\Ga$ 
\footnote{The space ${\mathcal L}$ of discrete geodesic lines $l:\Z\to\Ga$ is equipped 
with the topology of pointwise convergence.
It is a locally compact Hausdorff space 
on which $\Ga$ acts properly discontinuously and cocompactly.}
sending 
$l\mapsto (\beta(l(-\infty)),\beta(l(+\infty)),l(0)x)$
is continuous.
Composing both, we see that the map 
$$ l\mapsto d\bigl(l(0)x,P(\beta(l(-\infty)),\beta(l(+\infty)))\bigr)$$
is continuous.
Since it is also $\Ga$-periodic,
the cocompactness of the action $\Ga\acts{\mathcal L}$ implies
that it is bounded,
whence the assertion. 
\qed

\medskip
From now on, we {\em assume} that the subgroup $\Ga<G$ is, in addition to being $\taumod$-boundary embedded,
also {\em $\taumod$-regular}.
This assumption will enable us to further restrict the orbit geometry
and will lead to information on the relation between the boundary embedding and the limit set.

We now analyze the position of the images of 
{\em rays} in $\Ga$ along the parallel sets.
Let $r:\N_0\to\Ga$ be a discrete geodesic ray
with ideal endpoint $\zeta:=r(+\infty)\in\geo\Ga$. 
There is a dichotomy for the position of the orbit path $rx:\N_0\to X$
relative to the Weyl cone $V(r(0)x,\st(\beta(\zeta)))$ with tip at its initial point,
namely the path must either drift away from the cone or dive deep into it:

\begin{lem}[Rays dive into Weyl cones or drift away]
\label{lem:raydivdrif}
There exist constants $\rho'=\rho'(\Ga,x)>0$ 
and for all $R>0$ numbers $n_0=n_0(\Ga,x,R)\in\N$ 
such that
the following holds:

For all $n\in\N$ with $n\geq n_0$,
the point $r(n)x$ either has 

(i) distance $\geq R$ from the Weyl cone $V(r(0)x,\st(\beta(\zeta)))$, 
or has

(ii) distance $\leq\rho'$ from this Weyl cone and distance $\geq R$ from its boundary.
\end{lem}
\proof
In a word hyperbolic group, discrete geodesic rays are contained in uniformly bounded neighborhoods of discrete geodesic lines.
Thus, $r$ is contained in a tubular neighborhood with uniform radius $c(\Ga)$
of a line $l:\Z\to\Ga$ asymptotic to $\zeta=r(+\infty)$ and some $\hat\zeta\in\geo\Ga-\{\zeta\}$.

It follows from the previous lemma that the path $rx$ is contained in a tubular neighborhood 
of the parallel set $P=P(\beta(\hat\zeta),\beta(\zeta))$
with uniform radius $\rho''(\Ga,x)$.
Let $x_0\in P$ be a point with $d(x_0,r(0)x)\leq\rho''$.
The Weyl cone 
$V(r(0)x,\st(\beta(\zeta)))$ 
is then $\rho''$-Hausdorff close to the asymptotic Weyl cone 
$V(x_0,\st(\beta(\zeta)))\subset P$.

Now we use that 
the interior of the Weyl cone $V(x_0,\st(\beta(\zeta)))$ is {\em open} in the parallel set $P$
and the boundary $\D V(x_0,\st(\beta(\zeta)))$ of the cone {\em disconnects} the parallel set, see Lemma~\ref{lem:open-cone}.
The $\taumod$-regularity of $\Ga$ implies 
(along with the triangle inequality for $\De$-lengths)
that the path $rx$ drifts away from $\D V(x_0,\st(\beta(\zeta)))$
at a {\em uniform} rate,
$$ d\bigl(r(n)x,\D V(x_0,\st(\beta(\zeta)))\bigr) \geq\phi(n) $$
with a function $\phi(n)\to+\infty$ as $n\to+\infty$
independent of the ray $r$.
The assertion follows.
\qed

\medskip
For all rays in $\Ga$ the same of the two alternatives must occur:
\begin{lem}[Dichotomy]
\label{lem:dichot}
For all discrete geodesic rays $r:\N_0\to\Ga$, 
either 

(i) 
$rx$ drifts away from the Weyl cone 
$V(r(0)x,\st(\beta(\zeta)))$, $\zeta=r(+\infty)$, at a uniform rate,
$$ d\bigl(r(n)x,V(r(0)x,\st(\beta(\zeta)))\bigr)\to+\infty$$
uniformly as $n\to+\infty$, or

(ii) $rx$ is contained in the tubular $\rho'(\Ga,x)$-neighborhood of the cone 
$V(r(0)x,\st(\beta(\zeta)))$ 
and drifts away from its boundary at a uniform rate,
$$ d\bigl(r(n)x,\D V(r(0)x,\st(\beta(\zeta)))\bigr)\to+\infty$$
uniformly as $n\to+\infty$.
\end{lem}
\proof
We give two arguments.
The first one is restricted to the nonelementary case:
As a consequence of the previous lemma,
for every ray $r$ one of the alternatives (i) and (ii) occurs with growth rates independent of the ray.
Which alternative occurs, depends only on the asymptote class 
$\zeta=r(+\infty)$ 
of the ray,
and depends on it continuously,
i.e.\ the subsets of endpoints for either alternative are open in $\geo\Ga$.
Since they are also $\Ga$-invariant,
if $|\geo\Ga|\geq3$, 
the minimality of the action $\Ga\acts\geo\Ga$ implies that one of the subsets must be empty.

The second argument works in the general case: 
Again we use that it depends only on the asymptote class of the ray, which alternative occurs.
We show that the same alternative occurs for any two distinct asymptote classes $\zeta,\hat\zeta\in\geo\Ga$.
After replacing a ray $r$ asymptotic to $\zeta$ with a subray,
we may assume that we are in the situation of the proof of the previous lemma
(whose notation we adopt),
i.e.\ that $r$ lies in a uniform tubular neighborhood of a 
line $l:\Z\to\Ga$ asymptotic to $\hat\zeta$ and $\zeta$.
Moreover, we assume that alternative (ii) holds for $\zeta$
and claim that it holds for $\hat\zeta$, as well.

To see this, fix $R>>\rho',\rho''$ and $n>>n_0$. 
Let $x_n\in P=P(\beta(\hat\zeta),\beta(\zeta))$ be a point with $d(x_n,r(n)x)\leq\rho''$.
Since (ii) holds for $r$,
the point $x_n$ must lie deep inside the cone $V(x_0,\st(\beta(\zeta)))\subset P$.
This is equivalent to $x_0$ lying deep inside the cone $V(x_n,\st(\beta(\hat\zeta)))\subset P$ 
opening towards the opposite direction.
This however implies that $r(0)x$ is uniformly close (with distance $\leq 2\rho''<<R$)
to the cone $V(r(n)x,\st(\beta(\hat\zeta)))$.
Thus alternative (ii) holds for the subray $l|_{(-\infty,n]\cap\Z}$ of $l$,
and hence also for its ideal endpoint $\hat\zeta$.
\qed

\medskip
On the other hand,
in the nonelementary case,
the ray images always drift away (at non-uniform rates) from ``opposite'' Weyl cones:

\begin{lem}[Drifting away from opposite cones]
\label{lem:drftawopp}
Suppose that $|\geo\Ga|\geq3$.
Then for every discrete geodesic ray $r:\N_0\to\Ga$ and ideal point $\hat\zeta\in\geo\Ga-\{\zeta\}$,
$\zeta=r(+\infty)$,
it holds that 
$$ d\bigl(r(n)x,V(r(0)x,\st(\beta(\hat\zeta)))\bigr)\to+\infty$$
as $n\to+\infty$.
\end{lem}
\proof
The ray $r$ is contained in a (non-uniform) tubular neighborhood of a line $l:\Z\to\Ga$ asymptotic to $\hat\zeta$ and $\zeta$.
The line image $lx$, and therefore also the ray image $rx$ is contained in a tubular neighborhood of the parallel set 
$P=P(\beta(\hat\zeta),\beta(\zeta))$. 

It follows 
that the accumulation set $\acct(r)\subset\Flagt$ of $r$ 
(with respect to $\taumod$-flag convergence, compare section~\ref{sec:conv})
consists of simplices contained in $\geo P$:
Indeed,
the nearest point projections $x_n\in P$ of $r(n)x$ lie in euclidean Weyl chambers $V(x_0,\si_n)\subset P$.
Therefore, 
in view of Lemma~\ref{lem:regseqindepch},
$\acct(r)$
equals the accumulation set of the sequence $(\tau_n)$ in $\Flagt$
consisting of the type $\taumod$ faces $\tau_n\subseteq\si_n\subset\geo P$.

Now we use nonelementarity and vary the ideal point opposite to $\zeta$.
Since $|\geo\Ga|\geq3$,
there exists a third ideal point $\hat\zeta'\in\geo\Ga-\{\zeta,\hat\zeta\}$.
It determines another parallel set $P'=P(\beta(\hat\zeta'),\beta(\zeta))$,
and the simplices in $\acct(r)$ must also be contained in $\geo P'$.
In view of $\beta(\hat\zeta)\not\subset\geo P'$,
it follows that $\beta(\hat\zeta)\not\in\acct(r)$.

Since $rx$ is contained in a tubular neighborhood of $P$,
we also again have the dichotomy, analogous to the previous lemma,
that $rx$ either drifts away from the Weyl cone $V(r(0)x,\st(\beta(\hat\zeta)))$ at a uniform rate,
as claimed,
or stays in a tubular neighborhood of it and drifts away only from its boundary. 
However, in the latter case, we would have (conical) flag convergence $r(n)\to\beta(\hat\zeta)$ as $n\to+\infty$,
equivalently, $\acct(r)=\{\beta(\hat\zeta)\}$, a contradiction. 
\qed

\medskip
If $\Ga$ is virtually cyclic, i.e.\ if $|\geo\Ga|=2$, 
there is a trivial way of modifying the boundary embedding.
Namely, 
then the action $\Ga\acts\geo\Ga$ commutes with the transposition $t:\geo\Ga\to\geo\Ga$
exchanging the points,
and therefore $-\beta:=\beta\circ t$ is a boundary embedding as well.
Therefore the previous lemma may fail. 
However, if it fails for $\beta$, then it holds for $-\beta$,
because case (ii) of the dichotomy in Lemma~\ref{lem:dichot} arises.

\medskip
From the above observations on the orbit geometry 
we will now deduce information about the limit set and its position relative to the 
image of the boundary embedding.

Let 
\begin{equation}
\label{eq:extinfin} 
\bar o_x=o_x\sqcup\beta : \ol\Ga=\Ga\sqcup\geo\Ga \lra X\sqcup\Flagt
\end{equation}
denote the extension 
of the orbit map $o_x:\Ga\to\Ga x\subset X$ to the Gromov compactification $\ol\Ga$ of $\Ga$ by 
$\ol o_x|_{\geo\Ga}:=\beta$.
We say that the extension $\ol o_x$ is {\em continuous at infinity} 
if for all sequences $\ga_n\to\infty$ in $\Ga$ 
we have flag convergence $\ga_n\to\beta(\zeta)$ 
whenever $\ga_n\to\zeta\in\geo\Ga$ in $\ol\Ga$.

We obtain the following dichotomy
corresponding to the one in Lemma~\ref{lem:dichot}:
\begin{thm}[Boundary embedding and limit set]
\label{thm:bemblim}
Let $\Ga<G$ be a $\taumod$-regular $\taumod$-boundary embedded subgroup.
Then for every boundary embedding $\beta$ 
either 

(i) $\beta(\geo\Ga)\cap\LatGa=\emptyset$, 
and no simplex in $\beta(\geo\Ga)$ is opposite to a simplex in $\LatGa$,\footnote{Note that 
in view of  the antipodality of $\beta$
the second part of (i) implies the first part.}
or

(ii) $\beta(\geo\Ga)=\LatGa$.
Moreover, 
the extension $\bar o_x$ is continuous at infinity,
after replacing $\beta$ with $-\beta$ in the case $|\geo\Ga|=2$, if necessary.
\end{thm}
\proof
Assume first that case (ii) of Lemma~\ref{lem:dichot} occurs.
Consider a sequence $\ga_n\to\infty$ in $\Ga$.
There exist rays $r_n:\N_0\to\Ga$ starting in $r_n(0)=e$ and passing at uniformly bounded distance of $\ga_n$.
We denote their ideal endpoints by $\zeta_n:=r_n(+\infty)$.
Then the orbit points $\ga_nx$ lie in uniform tubular neighborhoods of the Weyl cones 
$V(x,\st(\beta(\zeta_n)))$.
If $\ga_n\to\zeta\in\geo\Ga$ in $\ol\Ga$,
equivalently, 
$\zeta_n\to\zeta$ in $\geo\Ga$, 
then $\beta(\zeta_n)\to\beta(\zeta)$ in $\Flagt$,
and it follows $\taumod$-flag convergence $\ga_nx\to\beta(\zeta)$.
This shows that $\bar o_x$ is continuous at infinity and $\beta(\geo\Ga)\subseteq\LatGa$.
To see the opposite inclusion,
suppose that $\ga_nx\to\la\in\LatGa$.
After extraction, we get convergence $\ga_n\to\zeta\in\geo\Ga$ and conclude from the above that $\la=\beta(\zeta)$.
Thus also $\LatGa\subseteq\beta(\geo\Ga)$,
and conclusion (ii) of the theorem is satisfied.

If $|\geo\Ga|=2$ and case (ii) of Lemma~\ref{lem:dichot} occurs for $-\beta$,
we reach the same conclusion after replacing $\beta$ with $-\beta$.

Assume now that we are in case (i) of Lemma~\ref{lem:dichot}.
After replacing $\beta$ with $-\beta$ in the case $|\geo\Ga|=2$, if necessary,
we may also assume that the conclusion of Lemma~\ref{lem:drftawopp} holds.
As before, we consider a sequence $\ga_n\to\infty$ in $\Ga$ and rays $r_n$.
Suppose that $\ga_n\to\zeta\in\geo\Ga$ and let $\hat\zeta\in\geo\Ga-\{\zeta\}$ be arbitrary.
Since $\zeta_n\to\zeta$, 
there exist for all large $n$ lines $l_n:\Z\to\Ga$ with ideal endpoints 
$l_n(-\infty)=\hat\zeta$ and $l_n(+\infty)=\zeta_n$.
The lines $l_n$ pass at uniformly bounded distance from $e$ and $\ga_n$, 
and they contain the rays $r_n$ in uniform tubular neighborhoods.
(For the rest of this argument, uniformity will mean that bounds are independent of $n$.)

By Lemma~\ref{lem:lnclpar}, 
the ray images $r_nx$ lie in uniform tubular neighborhoods of the parallel sets
$P_n=P(\beta(\hat\zeta),\beta(\zeta_n))$
and drift away from both Weyl cones 
$V(x,\st(\beta(\hat\zeta)))$ and $V(x,\st(\beta(\zeta_n)))$.
The drift is uniform in the latter case by Lemma~\ref{lem:dichot}(i),
and also in the former case since $r_n(0)x=x$ and $d(x,P_n)$ is bounded. 

The uniformity implies 
that the orbit points $\ga_nx$ lie in uniform tubular neighborhoods of Weyl cones 
$V(x,\st(\tau_n))$ for simplices $\tau_n\in\Flagt$ 
with $\tau_n\subset\geo P_n$
but $\tau_n\neq\beta(\hat\zeta),\beta(\zeta_n)$.
(Indeed, as in the proof of the previous lemma,
$\ga_nx$ is uniformly close to a euclidean Weyl chamber $V(x,\si_n)$
with visual boundary chamber $\si_n\subset\geo P_n$ but $\si_n\not\subset\st(\beta(\hat\zeta))\cup\st(\beta(\zeta_n))$,
and we let $\tau_n\subseteq\si_n$ be the type $\taumod$ face.)
In particular, $\tau_n$ is not opposite to both $\beta(\hat\zeta)$ and $\beta(\zeta_n)$.
The accumulation set of the sequence $(\tau_n)$ in $\Flagt$, 
which coincides with the $\taumod$-flag accumulation set of the sequence $(\ga_n)$,
therefore consists of simplices 
which are not opposite to both $\beta(\hat\zeta)$ and $\beta(\zeta)$,
because 
oppositeness 
is an open property.
Letting $\hat\zeta$ run through $\geo\Ga-\{\zeta\}$, it follows that these simplices are not opposite 
to {\em any} simplex in $\beta(\geo\Ga)$.

Every limit simplex in $\LatGa$ arises as the $\taumod$-flag limit of a sequence $(\ga_n)$ which converges at infinity in $\ol\Ga$.
We obtain that no simplex in $\LatGa$ is opposite to a simplex in $\beta(\geo\Ga)$.
In particular, $\LatGa\cap\beta(\geo\Ga)=\emptyset$.
Thus, conclusion (i) of the theorem holds.
\qed

\medskip
Consequently,
as soon as a boundary embedding hits the limit set, it identifies it with the Gromov boundary of the subgroup 
and moreover continuously extends the orbit maps: 

\begin{cor}
\label{cor:bembaemb}
Let $\Ga<G$ be a $\taumod$-regular $\taumod$-boundary embedded subgroup
with boundary embedding $\beta$.
If $\beta(\geo\Ga)\cap\LatGa\neq\emptyset$,
then $\beta(\geo\Ga)=\LatGa$.
Moreover, the extension $\bar o_x$ is continuous at infinity,
after replacing $\beta$ with $-\beta$ in the case $|\geo\Ga|=2$, if necessary.
\end{cor}

Otherwise, 
if the boundary embedding avoids the limit set, 
the image of the boundary embedding and the limit set must have special position:
\begin{lem}
\label{lem:subvr}
In case (i) of Theorem~\ref{thm:bemblim},
both $\beta(\geo\Ga)$ and $\LatGa$ are not Zariski dense in $\Flagt$.
In particular, $\Ga$ is not Zariski dense in $G$.
\end{lem}
\proof
Since no simplex in $\beta(\geo\Ga)$ is opposite to a simplex in $\LatGa$,
it follows that $\beta(\geo\Ga)$ is disjoint from the union of open Schubert strata $C(\la)$
over all limit simplices $\la\in\LatGa$.
In other words,
$\beta(\geo\Ga)$ is contained 
in the intersection of the proper subvarieties $\D C(\la)=\Flagt-C(\la)$.
Similarly, 
$\LatGa$ lies in the intersection of the $\D C(\tau)$ over all simplices $\tau\in\beta(\geo\Ga)$.
In particular, both are $\Ga$-invariant proper subvarieties,
which forces $\Ga$ to be non-Zariski dense. 
\qed

\medskip
Therefore, the first alternative in the theorem cannot occur
in the Zariski dense case,
compare \cite[Thm.\ 1.5]{GW}:
\begin{cor}
\label{cor:zdnsasemb}
Let $\Ga<G$ be a Zariski dense $\taumod$-regular $\taumod$-boundary embedded subgroup.
Then it admits a unique boundary embedding $\beta$,
and $\beta(\geo\Ga)=\LatGa$.
\end{cor}
\proof
By the lemma,
for any boundary embedding $\beta$,
only case (ii) in the theorem can occur.
It follows that $\beta(\geo\Ga)=\LatGa$.
Moreover, $\beta$ is uniquely determined
because,
due to the density of attractive fixed points of infinite order elements,
there are no $\Ga$-equivariant self homeomorphisms of $\geo\Ga$ besides the identity.
(Note that $|\geo\Ga|\geq3$ by Zariski density.)
\qed

\medskip
It is worth noting that in the case $\taumod=\simod$ the boundary embedding can always be modified 
so that it maps onto the limit set:

\begin{thm}
\label{thm:bemmdsi}
Let $\Ga<G$ be a $\simod$-regular $\simod$-boundary embedded subgroup.
Then 
there exists a boundary embedding $\beta$ with $\beta(\geo\Ga)=\LasGa$.
\end{thm}
\proof
In the case $\taumod=\simod$,
the parallel sets considered above are maximal flats 
and the Weyl cones are euclidean Weyl chambers.
What makes it possible to push the argument further,
is the fact that the walls in a maximal flat through a fixed point disconnect the flat 
into euclidean Weyl chambers.
Therefore, the above discussion now yields more precise information 
about the position of the paths $rx$:

Since the $rx$ are uniformly close to maximal flats
(provided by a boundary embedding $\beta'$ for $\Ga$, cf.\ Lemma~\ref{lem:lnclpar}),
$\simod$-regularity forces them to dive into (uniform tubular neighborhoods of) Weyl chambers inside these flats.
It follows that the paths $rx$ are contained 
in uniform tubular neighborhoods of euclidean Weyl chambers with tips at the initial points $r(0)x$.
Again by regularity,
the asymptote class of the Weyl chamber depends only on the asymptote class of the ray $r$.
We therefore obtain a new boundary map $\beta:\geo\Ga\to\Flags$ 
such that $rx$ is contained in the tubular $\rho'(\Ga,x)$-neighborhood of the euclidean Weyl chamber $V(r(0)x,\beta(\zeta))$
for $\zeta=r(+\infty)$. 
Clearly, $\beta(\geo\Ga)\subseteq\LasGa$ 
and $\beta$ is $\Ga$-equivariant.
An argument as in the last part of the proof of Lemma~\ref{lem:dichot} shows that $\beta$ is antipodal.

To verify that $\beta$ is continuous,
suppose that $\zeta_n\to\zeta$ in $\geo\Ga$ and $\beta(\zeta_n)\to\si$ in $\Flags$.
We must show that $\si=\beta(\zeta)$. 
Let $r_n,r:\N_0\to\Ga$ be rays starting in $e$ and asymptotic to $\zeta_n,\zeta$.
We note that for any sequence $m_n\to+\infty$ in $\N_0$, 
we have $\simod$-flag convergence $r_n(m_n)\to\si$,
because $r_n(m_n)x$ lies in a uniform tubular neighborhood of $V(x,\st(\beta(\zeta_n)))$.
On the other hand, 
if $m_n$ grows sufficiently slowly,
then the sequence $(r_n(m_n))$ in $\Ga$ is contained in a tubular neighborhood of $r$,
and hence $r_n(m_n)\to\beta(\zeta)$.
This shows that $\si=\beta(\zeta)$, as desired.

Thus, $\beta$ is a boundary embedding.
Since also $\beta(\geo\Ga)\subseteq\LasGa$,
we conclude using Theorem~\ref{thm:bemblim} that $\beta(\geo\Ga)=\LasGa$.
\qed

\subsection{Asymptotic embeddings and coarse extrinsic geometry}
\label{sec:asyemb} 

The discussion in the previous section,
notably part (ii) of the conclusion of Theorem~\ref{thm:bemblim},
motivates the following strengthening 
of the notion of boundary embeddedness:

\begin{dfn}[Asymptotically embedded]
\label{def:asyemb}
A subgroup $\Ga<G$ is {\em $\taumod$-asymptotically embedded} 
if it is $\taumod$-regular, $\taumod$-antipodal,
intrinsically word hyperbolic 
and there is a $\Ga$-equivariant homeomorphism 
\begin{equation*}
\label{eq:mapalphatau}
\alpha: \geo \Ga \stackrel{\cong}{\lra}
\LatGa\subset \Flagt
\end{equation*}
from its Gromov boundary onto its $\tau_{mod}$-limit set. 
\end{dfn}

The definition can also be 
phrased 
purely {\em dynamically} in terms of the 
$\Ga$-action on $\Flagt$,
by replacing $\taumod$-regularity with the $\taumod$-convergence condition. 

Note that $\taumod$-asymptotically embedded subgroups are necessarily {\em discrete} by $\taumod$-regularity.
We also keep assuming that $\taumod$ is {\em $\iota$-invariant};
this is implicit in $\taumod$-antipodality.

We observe that the boundary map $\al$ is antipodal, 
because it is injective with antipodal image.
It is therefore a boundary embedding for $\Ga$,
i.e.\ $\taumod$-asymptotically embedded implies $\taumod$-boundary embedded.
According to Corollary~\ref{cor:bembaemb}, 
the extension
\begin{equation}
\label{eq:extinfinal}
\bar o_x=o_x\sqcup\al : \ol\Ga=\Ga\sqcup\geo\Ga \lra X\sqcup\Flagt
\end{equation}
cf.\ (\ref{eq:extinfin}),
is continuous, 
after replacing $\al$ with $-\al$ in the case $|\geo\Ga|=2$, if necessary.
We will refer to $\al$ then as the {\em asymptotic embedding} for $\Ga$.

\medskip
We rephrase the {\em criteria} for asymptotic embeddedness
obtained in the previous section 
(cf.\ Corollaries~\ref{cor:bembaemb}, \ref{cor:zdnsasemb} and Theorem~\ref{thm:bemmdsi}):
\begin{thm}
\label{thm:bembaembref}
Let $\Ga<G$ be a $\taumod$-regular $\taumod$-boundary embedded subgroup 
with boundary embedding $\beta$.
If $\beta(\geo\Ga)\cap\LatGa\neq\emptyset$,
then $\Ga$ is $\taumod$-asymptotically embedded,
and $\beta$ is the asymptotic embedding, 
after replacing it with $-\beta$ in the case $|\geo\Ga|=2$, if necessary.
\end{thm}

\begin{thm}
\label{thm:zdnsasembref}
Zariski dense $\taumod$-regular $\taumod$-boundary embedded subgroups 
are $\taumod$-asymp\-to\-ti\-cally embedded 
and admit no other boundary embedding besides their asymptotic embedding.
\end{thm}

\begin{thm}
\label{thm:bemmdsiref}
$\simod$-Regular $\simod$-boundary embedded subgroups
are $\simod$-asymptotically embedded.
(But they may admit boundary embeddings different from the asymptotic embedding.)
\end{thm}

We also summarize 
what the discussion in the previous section yields for the {\em orbit geometry} 
of asymptotically embedded subgroups.
In addition to the continuity at infinity (\ref{eq:extinfinal}) of the orbit maps $o_x$, $x\in X$,
we obtained (cf.\ Lemmas~\ref{lem:lnclpar} and~\ref{lem:dichot}):

\begin{prop}[Orbit geometry of asymptotically embedded subgroups]
\label{prop:asyembprp}
Let $\Ga<G$ be a $\taumod$-asymptotically embedded subgroup
with asymptotic embedding $\al$.
Then: 

(i) 
For every discrete geodesic line $l:\Z\to\Ga$,
the path $lx$ is contained in a tubular neighborhood 
of uniform radius $\rho(\Ga,x)$
of the parallel set $P(\al(\zeta_-),\al(\zeta_+))$, 
where  $\zeta_{\pm}:=l(\pm\infty)\in\geo\Ga$.

(ii)
For every discrete geodesic ray $r:\N_0\to\Ga$,
the path $rx$ is contained in a tubular neighborhood 
of uniform radius $\rho'(\Ga,x)$ 
of the Weyl cone $V(r(0)x,\st(\al(\zeta)))$,
where  $\zeta:=r(+\infty)\in\geo\Ga$,
and drifts away from its boundary at a uniform rate,
\begin{equation}
\label{eq:awbdcn}
d\bigl(r(n)x,\D V(r(0)x,\st(\al(\zeta)))\bigr)\to+\infty
\end{equation}
uniformly as $n\to+\infty$.
\end{prop}

These properties motivate the {\em Morse property} 
to be introduced and discussed below.
Let us first draw some further immediate consequences for the coarse extrinsic geometry of subgroups
and see how property (ii) leads to {\em undistortion} and {\em uniform regularity}. 

We consider 
the orbit path $rx$ for a discrete ray $r$.
According to property (ii),
the path $rx$ must stay uniformly close to the Weyl cone $V(r(0)x,\st(\al(\zeta)))$ 
predicted by the boundary map and
drift away from the boundary of the cone at a uniform rate. 
Since the same applies to all subrays of $r$,
it follows that the cones $V(r(n)x,\st(\al(\zeta)))$ must, up to bounded perturbation,
be {\em uniformly nested}.
This forces the orbit path $rx$ to have a {\em linear} drift 
away from the boundary of the Weyl cone and in particular towards infinity,
i.e.\ $rx$ is uniformly $\taumod$-regular and undistorted.

We combine these properties in the following notion:

\begin{dfn}[URU]
\label{dfn:uru}
A finitely generated subgroup $\Ga<G$ is {\em $\taumod$-URU},
if it is 

(i) uniformly $\taumod$-regular, and 

(ii) undistorted,
i.e.\ the inclusion $\Ga\subset G$, equivalently,
the orbit maps $\Ga\to\Ga x\subset X$, are {\em quasiisometric embeddings} 
with respect to a word metric on $\Ga$.
\end{dfn}

Note that URU subgroups cannot contain parabolic elements.

The above discussion before the definition thus leads to:
\begin{thm}
\label{thm:asymbur}
$\taumod$-Asymptotically embedded subgroups $\Ga<G$ are $\taumod$-URU.
\end{thm}
\proof
We add some details to the discussion above:

Let $x_n\in V(r(0)x,\st(\al(\zeta)))$ be the nearest point projections of the points $r(n)x$, $n\in\N_0$.
Then $d(r(n)x,x_n)\leq\rho'=\rho'(\Ga,x)$ by part (ii) of the proposition.
We consider the sequence of Weyl cones $V(x_n,\st(\al(\zeta)))\subset V(r(0)x,\st(\al(\zeta)))$. 
Note that 
the cones $V(r(n)x,\st(\al(\zeta)))$ and $V(x_n,\st(\al(\zeta)))$ are asymptotic to each other 
and have Hausdorff distance $\leq d(r(n)x,x_n)\leq\rho'$,
as do their boundaries.
Applying (ii) to the subrays of $r$,
it follows that the paths $m\mapsto r(n+m)x$ are contained in uniform neighborhoods 
of the cones $V(x_n,\st(\al(\zeta)))$ and drift away from their boundaries at uniform rates.
Thus, 
for every $d_0>0$ there exists a number $m_0=m_0(\Ga,x,d_0)\in\N$
such that 
$$x_{n+m}\in V(x_n,\st(\al(\zeta)))$$
and 
$$d\bigl(x_{n+m},\D V(x_n,\st(\al(\zeta)))\bigr)\geq d_0$$
for all $n\geq0$ and $m\geq m_0$.
The latter inequality implies that the boundaries of the Weyl cones $V(x_n,\st(\al(\zeta)))$ and $V(x_{n+m},\st(\al(\zeta)))$
have (nearest point) distance $\geq d_0$,
cf.\ Proposition~\ref{prop:hdstwcones}(ii).
From the uniform nestedness 
of the cones $V(x_{km_0},\st(\al(\zeta)))$ for $k\in\N_0$,
it follows that the drift (\ref{eq:awbdcn}) away from the boundary of the Weyl cone is {\em uniformly linear}.
Consequently, 
the ray images $rx$ are uniformly undistorted and uniformly $\taumod$-regular.
Since any pair of elements in $\Ga$ lies in a uniform tubular neighborhood of some discrete geodesic ray,
our assertion follows.
\qed

\begin{rem} 
(i) That, conversely, URU implies asymptotic embeddedness is proven in \cite{mlem}.
In particular, URU subgroups are necessarily word hyperbolic. 

(ii) In \cite{bordif} we prove that URU subgroups $\Ga<G$ satisfy the even stronger coarse geometric property
of being {\em coarse Lipschitz retracts} of $G$.
\end{rem}

Similarly,
we also derive a version of Proposition~\ref{prop:asyembprp}
for discrete geodesic {\em segments} in $\Ga$: 

Consider a line $l:\Z\to\Ga$ 
and denote $\zeta_{\pm}=l(\pm\infty)$. 
Let $x_n\in P(\al(\zeta_-),\al(\zeta_+))$ be 
the nearest point projections of the points $l(n)x$, $n\in\Z$. 
As in the proof of the previous theorem,
we see using Proposition~\ref{prop:asyembprp}(i+ii),
that for any $d_0>0$ there exists $m'_0=m'_0(\Ga,x,d_0)\in\N$
such that 
$$x_{n\pm m}\in V(x_n,\st(\al(\zeta_{\pm})))$$
and 
$$d\bigl(x_{n\pm m},\D V(x_n,\st(\al(\zeta_{\pm})))\bigr)\geq d_0$$
for all $n$ and $m\geq m'_0$.
It follows that, 
for $n_{\pm}\in\Z$ with $n_+-n_-\geq m'_0$,
the diamond 
$$\diamot(x_{n_-},x_{n_+})=
V(x_{n_-},\st(\al(\zeta_+)))\cap V(x_{n_+},\st(\al(\zeta_-)))\subset P(\al(\zeta_-),\al(\zeta_+)) $$
is defined and,
using Proposition~\ref{prop:asyembprp}(ii) again,
contains the finite subpath
$l|_{[n_-,n_+]\cap\Z}x$
in a uniform tubular neighborhood.

Our discussion yields the following complement to, respectively, strengthening of Proposition~\ref{prop:asyembprp},
saying that the images of discrete geodesic segments in $\Ga$ are contained in uniform neighborhoods 
of diamonds with tips at uniform distance from the endpoints:
\begin{prop}[Segments go close to diamonds]
\label{prop:cldm}
Let $\Ga<G$ be a $\taumod$-asymptotically embedded subgroup.
Then 
for every discrete geodesic segment $s:[n_-,n_+]\cap\Z\to\Ga$,
the path $sx$ is contained in a tubular neighborhood of uniform radius $\rho''=\rho''(\Ga,x)$
of a diamond $\diamot(x_-,x_+)$ with 
$d(x_{\pm},s(n_{\pm})x)\leq\rho''$.
\end{prop}
\proof
This is a consequence of the above discussion, because every discrete geodesic segment in $\Ga$ 
lies in a uniform neighborhood of a discrete geodesic line.
\qed

\subsection{Morse property}
\label{sec:morse}

The {\em Morse Lemma} for Gromov hyperbolic spaces
asserts that quasigeodesic segments are uniformly close to geodesic segments with the same endpoints. 
Proposition~\ref{prop:cldm} along with Proposition~\ref{prop:asyembprp}
in the previous section 
can be interpreted 
as saying that,
for asymptotically embedded subgroups $\Ga<G$, 
the images of discrete geodesic segments, rays and lines in $\Ga$ 
under the orbit maps into $X$
satisfy a higher rank version of the Morse Lemma,
with geodesic segments 
replaced by diamonds.

This motivates the following notion
(we keep assuming that $\taumod$ is $\iota$-invariant):

\begin{dfn}[Morse]
\label{def:mrs}
A subgroup $\Ga<G$ is {\em $\taumod$-Morse}
if it is $\taumod$-regular, intrinsically word hyperbolic 
and satisfies the following property: 

For every discrete geodesic segment $s:[n_-,n_+]\cap\Z\to\Ga$,
the path $sx$ is contained in a tubular neighborhood of uniform radius $\rho''=\rho''(\Ga,x)$
of a diamond $\diamot(x_-,x_+)$ with 
tips at distance $d(x_{\pm},s(n_{\pm})x)\leq\rho''$ from the endpoints.
\end{dfn}

Note 
that the definition does not a priori assume the existence of a boundary map,
neither does it assume undistortion.
These will be consequences.

As we saw, 
asymptotically embedded subgroups are Morse.
We will now show that, conversely,
asymptotic embeddedness follows from the Morse property,
in fact from an a priori weaker version of it for rays in $\Ga$ 
(instead of segments):

\begin{thm}
\label{thm:mrsasembeq}
For a 
subgroup $\Ga<G$ the following properties are equivalent:

(i) $\Ga$ is $\taumod$-asymptotically embedded.

(ii) $\Ga$ is $\taumod$-Morse.

(iii) $\Ga$ is $\taumod$-regular, intrinsically word hyperbolic 
and satisfies the following property: 
For every discrete geodesic ray $r:\N_0\to\Ga$,
the path $rx$
is contained in a tubular neighborhood of uniform radius $\rho'''=\rho'''(\Ga,x)$
of a $\taumod$-Weyl cone with tip at the initial point $r(0)x$.

The $\taumod$-Weyl cone in (iii) is then the cone $V(r(0)x,\al(r(+\infty)))$
where $\al$ is the asymptotic embedding for $\Ga$.
\end{thm}
\proof
The implication (i)$\Ra$(ii) is Proposition~\ref{prop:cldm}.
The implication (ii)$\Ra$(iii) is immediate by a limiting argument.
It remains to show that (iii)$\Ra$(i).

We first observe that 
the $\taumod$-Weyl cone $V(r(0)x,\st(\tau))$
containing the path $rx$ in a tubular neighborhood
is uniquely determined.
This follows from the $\taumod$-flag convergence $r(n)\to\tau$.
Moreover, $\tau$ depends only on the asymptote class $r(+\infty)$ of the ray $r$.
Hence there is a well-defined map at infinity
$$ \check\al:\geo\Ga \to \Flagt $$
such that for every ray $r$ the path $rx$ is contained in a uniform tubular neighborhood of 
the Weyl cone 
$V(r(0)x,\st(\check\al(r(+\infty))))$.
Our goal is to show that $\check\al$ is an asymptotic embedding.

\begin{lem}
$\check\al$ is continuous and 
continuously extends the orbit maps $o_x$ at infinity.
\end{lem}
\proof
We proceed as in the proof of Theorem~\ref{thm:bemmdsi} (continuity of $\beta$).
Consider a converging sequence $\zeta_n\to\zeta$ in $\geo\Ga$.
Let $r_n,r:\N_0\to\Ga$ be rays starting in $e$ and asymptotic to $\zeta_n,\zeta$.
We note that for any sequence $m_n\to+\infty$ in $\N_0$, 
the flag accumulation set of the sequence $(r_n(m_n))$ in $\Flagt$ 
equals 
the accumulation set of the sequence $(\check\al(\zeta_n))$ in $\Flagt$,
and in particular does not depend on the sequence $(m_n)$.
On the other hand, 
if $(m_n)$ grows sufficiently slowly,
then the sequence $(r_n(m_n))$ in $\Ga$ is contained in a tubular neighborhood of $r$,
and hence flag converges to $\check\al(\zeta)$.
It follows that $\check\al(\zeta_n)\to\check\al(\zeta)$.
This shows that $\check\al$ is continuous.

Proceeding as in the first part of the proof of Theorem~\ref{thm:bemblim},
we then see that, for a sequence $\ga_n\to\infty$ in $\Ga$, 
convergence $\ga_n\to\zeta\in\geo\Ga$ in $\ol\Ga$ implies 
flag convergence $\ga_n\to\check\al(\zeta)$,
i.e.\ $\check\al$ continuously extends $o_x$ at infinity.
\qed

\medskip

The continuous extension part of the lemma implies 

\begin{cor}\label{cor:Morse-limit}
 $\check\al(\geo\Ga)=\LatGa$.
\end{cor}

In order to see that $\LatGa$ is antipodal and 
$\check\al$ is an asymptotic embedding for $\Ga$, 
it remains to verify:
\begin{lem}
The map $\check\al$ is antipodal.
\end{lem}
\proof
Let $\zeta_{\pm}\in\geo\Ga$ be distinct,
and let $l:\Z\to\Ga$ be a line with $l(\pm\infty)=\zeta_{\pm}$.
Applying property (iii) to the subrays $l|_{[-n,+\infty)}$ for large $n\in\N$,
we get that the point $l(0)x$ is uniformly close to the cones $V(l(-n)x,\st(\check\al(\zeta_+)))$,
equivalently,
there exists a bounded sequence of points $y_n\in V(l(-n)x,\st(\check\al(\zeta_+)))$.
By $\taumod$-regularity, $d(y_n,\D V(l(-n)x,\st(\check\al(\zeta_+))))\to+\infty$ as $n\to+\infty$.
We denote by $\tau_n^-\in\Flagt$ the simplex $l(-n)x$-opposite to $y_n$.\footnote{
I.e.\ $l(-n)x\in V(y_n,\st(\tau_n^-))$.
Then $l(-n)x,y_n\in P(\tau_n^-,\check\al(\zeta_+))$.}
Then $l(-n)x\in V(y_n,\st(\tau_n^-))$, and hence $l(-n)x$ is uniformly close to $V(l(0)x,\st(\tau_n^-))$.
In view of the flag convergence $l(-n)\to\check\al(\zeta_-)$, it follows that $\tau_n^-\to\check\al(\zeta_-)$ in $\Flagt$.
Since the parallel sets $P(\tau_n^-,\check\al(\zeta_+))$ lie at bounded distance from $l(0)x$,
as they contain the points $y_n$,
the sequence $(\tau_n^-)$ is relatively compact in 
the open Schubert stratum
$C(\check\al(\zeta_+))$.
Hence $\check\al(\zeta_-)\in C(\check\al(\zeta_+))$,
i.e.\ $\check\al(\zeta_-)$ is opposite to $\check\al(\zeta_+)$.
\qed

\medskip
This concludes the proof of the theorem.
\qed

\medskip
Note that the theorem implies in particular that 
{\em $\taumod$-Morse subgroups are $\taumod$-URU},
because asymptotically embedded subgroups are URU 
by Theorem~\ref{thm:asymbur}.

\begin{rem}
We restricted our definition of the Morse property to {\em word hyperbolic} subgroups 
because, as shown in \cite{mlem}, URU subgroups are always word hyperbolic.
This had been unknown at the time of writing the first version of \cite{morse}.
\end{rem}

\subsection{Conicality}
\label{sec:conic}

The condition for discrete subgroups which we study
in this section 
concerns the asymptotic geometry of their orbits,
i.e.\ how they approach infinity.
To state it, 
we first need to elaborate on 
our discussion of convergence at infinity for sequences from section~\ref{sec:conv}.

For arbitrary $\taumod$,
consider a $\taumod$-flag converging sequence $(x_n)$ in $X$,
$$ x_n\to\tau\in\Flagt .$$
The following notion of 
going ``straight'' to the limit simplex
generalizes {\em conical} or {\em radial} convergence at infinity 
in rank one symmetric spaces
where one requires the sequence to stay in a tubular neighborhood of a geodesic ray.
Working with rays also in higher rank turns out to be too restrictive,\footnote{From 
our construction of Anosov Schottky subgroups, see \cite{morse}, it immediately follows 
that in higher rank they are generically not ray conical,
for instance never in the Zariski dense case.
This implies furthermore
that Zariski dense Anosov subgroups are never ray conical.}
and we replace the rays with Weyl cones,
compare \cite[Def.\ 5.2]{Albuquerque}:

\begin{dfn}[Conical convergence]
A $\taumod$-flag converging sequence $x_n\to\tau\in\Flagt$ 
converges {\em $\taumod$-conically},
$$x_n\stackrel{con}{\to}\tau ,$$
if it is contained in a tubular neighborhood of a Weyl cone $V(x,\st(\tau))$ for some point $x\in X$.

Accordingly, 
$\taumod$-flag converging sequences in $G$ are said to converge {\em $\taumod$-conically}
if their orbit sequences in $X$ do.
\end{dfn}

Note that the Weyl cones $V(x,\st(\tau))$ for different points $x\in X$ are Hausdorff close to each other,
and the conical convergence condition is therefore independent of the choice of $x$.

The next result describes a situation for sequences close to 
parallel sets
where flag convergence already implies the stronger form of conical convergence:

\begin{lem}
\label{lem:frconcnv}
Suppose that a sequence $(x_n)$ in $X$ $\taumod$-flag converges,
$x_n\to\tau\in\Flagt$.

(i) 
If $(x_n)$ is contained in a tubular neighborhood of a parallel set $P(\hat\tau,\tau)$
for some $\hat\tau\in C(\tau)$,

(ii) or if, more generally,
there exists a relatively compact sequence $(\hat\tau_n)$ in $C(\tau)$ such that 
\begin{equation*}
\sup_n d(x_n,P(\hat\tau_n,\tau))<+\infty ,
\end{equation*}
then $x_n\stackrel{con}{\to}\tau$.
\end{lem}
\proof
Suppose first that 
the stronger 
condition (i) holds and 
that $x_n\stackrel{con}{\nrightarrow}\tau$.
Let $x\in P(\hat\tau,\tau)$.
As in the proof of Lemma~\ref{lem:raydivdrif},
it follows from the openness of the cone $V(x,\st(\tau))$ in the parallel set $P(\hat\tau,\tau)$ 
that, after extraction, 
the sequence $(x_n)$ drifts away from $V(x,\st(\tau))$.
As in the proof of Theorem~\ref{thm:bemblim},
the points $x_n$ are then contained in uniform neighborhoods 
of cones $V(x,\st(\tau_n))$
with simplices $\tau_n\in\Flagt$ satisfying 
$\tau_n\subset\geo P(\hat\tau,\tau)$ but $\tau_n\neq\tau$.
Since $\tau$ is the only simplex in $C(\hat\tau)$ which lies in $P(\hat\tau,\tau)$,
see \eqref{eq:oppinsusp} and the discussion preceding Lemma \ref{lem:opst}, 
the sequence $(\tau_n)$ is contained in the closed set $\Flagt-C(\hat\tau)$,
and hence so is its accumulation set.
In particular,
$\tau$ does not belong to the accumulation set of $(\tau_n)$ in $\Flagt$.
Since the latter set equals the flag accumulation set of the sequence $(x_n)$ in $\Flagt$,
it follows in particular that 
$x_n\nrightarrow\tau$,
a contradiction. 

Suppose now that the weaker condition (ii) holds. 
Since $C(\tau)$ is a homogeneous $P_{\tau}$-space, 
there exist $\hat\tau\in C(\tau)$ and a bounded sequence $(b_n)$ in $P_{\tau}$ such that $\hat\tau_n=b_n\hat\tau$.
The sequence $(b_n^{-1}x_n)$ is then contained in a tubular neighborhood of $P(\hat\tau,\tau)$, 
i.e.\ it satisfies condition (i).
Moreover, we also have flag convergence $b_n^{-1}x_n\to\tau$.\footnote{Because
the $b_n$ are bounded and fix $\tau$ on $\Flagt$.}
Hence, by the above, 
it follows that $b_n^{-1}x_n\stackrel{con}{\to}\tau$.
By the definition of conical convergence, 
this means that the sequence $(b_n^{-1}x_n)$ lies in a tubular neighborhood of the cone $V(x,\st(\tau))$
for some point $x\in X$, equivalently,
that 
$$\sup_n d(x_n,V(b_nx,\st(\tau)))<+\infty .$$
Now the cones $V(b_nx,\st(\tau))$ 
are asymptotic to 
$V(x,\st(\tau))$ and have finite Hausdorff distance $\leq d(x,b_nx)$ from it. 
This Hausdorff distance is uniformly bounded and it also follows that the sequence 
$(x_n)$
lies in a tubular neighborhood of $V(x,\st(\tau))$,
i.e.\ $x_n\stackrel{con}{\to}\tau$.
\qed

\medskip
As we did with regularity and flag convergence, 
we will now also rephrase
conical convergence for sequences in $G$ in terms of their dynamics on flag manifolds.

For a flag convergent sequence, 
conical convergence is reflected as follows by the dynamics on the space of parallel sets,
equivalently, on 
the space of pairs of opposite simplices,
cf.\ (\ref{eq:oppfl}):

\begin{lem}
Suppose that a sequence $(g_n)$ in $G$ $\taumod$-flag converges,
$g_n\to\tau\in\Flagt$.
Then for a relatively compact sequence $(\hat\tau_n)$ in $C(\tau)$,
the following are equivalent:

(i) 
$g_n\stackrel{con}{\to}\tau$.

(ii)
The 
parallel sets $g_n^{-1}P(\hat\tau_n,\tau)$ 
all intersect a fixed bounded subset in $X$.

(ii')
The sequence of pairs 
$g_n^{-1}(\hat\tau_n,\tau)$ is relatively compact in $(\Flagit\times\Flagt)^{opp}$.

\end{lem}
\proof
We first note that conditions (ii) and (ii') are equivalent 
as a consequence of:
\begin{sublem}
A subset $A\subset(\Flagit\times\Flagt)^{opp}$ is relatively compact 
iff  the corresponding parallel sets $P(\tau_-,\tau_+)$ for $(\tau_-,\tau_+)\in A$
all intersect a fixed bounded subset of $X$, i.e.\ 
$$ \sup_{(\tau_-,\tau_+)\in A} d(x,P(\tau_-,\tau_+))<+\infty$$
for a base point $x\in X$.
\end{sublem}
\proof
The forward direction follows from the continuity of the function (\ref{eq:dstcnt}).\footnote{Since here $\taumod$
is not required to be $\iota$-invariant, we consider the function on $(\Flagit\times\Flagt)^{opp}\times X$.}

For the converse direction
we note that for a pair $(\tau_-,\tau_+)\in (\Flagit\times\Flagt)^{opp}$
the intersection of parabolic subgroups $P_{\tau_-}\cap P_{\tau_+}$
preserves the parallel set $P(\tau_-,\tau_+)$ and acts transitively on it. 
Consequently,
the set of triples $(\tau_-,\tau_+,x')\in(\Flagit\times\Flagt)^{opp}\times X$
such that $x'\in P(\tau_-,\tau_+)$ is still a homogeneous $G$-space.
Let us fix in it a reference triple $(\tau_0^-,\tau_0^+,x)$.
Then the parallel sets $P(\tau_-,\tau_+)$ 
intersecting a closed ball $\ol B(x,R)$ 
are of the form $gP(\tau_0^-,\tau_0^+)$ with $g\in G$ such that $d(x,gx)\leq R$.
It follows that the set of these pairs $(\tau_-,\tau_+)=g(\tau_0^-,\tau_0^+)$
is compact.
\qed

\medskip
Continuing with the proof of the lemma, 
let $x\in X$ be a base point. 
In view of 
$$ d(x,g_n^{-1}P(\hat\tau_n,\tau))=d(g_nx,P(\hat\tau_n,\tau)) $$
condition (ii) is equivalent to 
\begin{equation}
\label{eq:gxparn}
\sup_n d(g_nx,P(\hat\tau_n,\tau)) <+\infty .
\end{equation}
The implication (ii)$\Ra$(i)
thus follows from the previous lemma. 
The 
reverse implication (i)$\Ra$(ii) is easy:
Since 
$\sup_n d(x,P(\hat\tau_n,\tau))<+\infty$,
compare the sublemma,
the cone $V(x,\st(\tau))$ is contained in uniform tubular neighborhoods 
of all parallel sets $P(\hat\tau_n,\tau)$,
and conical convergence implies the same for the sequence $(g_nx)$,
i.e.\ (\ref{eq:gxparn}) is satisfied.
\qed

\medskip
Combining the lemma
with our earlier dynamical characterization of flag convergence,
see Lemma~\ref{lem:flconvitcontr},
we obtain:
\begin{prop}[Dynamical characterization of conical convergence]
A sequence $(g_n)$ in $G$ is $\taumod$-regular and $g_n\stackrel{con}{\to}\tau\in\Flagt$ 
iff  there exists a bounded sequence $(b_n)$ in $G$ 
and a simplex $\tau_-\in\Flagit$
such that the following conditions are satisfied:

(i) 
$b_ng_n^{-1}|_{C(\tau)}\to\tau_-$ uniformly on compacts.

(ii)
The accumulation set of the sequence $(b_ng_n^{-1}\tau)$ in $\Flagt$ is contained in $C(\tau_-)$.
\end{prop}
\proof
Suppose first that $(g_n)$ is $\taumod$-regular and $g_n\stackrel{con}{\to}\tau\in\Flagt$.
Then we have in particular flag convergence $g_n\to\tau$,
and Lemma~\ref{lem:flconvitcontr} yields $(b_n)$ and $\tau_-$ with (i).
The conical convergence $g_n\stackrel{con}{\to}\tau$ is equivalent to 
$g_nb_n^{-1}\stackrel{con}{\to}\tau$,
and so the previous lemma implies for any $\hat\tau\in C(\tau)$ that 
the sequence $b_ng_n^{-1}(\hat\tau,\tau)$ is relatively compact in $(\Flagit\times\Flagt)^{opp}$.
Since $b_ng_n^{-1}\hat\tau\to\tau_-$ by (i),
the sequence $(b_ng_n^{-1}\tau)$
therefore cannot accumulate at points outside $C(\tau_-)$.

Suppose now vice versa that $(b_n)$ and $\tau_-$ with (i+ii) are given.
By Lemma~\ref{lem:flconvitcontr}, (i) implies that $(g_n)$ is $\taumod$-regular and $g_n\to\tau$,
and the same follows for the sequence $(g_nb_n^{-1})$.
Furthermore,
(i+ii) imply that
for any $\hat\tau\in C(\tau)$
the sequence $b_ng_n^{-1}(\hat\tau,\tau)$ is relatively compact in $(\Flagit\times\Flagt)^{opp}$.
Thus $g_nb_n^{-1}\stackrel{con}{\to}\tau$ by the previous lemma,
and hence $g_n\stackrel{con}{\to}\tau$.
\qed

\medskip
We deduce the following criterion for being a conical limit simplex of a subsequence:
\begin{cor}
\label{cor:recogncnlim}
A sequence $(g_n)$ in $G$ has a 
$\taumod$-regular 
subsequence $\taumod$-conically converging to $\tau\in\Flagt$
iff  there exists a subsequence $(g_{n_k})$ and a simplex $\tau_-\in\Flagit$
such that the following conditions are satisfied:

(i) 
$g_{n_k}^{-1}|_{C(\tau)}\to\tau_-$ uniformly on compacts.

(ii)
$(g_{n_k}^{-1}\tau)$ converges to a simplex in $C(\tau_-)$.
\end{cor}
\proof
Suppose that there is a $\taumod$-regular subsequence $(g_{n_k})$ with $g_{n_k}\stackrel{con}{\to}\tau$.
The proposition yields a bounded sequence $(b_k)$ and $\tau_-$ such that properties (i+ii) in the proposition 
are satisfied for the sequence $(b_kg_{n_k}^{-1})$. 
After extraction, we obtain convergence $b_k\to b$ in $G$ 
and $b_kg_{n_k}^{-1}\tau\to\hat\tau_-\in C(\tau_-)$ in $\Flagt$.
The asserted properties (i+ii) then result from replacing $\tau_-$ with $b^{-1}\tau_-$. 
The converse is immediate in view of the proposition. 
\qed

\medskip
Now we turn to subgroups.

\begin{dfn}[Conical limit set]
For a subgroup $\Ga<G$,
a limit simplex $\la\in\LatGa$ is {\em $\taumod$-conical}
if there exists a $\taumod$-regular sequence $(\ga_n)$ in $\Ga$
such that $\ga_n\stackrel{con}{\to}\la$.
The {\em conical $\taumod$-limit set} 
$\LatGacon\subseteq\LatGa$ 
is the subset of conical limit simplices.
The subgroup $\Ga$ {\em has conical $\taumod$-limit set} or {\em is $\taumod$-conical} 
if all limit simplices are conical, $\LatGacon=\LatGa$.
\end{dfn}

We restrict ourselves to $\taumod$-antipodal $\taumod$-regular subgroups
and assume in particular that $\taumod$ is $\iota$-invariant.
Recall that then the action $$\Ga\acts\LatGa$$ is a {\em convergence action},
see section~\ref{sec:antip}.
This raises the question how the $\taumod$-conicality of limit simplices 
compares to their intrinsic conicality with respect to this convergence action,
cf.\ section~\ref{sec:convdy}.
We show that these properties are equivalent:
\begin{prop}
[Conical versus intrinsically conical limit simplex]
\label{prop:concon}
Let $\Ga<G$ be a $\taumod$-antipodal $\taumod$-regular
subgroup with $|\LatGa|\geq3$.
Then a limit simplex in $\LatGa$ is conical 
iff it is intrinsically conical 
for the convergence action $\Ga\acts\LatGa$.
\end{prop}
\proof
That conicality implies intrinsic conicality is,
in view of the corollary,
an immediate consequence of antipodality and 
Lemma~\ref{lem:recogncnlimcv}.

Suppose that, conversely, $\la\in\LatGa$ is intrinsically conical.
Again invoking Lemma~\ref{lem:recogncnlimcv},
this means that there exist a sequence $(\ga_n)$ in $\Ga$ and a limit simplex $\la_-\in\LatGa$ 
such that 
$\ga_n^{-1}|_{\LatGa-\{\la\}}\to\la_-$ uniformly on compacts and 
$\ga_n^{-1}\la\to\hat\la_-\in\LatGa-\{\la_-\}\subset C(\la_-)$.
On the other hand,
since $\Ga$ is a $\taumod$-convergence subgroup,
after extraction, the sequence $(\ga_n^{-1})$ becomes $\taumod$-contracting 
and there are limit simplices $\la',\la'_-\in\LatGa$ such that 
$\ga_n^{-1}|_{C(\la')}\to\la'_-$ uniformly on compacts.
In view of antipodality,
$C(\la')$ contains $\LatGa-\{\la'\}$.
Since $|\LatGa|\geq3$, it follows that 
$C(\la')$ intersects $\LatGa-\{\la\}$ and therefore $\la'_-=\la_-$.
Moreover,
from $\ga_n^{-1}\la\to\hat\la_-\neq\la_-$
it follows that $\la\not\in C(\la')$ and hence also $\la'=\la$.
We conclude that 
$\ga_n^{-1}|_{C(\la)}\to\la_-$ uniformly on compacts and 
$\ga_n^{-1}\la\to\hat\la_-\in C(\la_-)$.
Corollary~\ref{cor:recogncnlim} now yields that the limit simplex $\la$ is $\taumod$-conical. 
\qed

\begin{cor}[Conical versus intrinsically conical subgroup]
\label{cor:concon}
Let $\Ga<G$ be a $\taumod$-antipodal $\taumod$-regular
subgroup with $|\LatGa|\geq3$.
Then $\Ga$ is $\taumod$-conical 
iff  all simplices in $\LatGa$ are conical limit points 
for the convergence action $\Ga\acts\LatGa$.
\end{cor}

We introduce the following asymptotic condition on the orbit geometry of
subgroups:
\begin{dfn}[RCA]
\label{def:rca}
A subgroup $\Ga<G$ is {\em $\taumod$-RCA}
if it is $\taumod$-regular, $\taumod$-conical and $\taumod$-antipodal.
\end{dfn}

From the corollary we deduce, 
using the dynamical characterization of word hyperbolic groups and their boundary actions,
the following equivalence:

\begin{thm}
\label{thm:rcasmbeq}
For a 
subgroup $\Ga<G$ with $|\LatGa|\geq3$
the following properties are equivalent:

(i) $\taumod$-RCA

(ii) $\taumod$-asymptotically embedded

The implication (ii)$\Ra$(i) holds without restriction on the size of the limit set.
\end{thm}
\proof
Since this is part of both conditions,
we assume that $\Ga$ is $\taumod$-regular and $\taumod$-antipodal. 

The implication (ii)$\Ra$(i) follows, without restriction on the size of $\LatGa$,
from the implication (i)$\Ra$(iii) of Theorem~\ref{thm:mrsasembeq}.

Suppose now that $|\LatGa|\geq3$.
According to the previous corollary,
the subgroup $\Ga$ is $\taumod$-RCA if and only if  the convergence action $\Ga\acts\LatGa$ is (intrinsically) conical.
In view of Theorems~\ref{thm:bowdchr} and~\ref{thm:hypgpbdac}
this is equivalent to $\Ga$ being word hyperbolic 
and $\LatGa$ being $\Ga$-equivariantly homeomorphic to $\geo\Ga$,
i.e.\ to $\Ga$ being $\taumod$-asymptotically embedded.
\qed

\subsection{Subgroups with two-point limit sets}

For antipodal regular subgroups with two-point limit sets,
some of our conditions are automatically satisfied:

\begin{lem}
\label{lem:2limrca}
Suppose that $\Ga<G$ is $\taumod$-antipodal $\taumod$-regular
with $|\LatGa|=2$.
Then: 

(i) $\Ga$ is $\taumod$-RCA,

(ii) $\Ga$ is virtually cyclic,

(iii) The orbit maps $o_x:\Ga\to\Ga x\subset X$ extend continuously to infinity by an asymptotic embedding.
In particular, $\Ga$ is $\taumod$-asymptotically embedded.
\end{lem}
\proof
(i) By antipodality, $\LatGa$ consists of a pair of opposite simplices $\la_{\pm}\in\Flagt$.
The subgroup $\Ga$ therefore preserves the parallel set $P(\la_-,\la_+)$.
The limit simplices $\la_{\pm}$ must be conical by Lemma~\ref{lem:frconcnv}.
Hence $\Ga$ is $\taumod$-RCA.

(ii) Pick a point $x\in P(\la_-,\la_+)$.
By conicality,
there exists an element $\ga_0\in\Ga$ which fixes $\la_{\pm}$ 
and so that $\ga_0x$ lies in the interior of the Weyl cone $V=V(x,\st(\la_+))\subset P(\la_-,\la_+)$.
We consider the biinfinite nested sequence of Weyl cones $\ga_0^nV$ for $n\in\Z$.
The cones $\ga_0^nV$ cover $P(\la_-,\la_+)$,
cf.\ Proposition~\ref{prop:hdstwcones}. 
Moreover,
$\ga_0^{n+1}V$ is contained in the interior of $\ga_0^nV$ 
and has finite Hausdorff distance from it.
By regularity, 
the difference of cones $V-\ga_0V$
can only contain finitely many points of the orbit $\Ga x$.
The corresponding elements in $\Ga$
form a set of representatives for the cosets of the infinite cyclic subgroup $\Ga_0$ generated by $\ga_0$ in $\Ga$.
Hence $\Ga$ is virtually cyclic.

(iii) Since $\ga_0^{\pm n}\to\la_{\pm}$ as $n\to+\infty$,
the restrictions of the orbit maps to $\Ga_0$ extend continuously to $\geo\Ga_0\cong\geo\Ga$
by an asymptotic embedding $\al$. 
Since $\Ga_0$ has finite index in $\Ga$, 
the map $\al$ is a continuous extension also of the orbit maps of $\Ga$ itself. 
Moreover, it is $\Ga$-equivariant.
\qed

\subsection{Expansion}
\label{sec:expa}

We define another purely dynamical condition for subgroups,
inspired by Sullivan's notion of {\em expanding actions} \cite{Sullivan},
namely that their action on the appropriate flag manifold is expanding at the limit set
in the sense of Definition~\ref{def:metexpan}. 
As before, 
we equip the flag manifolds with auxiliary Riemannian metrics.

\begin{dfn}[CEA]
\label{def:cea}
A subgroup $\Ga<G$ 
is {\em $\taumod$-CEA}
(convergence, expanding, antipodal)
if it is $\taumod$-convergence, $\taumod$-antipodal 
and the action $\Ga\acts\Flagt$ is expanding at $\LatGa$.
\end{dfn}

The next result relates conicality to infinitesimal expansion,
cf.\ Definition~\ref{def:infexpan}.
For smooth actions on Riemannian manifolds,
metric and infinitesimal expansion are equivalent.
\begin{lem}[Expansion at conical limit simplices]
\label{lem:expconlim}
Let $(g_n)$ be a $\taumod$-regular sequence in $G$ such that 
$g_n\stackrel{con}{\to}\tau\in\Flagt$.
Then the inverse sequence $(g_n^{-1})$ has diverging infinitesimal expansion on $\Flagt$ at $\tau$,
i.e.\ 
$$ \eps(g_n^{-1},\tau)\to+\infty$$
\end{lem}
\proof
This follows from the expansion estimate in Theorem~\ref{thm:expand}.
\qed

\medskip
Applied to subgroups, the lemma yields:
\begin{prop}[Conical implies expansive]
\label{prop:conicimplexpatau}
Let $\Ga<G$ be a subgroup. 
If $\la\in\LatGacon$,
then the action $\Ga\acts\Flagt$ 
has diverging infinitesimal expansion at $\la$. 

In particular, 
if $\Ga$ is $\taumod$-conical,
then $\Ga\acts\Flagt$ is expanding at $\LatGa$.
\end{prop}
\proof
This is a direct consequence of the lemma,
together with the fact that infinitesimal expansion implies metric expansion. 
\qed

\medskip
We obtain the equivalence of conditions:
\begin{thm}
\label{thm:rceaq}
For a 
subgroup $\Ga<G$ with $|\LatGa|\geq2$,
the following properties are equivalent:

(i) $\taumod$-RCA

(ii) $\taumod$-CEA

The implication (i)$\Ra$(ii) holds without restriction on the size of the limit set.
\end{thm}
\proof
We recall that $\taumod$-regularity is equivalent to the $\taumod$-convergence property,
cf.\ Theorem~\ref{thm:regimplcontrgp}.
Thus either condition implies that $\Ga$ is $\taumod$-regular and $\taumod$-antipodal. 

The implication (i)$\Ra$(ii) is the previous proposition. 
(We do not need that $|\LatGa|\geq2$.)

For the direction 
(ii)$\Ra$(i)
we first assume that $|\LatGa|\geq3$ and 
consider the convergence action $\Ga\acts\LatGa$.
Since $\LatGa$ contains at least three points, it must be perfect\footnote{I.e.\ has no isolated points.}
(see \cite[Thm.\ 2S]{Tukia_convgps}).
By assumption, the action $\Ga\acts\LatGa$ is expanding.
Therefore all points $\la\in\LatGa$ are intrinsically conical, cf.\ Lemma~\ref{lem:conical}, 
and hence (extrinsically) conical, 
i.e.\ $\Ga$ is $\taumod$-conical, cf.\ Corollary~\ref{cor:concon}.

In the case $|\LatGa|=2$,
the assertion follows from Lemma~\ref{lem:2limrca}.
\qed

\subsection{Anosov property}
\label{sec:anosov}

The Anosov condition combines boundary embeddedness 
with an infinitesimal expansion condition at the image of the boundary embedding:

\begin{dfn}[Anosov]
\label{dfn:ouranosov}
A subgroup $\Ga<G$ is {\em $\taumod$-Anosov} if:

(i) $\Ga$ is $\taumod$-boundary embedded with boundary embedding $\beta$.

(ii) For every ideal point $\zeta\in \geo \Ga$
and every normalized 
(by $r(0)=e\in \Ga$)
discrete geodesic ray $r: \N\to \Ga$ asymptotic to $\zeta$,
the action $\Ga\acts\Flagt$ satisfies
\begin{equation*}
\eps(r(n)^{-1}, \beta(\zeta))\ge A e^{Cn}
\end{equation*}
for $n\geq 0$ with constants $A, C>0$  independent of $r$. 
\end{dfn}

We recall that boundary embedded subgroups are {\em discrete}.

Our notion of $\taumod$-Anosov is equivalent to the notion of $P$-Anosov in \cite{GW}
where $P<G$ is a parabolic subgroup in the conjugacy class corresponding to $\taumod$,
see section~\ref{sec:onos}.
We note also that the study of $(P_+,P_-)$-Anosov subgroups quickly reduces 
to the case of $P$-Anosov subgroups by intersecting parabolic subgroups,
cf.\ \cite[Lemma 3.18]{GW}.

In both our and the original definition uniform exponential expansion rates are required.
We will see that the conditions can be relaxed without altering the class of subgroups.
Uniformity can be dropped, and instead of exponential divergence 
the mere unboundedness of the expansion rate suffices.
\begin{dfn}[Non-uniformly Anosov]
\label{def:nuanos}
A subgroup $\Ga<G$ is {\em non-uniformly $\taumod$-Anosov} if:

(i) $\Ga$ is $\taumod$-boundary embedded with boundary embedding $\beta$.

(ii) For every ideal point $\zeta\in \geo \Ga$
and every normalized\footnote{Here, the normalization can be dropped because no {\em uniform} growth is required.}
discrete geodesic ray $r: \N_0\to \Ga$ asymptotic to $\zeta$,
the action $\Ga\acts\Flagt$ satisfies 
\begin{equation}
\label{eq:nuanosexp}
\sup_{n\in\N}\, \eps(r(n)^{-1}, \beta(\zeta))=+\infty. 
\end{equation}
\end{dfn}
In other words,
we require that for every ideal point $\zeta\in\geo\Ga$
the expansion rate $\eps(\ga_n^{-1},\beta(\zeta))$
non-uniformly diverges along some sequence $(\ga_n)$ in $\Ga$
which converges to $\zeta$ {\em conically}.

\medskip
We relate the Anosov to the Morse property,
building on our discussion of the coarse extrinsic geometry of subgroups in sections~\ref{sec:asyemb} and~\ref{sec:morse}.

\begin{thm}[Non-uniformly Anosov implies Morse]
\label{thm:nuanosasemb}
Each non-uniformly $\taumod$-Anosov subgroup $\Ga< G$ 
is $\taumod$-Morse. 

Moreover, the boundary 
embedding $\beta$ of $\Ga$ sends $\geo \Ga$ homeomorphically 
onto $\Lat(\Ga)$. 
\end{thm}
\proof
Let $\Ga<G$ be non-uniformly $\taumod$-Anosov.
Since non-uniformly Anosov subgroups are boundary embedded by definition,
discrete geodesic lines in $\Ga$ are mapped into uniform neighborhoods of $\taumod$-parallel sets
prescribed by the boundary embedding,
see Lemma~\ref{lem:lnclpar}.
The same follows for discrete geodesic rays in $\Ga$ because they lie in uniform neighborhoods of lines,
compare the proof of Lemma~\ref{lem:raydivdrif}: 
For every ray $r:\N_0\to\Ga$ asymptotic to $\zeta=r(+\infty)$ 
there exists an ideal point $\hat\zeta\in\geo\Ga-\{\zeta\}$ 
such that the path $rx$ lies in the $\rho''(\Ga,x)$-neighborhood of the parallel set 
$P=P(\beta(\hat\zeta),\beta(\zeta))$.
Here, as usual, $x\in X$ is some fixed base point. 

The expansion condition (\ref{eq:nuanosexp}) further restricts the position of the path $rx$
along the parallel set: 
Let $x_n\in P$ denote points at distance $\leq\rho''$ from the points $r(n)x$,
e.g.\ their nearest point projections to $P$.
For a strictly increasing sequence $n_k\to+\infty$ with diverging expansion rate 
$$\eps(r(n_k)^{-1}, \beta(\zeta))\to+\infty$$
we have in view of Proposition \ref{prop:expand} and Theorem~\ref{thm:expand} that 
$x_{n_k}\in V(x_0,\st(\beta(\zeta)))$ for large $k$ and 
$$ d\bigl(x_{n_k},\D V(x_0,\st(\beta(\zeta)))\bigr)\to+\infty$$
(non-uniformly) as $k\to+\infty$. 
Fix a constant $d>>\rho''$.
It follows that there exists a smallest ``entry time'' $T=T(r)\in\N$ such that 
the point $r(T)x$ lies in the open $3\rho''$-neighborhood of the cone $V(r(0)x,\st(\beta(\zeta)))$ 
and has distance $>d$ from its boundary.

We observe next that $T(r')\leq T(r)$ for rays $r'$ sufficiently close to $r$,
because $\zeta$ varies continuously with $r$,
and rays sufficiently close to $r$ agree with $r$ up to time $T(r)$.
Thus, $T$ is locally bounded above as a function of $r$.
Since $\Ga$ acts cocompactly on rays, equivalently,
since the space of rays with fixed initial point is compact,
we conclude that $T$ is {\em bounded above globally}, 
i.e.\ there exists a number $T_0=T_0(\Ga,x,d)$ such that $T(r)\leq T_0$ for all rays $r$.

As a consequence, 
for every ray $r$ the above sequence of natural numbers $(n_k)$ can be chosen 
with bounded increase 
$n_{k+1}-n_k\leq T_0$
and so that 
$$x_{n_{k+1}}\in V(x_{n_k},\st(\beta(\zeta)))$$
and 
$$ d\bigl(x_{n_{k+1}},\D V(x_{n_k},\st(\beta(\zeta)))\bigr)>\frac{d}{2}$$
for all $k$,
i.e.\ the sequence $(n_k)$ increases uniformly linearly 
and the Weyl cones $V(x_{n_k},\st(\beta(\zeta)))$ are uniformly nested,
compare the proof of Theorem~\ref{thm:asymbur}.

It follows that the paths $rx$ are uniformly $\taumod$-regular and undistorted,
and are contained in uniform neighborhoods of the cones $V(r(0),\st(\beta(r(+\infty))))$.
In particular, $\Ga$ satisfies property (iii) of Theorem~\ref{thm:mrsasembeq},
and therefore is $\taumod$-Morse. 
It also follows that $\beta(\geo \Ga)\subseteq \Lat(\Ga)$. 
The equality $\beta(\geo \Ga)=\Lat(\Ga)$ follows from Theorem \ref{thm:bemblim}. 
\qed

\medskip
A converse readily follows from our earlier results:
\begin{thm}
\label{thm:mrsimplanos}
$\taumod$-Morse subgroups $\Ga<G$ are $\taumod$-Anosov. 
\end{thm}
\proof
Let $\Ga<G$ be $\taumod$-Morse.
By Theorems~\ref{thm:mrsasembeq} and~\ref{thm:asymbur},
$\Ga$ is then also $\taumod$-asymptotically embedded and uniformly $\taumod$-regular.
Furthermore, 
denoting the asymptotic embedding by $\al$ and fixing a point $x\in X$, 
we know that for every ray $r:\N_0\to\Ga$ 
the path $rx$ is contained in a uniform neighborhood of the Weyl cone $V(r(0)x,\al(r(+\infty)))$
and drifts away from its boundary at a uniform linear rate.
With Theorem~\ref{thm:expand} it follows that 
the infinitesimal expansion factor $\eps(r(n)^{-1},\al(r(+\infty)))$
for the action $\Ga\acts\Flagt$ 
grows at a uniform exponential rate.
Thus, $\Ga$ is $\taumod$-Anosov.
\qed

\subsection{Equivalence of conditions}

Combining our results 
comparing the various geometric and dynamical conditions for discrete subgroups,
we obtain:
\begin{thm}[Equivalence]
\label{thm:eqv}
The following properties for subgroups $\Ga<G$ 
are equivalent in the nonelementary\footnote{Meaning 
that  $|\LatGa|\geq3$ in (i), (ii), (v), (vi) and that $\Ga$ is  word hyperbolic  with $|\geo \Ga|\ge 3$ in (iii), (iv).} 
case:

(i)  $\taumod$-asymptotically embedded

(ii) $\taumod$-CEA

(iii) $\taumod$-Anosov

(iv) non-uniformly $\taumod$-Anosov

(v) $\taumod$-RCA

(vi) $\taumod$-Morse. 

These properties imply $\taumod$-URU. 

Moreover, the boundary maps in  (i), (iii) and (iv) coincide.
\end{thm}

\proof
By Theorem~\ref{thm:mrsasembeq}, (i) and (vi) are equivalent. 
By Theorems~\ref{thm:nuanosasemb} and~\ref{thm:mrsimplanos},
conditions (iii), (iv) and (vi) are equivalent. 
The fact that the boundary maps in (i), (iii) and (iv) coincide
follows from the second part of Theorem \ref{thm:nuanosasemb}. 

By Theorem~\ref{thm:asymbur}, 
(i) implies $\taumod$-URU.
By Theorem~\ref{thm:rcasmbeq}, (i) and (v) are equivalent.
By Theorem~\ref{thm:rceaq}, (ii) and (v) are equivalent. 
\qed

\begin{rem}
(i)
The equivalence of the conditions (i), (iii), (iv) and (vi), the fact that they imply $\taumod$-URU,
and the implications (i)$\Ra$(v)$\Ra$(ii)
hold without restriction on the size of the limit set.

(ii) 
It is shown in \cite{mlem} that, conversely, $\taumod$-URU implies $\taumod$-Morse.
\end{rem}

For subgroups with small limit sets we have the following additional information,
see Lemma~\ref{lem:2limrca}:

\begin{add}
For a $\taumod$-antipodal $\taumod$-regular subgroup $\Ga<G$ with $|\LatGa|=2$,
properties (i)-(vi) and $\taumod$-URU are always satisfied.
\end{add}

We are unaware of examples of $\taumod$-RCA or $\taumod$-CEA subgroups with one limit point
in higher rank.
Note that such subgroups cannot be $\taumod$-asymptotically embedded.

\subsection{Morse quasigeodesics}
\label{sec:mqg}

When studying the coarse geometry of Anosov subgroups 
in sections~\ref{sec:asyemb} and~\ref{sec:morse},
we were lead to the Morse and URU properties.
We also saw that Morse implies URU.
(The converse is true as well, but harder to prove, see \cite{mlem}.)

Thus, for Morse subgroups $\Ga<G$, 
the images of the discrete geodesics in $\Ga$ under an orbit map are uniform quasigeodesics in $X$ 
which are uniformly regular 
and satisfy a Morse type property involving closeness of subpaths to diamonds. 
Leaving the group-theoretic context, we will now make this class of quasigeodesics precise 
and study some of its geometric properties.
(See also \cite{morse} for further discussion.)
We will build in the uniform regularity into the Morse property 
by replacing the diamonds with smaller ``uniformly regular'' $\Theta$-diamonds.

In the following, 
$\Theta\subset\inte_{\taumod}(\simod)$ denotes 
an $\iota$-invariant $\taumod$-Weyl convex
compact subset
which is used to quantify uniform regularity.
We work with discrete paths;
$I\subseteq\R$ denotes an interval and 
$n_{\pm}$ integers.

\begin{definition}[Morse quasigeodesic]
\label{dfn:mqg}
A quasigeodesic $q:I\cap\Z\to X$ is {\em $(\Theta,\rho)$-Morse}
if for every subinterval $[n_-,n_+]\subseteq I$
the subpath $q|_{[n_-,n_+]\cap\Z}$
is contained in the $\rho$-neighborhood 
of a diamond $\diamoTh(x_-,x_+)$ 
with tips at distance $d(x_{\pm},q(n_{\pm}))\leq\rho$ from the endpoints. 

We say that an infinite quasigeodesic is {\em $\Theta$-Morse} if it is $(\Theta,\rho)$-Morse for some $\rho$,
and we say that it is {\em $\taumod$-Morse} if it is $\Theta$-Morse for some $\Theta$.
\end{definition}

The $\Theta$-Morse property for quasigeodesics is clearly stable under bounded perturbation.

We say that some paths are {\em uniform $\taumod$-Morse quasigeodesics}
if they are uniform quasigeodesics\footnote{I.e.\ quasigeodesics with the same quasiisometry constants.}
and $(\Theta,\rho)$-Morse with the same $\Theta,\rho$.

\medskip
We can now interpret the Morse subgroup property in terms of Morse quasigeodesics:

\begin{prop}
\label{prop:msbgpmqg}
An intrinsically word hyperbolic subgroup $\Ga< G$ is $\taumod$-Morse
if and only if an orbit map $o_x:\Ga\to\Ga x\subset X$ sends uniform quasigeodesics in $\Ga$
to uniform $\taumod$-Morse quasigeodesics in $X$.
\end{prop}
\proof Suppose that $\Ga$ is $\taumod$-Morse.
We fix a word metric on $\Ga$.
In view of the Morse Lemma for word hyperbolic groups (Gromov hyperbolic spaces) 
it suffices to prove that $o_x$
sends discrete geodesics in $\Ga$ to uniform $\taumod$-Morse quasigeodesics in $X$.

First of all,
since Morse subgroups are URU, 
we know that $\Ga$ is undistorted in $G$, i.e.\ $o_x$ is a quasiisometric embedding. 
Equivalently, the $o_x$-images of discrete geodesics in $\Ga$ are uniform quasigeodesics.
We need to show that they are uniformly $\taumod$-Morse.

Consider a discrete geodesic segment $s: [n_-, n_+]\cap\Z\to \Ga$.
According to the Morse subgroup property of $\Ga$, 
the image path $sx=o_x\circ s$
is contained in a tubular neighborhood of uniform radius $\rho''=\rho''(\Ga,x)$
of a diamond $\diamot(x_-,x_+)$ with 
$d(x_{\pm},s(n_{\pm})x)\leq\rho''$.
It will be enough to verify 
that $sx$ is also contained in a uniform tubular neighborhood of the smaller $\Theta$-diamond $\diamo_\Theta(x_-,x_+)$ 
for some $\Theta$ independent of $s$. 

For $n_-\leq n\leq n_+$,
let $p_n\in \diamot(x_-, x_+)$ 
denote the nearest point projection of $s(n)x$.
In view of the uniform upper bound $\rho''$ for the distances 
$d(x_{\pm},s(n_{\pm})x)$ and $d(p_n,s(n)x)$,
the uniform regularity of $\Ga$ implies:
If $n- n_-, n_+ -  n \geq C_0$ (with a uniform constant $C_0$), 
then 
$$
d_\Delta(x_\pm, p_n)\in V(0, \Theta)
$$
with a compact $\Theta\subset \inte_{\taumod}(\simod)$ independent of $s$.
Moreover,
after enlarging $\Theta$,
we may assume that it is $\iota$-invariant and $\taumod$-Weyl convex.
It follows that the diamond $\diamoTh(x_-,x_+)$ is defined and 
$p_n\in \diamoTh(x_-,x_+)$.
Hence, $sx$ is contained in a uniform tubular neighborhood of $\diamoTh(x_-,x_+)$.

Conversely, 
suppose that $o_x$ sends discrete geodesics in $\Ga$ to uniform $\taumod$-Morse quasigeodesics in $X$.
Then $\Ga$ is undistorted and the geodesic segments with endpoints in the orbit $\Ga x$
are uniformly close to $\Theta$-regular segments, 
equivalently,
the $\De$-distances $d_{\De}(x,\ga x)$ between orbit points 
are contained in a tubular neighborhood of the cone $V(0,\Theta)$.
It follows that $\Ga$ is (uniformly) $\taumod$-regular,
and hence $\taumod$-Morse.
\qed

\medskip
Next, we briefly discuss the {\em asymptotics} of infinite Morse quasigeodesics.
There is much freedom for the asymptotic behavior 
of arbitrary quasigeodesics in euclidean spaces,
and therefore also in symmetric spaces of higher rank.  
However, 
the asymptotic behavior of Morse quasigeodesics 
is as restricted as for quasigeodesics in rank one symmetric spaces.

Morse quasirays satisfy a version of the defining property for Morse quasigeodesic segments,
with diamonds replaced by cones.
As a consequence,
although Morse quasirays in general do not 
converge at infinity in the visual compactification, 
they {\em flag converge}:

\begin{lem}[Conicality of Morse quasirays]
\label{lem:morsecon}
A $(\Theta,\rho)$-Morse quasiray $q:\N_0\to X$ 
is contained in the $\rho$-neighborhood of a $\Theta$-cone $V(x,\stTh(\tau))$ with $d(x,q(0))\leq\rho$ 
for a unique simplex $\tau\in\Flagt$.
Furthermore, $q(n)\to\tau$ conically.
\end{lem}
\proof
The existence of the cone $V(x,\stTh(\tau))$ 
follows from the definition of Morse quasigeodesics by a limiting argument.
Obviously, we have conical $\taumod$-flag convergence $q(n)\to\tau$,
which also implies the uniqueness of $\tau$.
\qed

\medskip
Now we give a {\em Finsler geometric characterization} of Morse quasigeodesics.
We show that they are the coarsification of (uniformly regular) Finsler geodesics (cf.\ Definition~\ref{def:finsgeo}).
Even though this is true in general,
we will give the proof only in the infinite case (of rays and lines),
since it is simpler and suffices for the purposes of this paper:

\begin{thm}
[Morse quasigeodesics are uniformly close to Finsler geodesics] 
\label{thm:mqgclfins}
Uniform $\taumod$-Morse quasigeodesic rays and lines are uniformly {Hausdorff} close 
to uniformly $\taumod$-regular $\taumod$-Finsler geodesic rays and lines.
\end{thm}
\proof 
It suffices to treat the ray case.
The line case follows by a limiting argument. 

Let $q:\N_0\to X$
be a $(\Theta,\rho)$-Morse quasigeodesic ray.
According to Lemma \ref{lem:morsecon}, 
$q$ is contained in a uniform tubular neighborhood of a Weyl cone $V=V(q(0),\st(\tau))$.
As in the proof that asymptotically embedded implies URU (Theorem~\ref{thm:asymbur}),
we consider the sequence of nearest point projections $x_n\in V$ of the points $q(n)$, $n\in\N_0$.
Again by Lemma \ref{lem:morsecon}, 
the point 
$x_{n+m}$ lies in a uniform tubular neighborhood of the $\Theta$-cone $V(x_n,\stTh(\tau))\subset V$ 
for all $n,m\geq0$.

We slightly enlarge $\Theta$ to $\Theta'$, 
such that $\Theta\subset\inte(\Theta')$ as subsets of $\inte_{\taumod}(\simod)$.
Then there exists $m_0\in\N$ depending on $\Theta,\Theta',\rho$ and the quasiisometry constants of $q$, 
such that 
$$x_{n+m}\in V(x_n,\st_{\Theta'}(\tau))$$
for all $n\geq0$ and $m\geq m_0$.
The piecewise geodesic path 
$$ x_0x_{m_0}x_{2m_0}x_{3m_0}\ldots$$
is then a $\Theta'$-regular $\taumod$-Finsler geodesic ray uniformly Hausdorff close to $q$.
\qed 

\medskip
We use the approximation of Morse quasigeodesics by Finsler geodesics 
to coarsify Theorem~\ref{thm:dedstalfins} 
and deduce an analogous result on the {\em $\De$-distance along Morse quasigeodesics}.
Again, we restrict ourselves to the infinite case of rays:
\begin{thm}[$\De$-projection of Morse quasirays]
\label{thm:dedstalmqg}
If $q:\N_0\to X$ is a $\taumod$-Morse quasiray, 
then so is
$$ \bar q_{\De}=d_{\De}(q(0),q):\N_0\to \De .$$
Moreover, uniform $\taumod$-Morse quasirays $q$ yield uniform $\taumod$-Morse quasirays $\bar q_{\De}$.
\end{thm}
\proof
Suppose that $q$ is a $(\Theta,\rho)$-Morse quasiray.
We enlarge $\Theta$ to $\Theta'$ such that $\Theta\subset\inte(\Theta')$.
According to the proof of Theorem~\ref{thm:mqgclfins},
there exists a $\Theta'$-regular $\taumod$-Finsler geodesic ray $c:[0,+\infty)\to X$
which is uniformly close to $q$ in terms of the data $\Theta,\Theta',\rho$ and the quasiisometry constants,
i.e.\ $d(c(n),q(n))$ is uniformly bounded.
In particular, $c$ is also a uniform quasiray.

For the $\De$-projections $\bar c_{\De}=d_{\De}(c(0),c)$ and $\bar q_{\De}$,
the pointwise distance $d(\bar c_{\De}(n),\bar q_{\De}(n))$ is also uniformly bounded.
According to Theorem~\ref{thm:dedstalfins}, 
$\bar c_{\De}$ is again a $\Theta'$-regular $\taumod$-Finsler geodesic ray
and a uniform quasiray.
It follows that $\bar q_{\De}$ is a $(\Theta',\rho')$-Morse quasiray with uniform $\rho'$
and uniform quasiisometry constants. 
\qed

\subsection{Appendix: The original Anosov definition}
\label{sec:onos}

A notion of Anosov representations of surface groups into $PSL(n,\R)$ was introduced by Labourie in \cite{Labourie},
and generalized to a notion of $(P_+,P_-)$-Anosov representations $\Ga\to G$ of word hyperbolic groups into semisimple Lie groups by 
Guichard and Wienhard in \cite{GW}. The goal of this section is to review this definition of Anosov representations $\Ga\to G$ using the language of {\em expanding and contracting flows} and then present a closely related and equivalent definition which avoids the language of flows. 

Let $\Gamma$ be a non-elementary (i.e.\ not virtually cyclic) word hyperbolic group with a fixed word metric $d_\Ga$ and Cayley graph $C_\Ga$.  Consider a {\em geodesic flow} $\widehat\Ga$ of $\Ga$;  such a  flow was originally constructed by Gromov \cite{Gromov_hypgps} and then improved by Champetier \cite{Champetier} and Mineyev \cite{Mineyev}, resulting in definitions with different properties. 
We note that the exponential convergence of asymptotic geodesic rays will not be used in our discussion; 
as we will see, it is also irrelevant whether the trajectories of the geodesic flow are geodesics or uniform quasigeodesics in $\widehat\Ga$. In particular, it will be irrelevant for us which definition of $\widehat\Ga$ is used. 
Only the following properties of $\widehat\Ga$ will be used in the sequel:

1. $\widehat\Ga$ is a proper metric space. 

2. There exists a properly discontinuous isometric action 
$\Ga\acts\widehat\Ga$. 

3. There exists a $\Ga$-equivariant quasi-isometry $\pi: \widehat\Ga\to \Ga$; in particular, the fibers of $\pi$ are relatively compact.

4. There exists a continuous action $\R\acts\widehat\Ga$, denoted $\phi_t$ and called the {\em geodesic flow}, 
whose trajectories are uniform quasigeodesics in $\widehat\Ga$, i.e. for each $\hat m\in \widehat\Ga$ the 
{\em flow line} 
$$
t\to \hat m_t := \phi_t(\hat m)
$$ 
is a uniform quasi-isometric embedding $\R\to \widehat\Ga$. 

5. The flow $\phi_t$ commutes with the action of $\Ga$. 

6. Each $\hat m\in\hat\Ga$ defines a uniform  
quasigeodesic $m:t\mapsto m_t$ in $\Ga$ by the formula: 
$$
m_t= \pi(\hat m_t) 
$$
Following the notation in section \ref{sec:convdy}, we let $(\geo\Ga\times\geo\Ga)^{dist}
$ denote the subset of $\geo\Ga\times\geo\Ga$ consisting of pairs of distinct points. 
The natural map 
$$
e=(e_-,e_+):\hat\Ga\to (\geo\Ga\times\geo\Ga)^{dist}
$$
assigning to $\hat m$ the pair of  ideal endpoints $(m_{-\infty},m_{+\infty})$ 
of $m$ is continuous and surjective. 
In particular,
every uniform quasigeodesic in $\hat\Ga$ 
is uniformly Hausdorff close to a flow line.

\medskip 
The reader can think of the elements of $\widehat\Ga$ {\em as parameterized geodesics in $C_\Ga$, so that $\phi_t$ acts on geodesics via reparameterization}. This was Gromov's original viewpoint, 
although not the one in \cite{Mineyev}.  

We say that $\hat m\in \widehat\Ga$ is {\em normalized} if $\pi(\hat m)=1\in\Ga$. 
Similarly, maps $q: \Z\to \Ga$, and $q: \N\to \Ga$ will be called {\em normalized} if $q(0)=1$. 
It is clear that every $\hat m\in \widehat\Ga$ can be sent to a normalized element of $\widehat\Ga$ 
via the action of $m_0^{-1}\in \Ga$.

Since trajectories of $\phi_t$ are uniform quasigeodesics, for 
each normalized $\hat m\in \widehat\Ga$ we have
\begin{equation}\label{eq:distance}
C_1^{-1} t - C_2\le d_\Ga(1, m_t) \le C_1 t + C_2
\end{equation}
for some positive constants $C_1, C_2$. 

Let ${\mathrm F}^\pm=\Flagpmt$ be a pair of {\em opposite} partial flag manifolds 
associated to the Lie group $G$, i.e.\ they are quotient manifolds of the form ${\mathrm F}^\pm=G/P_{\pm\taumod}$, 
see section \ref{sec:symmbas}. As usual, we will regard elements of ${\mathrm F}^\pm$ as simplices of type 
$\taumod, \iota\taumod$ in the Tits boundary of $X$.

Define the trivial bundles  
$$
E^\pm =\widehat \Ga \times {\mathrm F}^\pm \to \widehat \Ga. 
$$
For every representation $\rho: \Ga\to G$, the group $\Ga$ acts on both bundles via its natural action on $\widehat\Ga$ and via the representation $\rho$ on ${\mathrm F}^\pm$. 
Put a $\Ga$-invariant background Riemannian metric on the fibers of theses bundles, which varies  continuously with respect to $\hat m\in \widehat\Ga$. We will use the notation ${\mathrm F}^\pm_{\hat m}$ for the fiber above the point $\hat m$ equipped with this Riemannian metric. Since the subspace of $\widehat\Ga$ consisting of normalized elements is 
compact, it follows that for normalized $\hat m, \hat m'$ the identity map
$$
{\mathrm F}^\pm_{\hat m} \to {\mathrm F}^\pm_{\hat m'}
$$
is uniformly bilipschitz (with bilipschitz constant independent of $\hat m, \hat m'$). 
We will identify $\Ga$-equivariant (continuous) sections of the bundles $E^\pm$ with equivariant 
maps $s_\pm: \widehat\Ga\to {\mathrm F}^\pm$. 
These sections are said to be 
{\em parallel along flow lines}
if 
$$
s_\pm(\hat m)= s_\pm(\hat m_t)
$$
for all $t\in \R$ and $\hat{m} \in \widehat\Ga$. 

\begin{dfn}
Parallel sections $s_\pm$ are called 
{\em strongly parallel along flow lines} if for any two flow lines 
$\hat m,\hat m'$ with the same ideal endpoints, we have
$s_\pm(\hat m)=s_\pm(\hat m')$.
\end{dfn}

Note that this property is automatic 
for the geodesic flows constructed by Champetier and Mineyev
since (for their flows) 
any two flow lines which are at finite distance from each other 
are actually equal. 
Strongly parallel sections define $\Ga$-equivariant {\em boundary maps}
\begin{equation*}
\beta_\pm:\geo\Ga\to {\mathrm F}^\pm
\end{equation*}
from the Gromov boundary $\geo \Ga$ of the word hyperbolic group $\Ga$ by:
\begin{equation}
\label{eq:bdmapsect}
\beta_{\pm}\circ e_{\pm} = s_{\pm}~ . 
\end{equation}
\begin{lem}
The maps $\beta_{\pm}$ are continuous.
\end{lem}
\proof
Let $(\xi^n_-,\xi^n_+)\to(\xi_-,\xi_+)$
be a converging sequence in $(\geo\Ga\times\geo\Ga)^{dist}$.
There exists a bounded sequence $(\hat m^n)$ in $\hat\Ga$ 
such that $e_{\pm}(\hat m^n)=\xi^n_{\pm}$.
After extraction, the sequence $(\hat m^n)$ converges to some 
$\hat m\in \hat\Ga$.
Continuity of $s_{\pm}$ implies that 
$\beta_{\pm}(\xi^n_{\pm})=s_{\pm}(\hat m^n)\to s_{\pm}(\hat m)=\beta_{\pm}(\xi_{\pm})$.
This shows that no subsequence of $(\beta_{\pm}(\xi^n_{\pm}))$ can have a limit $\neq\beta_{\pm}(\xi_{\pm})$,
and the assertion follows from compactness of ${\mathrm F}^{\pm}$.
\qed

\medskip
Conversely, equivariant continuous maps $\beta_\pm$
define $\Ga$-equivariant sections $s_{\pm}$ strongly parallel along flow lines,  
by the formula (\ref{eq:bdmapsect}).

Consider the identity maps
$$
\Phi_{\hat m,t}: {\mathrm F}^\pm_{\hat m} \to {\mathrm F}^\pm_{\phi_t\hat m}. 
$$
These maps distort the Riemannian metric on the fibers. 
Using ~\eqref{eq:expfac},
we define the {\em infinitesimal expansion factor} 
of the flow $\phi(t)$ on the fiber  ${\mathrm F}^\pm_{\hat m}$ at the point $s_\pm(\hat m)$ as:
$$
\eps_\pm(\hat m, t) :=
\eps(\Phi_{\hat m,t}, s_\pm(\hat m))
$$

\begin{definition}\label{defn:uniformly expanding}
The geodesic flow $\phi_t$ is said to be {\em uniformly exponentially expanding} on the bundles $E^\pm$ with respect to the sections $s_\pm$ if there exist constants $a, c>0$ such that
$$
\eps_{\pm}(\hat m,\pm t) \ge a e^{ct}
$$
for all $\hat m\in \widehat \Ga$ and $t\geq0$.
\end{definition}

Our next goal is to give an alternative interpretation for the uniform expansion in this definition. 
First of all, since the metrics on the fibers are $\Ga$-invariant, it suffices to verify uniform exponential expansion only for normalized elements of $\widehat\Ga$. For a normalized element $\hat m\in \widehat\Ga$ and  $t\in\R$ consider the composition
$$
m_t^{-1} \circ \Phi_{\hat m,t}:  {\mathrm F}^\pm_{\hat m} \to {\mathrm F}^\pm_{m_t^{-1}\hat m_t}.
$$
Note that $\pi(m_t^{-1}\hat m_t)=m_t^{-1}m_t=1$,
i.e.\ both $\hat m$ and $m_t^{-1}\hat m_t$ are normalized.  
Since the group $\Ga$ acts isometrically on the fibers of the bundles $E^\pm$, the metric distortion of the above compositions is exactly the same as the distortion of $\Phi_{\hat m,t}$. Furthermore, since, as we noted above, the metrics on 
${\mathrm F}^\pm_{\hat m}$ and ${\mathrm F}^\pm_{m_t^{-1}\hat m_t}$ 
are uniformly bilipschitz to each other (via the ``identity'' map), the rate of expansion for the above composition (up to a uniform multiplicative error) is the same as the expansion rate for the map 
$$
\rho(m_t^{-1}): {\mathrm F}^\pm\to {\mathrm F}^\pm .
$$
(Here we are using fixed background Riemannian metrics on ${\mathrm F}^\pm$.) Thus, we get the estimate 
$$
C_3^{-1} \eps( \rho(m_t^{-1}), \beta_{\pm}(m_{\pm\infty}))\le    \eps_\pm(\hat m, t) \le 
C_3 \eps( \rho(m_t^{-1}), \beta_{\pm}(m_{\pm\infty})) 
$$
for some uniform constant $C_3>1$. By taking into account the equation \eqref{eq:distance}, we  obtain the following equivalent reformulation of Definition \ref{defn:uniformly expanding}:

\begin{lemma}\label{lem:expansion}
The geodesic flow is uniformly exponentially expanding with respect to the sections $s_\pm$ 
if and only if for every normalized uniform quasigeodesic 
$q: \Z\to \Ga$, which is asymptotic to points $\xi_\pm=q(\pm \infty)\in \geo \Ga$, the elements $\rho(q(\pm n))^{-1}$ act on $T_{\beta_\pm(\xi_\pm)} {\mathrm F}^\pm$ with uniform exponential expansion rate, i.e.
$$
\eps( \rho(q(\pm n))^{-1}, \beta_\pm(\xi_\pm))\ge Ae^{Cn}
$$
for all $\hat m\in \widehat \Ga$ and $n\ge 0$ 
with some fixed constants $A,C>0$.
\end{lemma}
\proof
There exists a normalized flow line $\hat m$ uniformly close to $q$, 
i.e.\ $q(n)$ is uniformly close to $m_{t_n}$ with 
$n\mapsto t_n$ being a uniform orientation-preserving quasiisometry $\Z\to\Z$.
Then $m_{\pm\infty}=\xi_\pm$, and
$\eps(\rho(q(\pm n))^{-1}, \beta_\pm(\xi_\pm))$
equals $\eps( \rho(m_{t_{\pm n}}^{-1}), \beta_{\pm}(m_{\pm\infty}))$
up to a uniform multiplicative error,
and hence also 
$\eps_\pm(\hat m, t_{\pm n})$.
\qed

\medskip
Since every uniform quasigeodesic ray in $\Ga$ 
extends to a uniform quasigeodesic line,
and in view of Morse lemma for hyperbolic groups, in the above definition it suffices to consider only normalized discrete 
geodesic rays $r: \N\to \Ga$.   

We can now  give the original and an alternative definition of Anosov representations. 
\begin{definition}
A pair of continuous maps $\beta_\pm: \geo \Ga \to {\mathrm F}^\pm$ is said to be {\em antipodal} if it satisfies the following conditions (called {\em compatibility} in \cite{GW}): 

(i) 
For every pair of distinct ideal points $\zeta, \zeta'\in \geo \Ga$, the simplices $\beta_+(\zeta)$, $\beta_-(\zeta')$ in the Tits boundary of $X$ are antipodal, equivalently, the corresponding parabolic subgroups of $G$ are opposite.  (In \cite{GW} this property is called {\em transversality}.) 

(ii) 
For every $\zeta\in \geo \Ga$, the simplices $\beta_+(\zeta), \beta_-(\zeta)$ belong to the same spherical Weyl chamber, i.e.\ the intersection of the corresponding parabolic subgroups of $G$ contains a minimal parabolic subgroup.
\end{definition}

Note that, as a consequence, the maps $\beta_{\pm}$ are {\em embeddings}, 
because antipodal simplices cannot be faces of the same chamber. 

\begin{definition}[\cite{GW}]
A representation $\rho: \Ga\to G$ is said to be $(P_{+\taumod}, P_{-\taumod})$-Anosov if there exists an antipodal pair of continuous $\rho$-equivariant 
maps $\beta_\pm: \geo \Ga \to {\mathrm F}^\pm$ such that the geodesic flow on the associated bundles $E^\pm$ satisfies the uniform expansion property with respect to the sections $s_\pm$  associated to the maps $\beta_\pm$. 
\end{definition}

The pair of maps $(\beta_+,\beta_-)$ in this definition is called {\em compatible} with the Anosov representation $\rho$.
Note that a $(P_{+\taumod}, P_{-\taumod})$-Anosov representation admits a unique compatible pair of maps. 
Indeed, 
the fixed points of infinite order elements $\ga\in\Ga$ 
are dense in $\geo\Ga$. 
The maps $\beta_{\pm}$ send the attractive and repulsive fixed points of $\ga$
to fixed points of $\rho(\ga)$ with contracting and expanding 
differentials,
and these fixed points are unique.
In particular, if $P_{+\taumod}$ is conjugate to $P_{-\taumod}$ (equivalently, $\iota\taumod=\taumod$)  
then $\beta_-=\beta_+$. 

We note that Guichard and Wienhard in  \cite{GW} use in their definition the uniform contraction property of the reverse flow $\phi_{-t}$ instead of the expansion property used above, but the two are clearly equivalent. Note also that  in the definition, it suffices to verify the uniform exponential expansion property only for the bundle $E_+$.
We thus obtain, as a corollary of Lemma \ref{lem:expansion}, the following alternative definition of Anosov representations:

\begin{prop}
[Alternative definition of Anosov representations] 
A representation $\rho: \Ga\to G$ is $(P_{+\taumod}, P_{-\taumod})$-Anosov if and only if there exists a pair 
of antipodal continuous $\rho$-equivariant maps $\beta_\pm: \geo \Ga \to {\mathrm F}^\pm$ such that for every normalized discrete geodesic ray 
$r: \N\to \Ga$ asymptotic to $\xi\in \geo \Ga$, the elements $\rho(r(n))^{-1}$ act on $T_{\beta_+(\xi)} {\mathrm F}_+$ with uniform exponential expansion rate, i.e.
\begin{equation}\label{eq:ano}
\eps(\rho(r(n))^{-1}, \beta_+(\xi))\ge A e^{Cn}
\end{equation}
for $n\geq0$
with  constants $A,C>0$ which are independent of $r$. 
\end{prop}

\medskip 
We now restrict to the case that the parabolic subgroups $P_{\pm\taumod}$ are conjugate to each other, 
i.e.\ the simplices $\iota\taumod=\taumod$. The $(P_{+\taumod}, P_{-\taumod})$-Anosov representations will in this case be called simply {\em $P_{\taumod}$-Anosov}, where $P_{\taumod}=P_{+\taumod}$, or simply 
$\tau_{mod}$-Anosov. Note that the study of general $(P_{+\taumod}, P_{-\taumod})$-Anosov representations quickly reduces to the case of $P$-Anosov representations by intersecting parabolic subgroups,
cf.\ \cite[Lemma 3.18]{GW}.
Now, 
$$
{\mathrm F}^{\pm}={\mathrm F}=G/P_{\taumod}=\Flagt
$$ 
and 
$$
\beta_{\pm}=\beta:\geo\Ga\to {\mathrm F} 
$$ 
is a single continuous embedding.
The compatibility condition 
reduces to the {\em antipodality} condition: 
For any two distinct ideal points $\xi,\xi'\in\geo \Ga$
the simplices $\beta(\xi)$ and $\beta(\xi')$ are antipodal to each other. 
In other words, $\beta$ is a {\em boundary embedding} in the sense of Definition~\ref{def:bdemb}.

We thus arrive to our definition,
compare Definition~\ref{dfn:ouranosov}:
\begin{dfn}[Anosov representation]
\label{defn:our anosov}
Let $\taumod$ be an $\iota$-invariant face of $\simod$.  
We call a representation $\rho: \Ga\to G$ {\em $P_{\taumod}$-Anosov} or $\tau_{mod}$-Anosov
if it is $\tau_{mod}$-boundary embedded 
with boundary embedding $\beta: \geo \Ga \to {\mathrm F}=\Flagt$ 
such that for every normalized discrete geodesic ray $r: \N\to \Ga$ asymptotic to $\zeta\in \geo \Ga$, 
the elements $\rho(r(n))^{-1}$ act on $T_{\beta(\zeta)} {\mathrm F}$ with uniform exponential expansion rate, i.e.\ 
$$
\eps(\rho(r(n))^{-1}, \beta(\zeta))\ge A e^{Cn}
$$ 
for $n\geq 0$ with constants $A, C>0$  independent of $r$. 
\end{dfn}

\noindent M.K.: Department of Mathematics, 
University of California, Davis, 
CA 95616, USA\\
email: kapovich@math.ucdavis.edu

\noindent B.L.: Mathematisches Institut, 
Universit\"at M\"unchen, 
Theresienstr. 39, 
D-80333 M\"unchen, Germany, 
email: b.l@lmu.de

\noindent J.P.: Departament de Matem\`atiques, 
Universitat Aut\`onoma de Barcelona, 
E-08193 Bellaterra, Spain, 
email: porti@mat.uab.cat


\begin{thebibliography}{BLP05}

\bibitem[Al]{Albuquerque}
P.\ Albuquerque, 
{\em Patterson-Sullivan theory in higher rank symmetric spaces}, 
Geom. Funct. Anal. Vol. {\bf 9} (1999), no. 1, p. 1--28. 

\bibitem[BGS]{BGS}
W. Ballmann, M. Gromov and V. Schroeder, ``Manifolds of nonpositive curvature'', Birkh\"auser Verlag, 1985. 

\bibitem[Be]{Benoist}
Y.\ Benoist,
{\em Propri\'et\'es asymptotiques des groupes lin\'eaires}, 
Geom. Funct. Anal. Vol. {\bf 7} (1997), no. 1, pp. 1--47.

\bibitem[Bo]{Bowditch_charhyp}
B.\ Bowditch, 
{\em A topological characterisation of hyperbolic groups}, 
J. Amer. Math. Soc. Vol. {\bf 11} (1998), no. 3, p. 643--667. 

\bibitem[Ch]{Champetier}
C.\ Champetier,
{\em Petite simplification dans les groupes hyperboliques},
Ann. Fac. Sci. Toulouse Math. (6) Vol. {\bf 3} (1994), no. 2, p. 161--221.

\bibitem[CP]{CP}
M. Coornaert, A. Papadopoulos, 
``Symbolic dynamics and hyperbolic groups'', 
Lecture Notes in Mathematics, Vol. {\bf 1539}, Berlin, 1993.

\bibitem[Eb]{Eberlein}
P. Eberlein, ``Geometry of nonpositively curved manifolds'', University of Chicago Press, 1997. 

\bibitem[Fra]{Frances} 
Ch. Frances, 
{\em Lorentzian Kleinian groups}, 
Comment. Math. Helv. Vol. {\bf 80} (2005), no. 4, p. 883--910.

\bibitem[Fre]{Freden}
E.\ Freden,
{\em Negatively curved groups have the convergence property I},
Ann. Acad. Sci. Fenn. Ser. A I Math. Vol. {\bf 20} (1995), no. 2, p. 333--348. 

\bibitem[GGKW]{GGKW}
F. Gu\'eritaud, O.\ Guichard, F.\ Kassel, A.\ Wienhard, 
{\em Anosov representations and proper actions}, 
Preprint arXiv:1502.03811, February 2015. 

\bibitem[Gr]{Gromov_hypgps}
M. Gromov,
{\em Hyperbolic groups}, 
in: ``Essays in group theory", Math. Sci. Res. Inst. Publ. 8, Springer, New York (1987), p. 75--263.

\bibitem[GW]{GW}
O.\ Guichard, A.\ Wienhard, 
{\em Anosov representations: Domains of discontinuity and applications}, 
Invent. Math. Vol. {\bf 190} (2012) no. 2, p. 357--438. 

\bibitem[Hel]{Helgason}
S.\ Helgason,
{\em Differential geometry and symmetric spaces},
Academic Press, 1962.

\bibitem[J]{Jantzen}
J. C. Jantzen, ``Representations of Algebraic Groups'', AMS Mathematical Surveys and Monographs, 2nd Edition, 2003. 

\bibitem[KLM]{ccm} M. Kapovich, B. Leeb and J. J. Millson, {\em Convex functions on symmetric spaces, side lengths of polygons and the stability inequalities for weighted configurations at infinity},  
Journal of Differential Geometry,  Vol. {\bf 81}, 2009, p. 297-- 354. 

\bibitem[KLP1a]{coco13}
M.\ Kapovich, B.\ Leeb, J.\ Porti,
{\em Dynamics at infinity of regular discrete subgroups of isometries of higher rank symmetric spaces},
arXiv e-print, June 2013.

\bibitem[KLP1b]{coco15}
M.\ Kapovich, B.\ Leeb, J.\ Porti,
{\em Dynamics on flag manifolds:
domains of proper discontinuity and cocompactness},
arXiv e-print, October 2015.

\bibitem[KLP2]{morse}
M.\ Kapovich, B.\ Leeb, J.\ Porti, 
{\em Morse actions of discrete groups on symmetric spaces}, 
arXiv e-print,  March 2014. 

\bibitem[KLP3]{mlem}
M.\ Kapovich, B.\ Leeb, J.\ Porti, 
{\em A Morse Lemma for quasigeodesics in symmetric spaces and euclidean buildings}, 
arXiv e-print, November 2014. 

\bibitem[KLP4]{anosov}
M.\ Kapovich, B.\ Leeb, J.\ Porti, 
{\em Some recent results on Anosov representations},
Transformation Groups Vol. {\bf 21} (2016), no. 4, p. 1105--1121.

\bibitem[KL1]{bordif}
M. Kapovich, B. Leeb, 
{\em Finsler bordifications of symmetric and certain locally symmetric spaces}, 
arXiv e-print, 
May 2015.

\bibitem[KL2]{manicures}
M. Kapovich, B. Leeb, 
{\em Discrete isometry groups of symmetric spaces}, 
MSRI Lecture Notes, 
Preprint, January 2017.



\bibitem[KlL1]{qirigid}
B.\ Kleiner, B.\ Leeb, 
{\em Rigidity of quasi-isometries for symmetric spaces 
and euclidean buildings}, 
Inst.\ Hautes \'Etudes Sci. Publ.\ Math.\ No.\ 86 (1997) p. 115--197.

\bibitem[KlL2]{convcoco}
B.\ Kleiner, B.\ Leeb, 
{\em Rigidity of invariant convex sets in symmetric spaces}, 
Invent. Math. Vol. {\bf 163}, No. 3, (2006) p. 657--676. 

\bibitem[La]{Labourie}	
F. Labourie, {\em Anosov flows, surface groups and curves in projective space}, 
Invent. Math. Vol. {\bf 165} (2006), no. 1, p. 51--114.


\bibitem[Le]{habil}
B.\ Leeb, 
{\em A characterization of irreducible symmetric spaces 
and euclidean buildings of higher rank by their asymptotic geometry}, 
Bonner Mathematische Schriften, Vol. {\bf 326} (2000), 
see also  arXiv:0903.0584 (2009).

\bibitem[Mi]{Mineyev}
I.\ Mineyev, 
{\em Flows and joins of metric spaces}, 
Geom. Topol. 9 (2005), p. 403--482. 

\bibitem[P]{Parreau}
A.\ Parreau, 
{\em La distance vectorielle dans les immeubles affines et les espaces sym\'et\-ri\-ques},
in preparation. 

\bibitem[Su]{Sullivan}
D.\ Sullivan, 
{\em Quasiconformal homeomorphisms and dynamics. II. Structural stability
implies hyperbolicity for Kleinian groups}, 
Acta Math. Vol. {\bf 155} (1985), no. 3-4, p. 243--260.

\bibitem[Tu]{Tukia_convgps}
P. Tukia, 
{\em Convergence groups and Gromov's metric hyperbolic spaces}, 
New Zealand J. Math. Vol. {\bf 23} (1994), no. 2, p. 157--187.

\end{thebibliography}
\end{document}